\documentclass[hidelinks,onefignum,onetabnum]{siamart251216} 

\usepackage{lipsum}
\usepackage{amsfonts}
\usepackage{graphicx}
\usepackage{epstopdf}
\usepackage{algorithmic}
\ifpdf
  \DeclareGraphicsExtensions{.eps,.pdf,.png,.jpg}
\else
  \DeclareGraphicsExtensions{.eps}
\fi

\usepackage{benstyle_siam}
\usepackage{cite}
\def\bn{\bm{n}}
\def\br{\bm{\omega}}
\def\bx{\bm{x}}
\usepackage{enumitem}
\usepackage{subcaption}

\usepackage{enumitem}
\setlist[enumerate]{leftmargin=.5in}
\setlist[itemize]{leftmargin=.5in}

\newsiamremark{remark}{Remark}
\newsiamremark{hypothesis}{Hypothesis}
\crefname{hypothesis}{Hypothesis}{Hypotheses}
\newsiamthm{claim}{Claim}
\newsiamremark{fact}{Fact}
\crefname{fact}{Fact}{Facts}

\headers{Stable extrapolation in learning}{B. Adcock, S. Brugiapaglia and X. Wang}

\title{Into the danger zone: stable extrapolation in high-dimensional function and operator learning\thanks{Submitted to the editors on September 30, 2026.
\funding{BA acknowledges support from the Natural Sciences and Engineering Research
Council of Canada (NSERC) through grants RGPIN-2026-04531 and DGDND-2026-04531. SB acknowledges support from NSERC through grant RGPIN-2020-06766. 
BA, SB \& XW acknowledge the support
of FRQ (Fonds de recherche du Qu\'ebec) – Nature et Technologies through grant 359708.}}}

\author{Ben Adcock\thanks{Simon Fraser University
  (\email{ben\_adcock@sfu.ca},\email{xuemeng\_wang\_2@sfu.ca}).}
\and Simone Brugiapaglia\thanks{Concordia University 
  (\email{simone.brugiapaglia@concordia.ca}).}
\and Xuemeng Wang\footnotemark[2]}

\usepackage{amsopn}

\ifpdf
\hypersetup{
  pdftitle={Into the danger zone: stable extrapolation in high-dimensional function and operator learning},
  pdfauthor={B. Adcock, S. Brugiapaglia and X. Wang}
}
\fi

\begin{document}

\maketitle

\begin{abstract}
Out-of-distribution (OOD) generalization is a central challenge in scientific machine learning. We study regression problems in which the test distribution differs from the training distribution and ask: under what assumptions on the target function or operator is stable extrapolation possible, and how far beyond the training domain can one extrapolate? Existing theory controls the test error through additive penalties measuring the discrepancy between the training and test distributions. Such guarantees show robustness to small distribution shifts, but can be very pessimistic in comparison to OOD performance observed empirically. We identify classes of holomorphic functions and operators for which the OOD generalization error converges at algebraic rates even in the presence of large distribution shifts. This phenomenon stems from the increasing smoothness of higher-index coordinates, leading to what we term a `blessing of high dimensionality'. For learning with either polynomials, deep neural networks or deep neural operators, we derive explicit rates for arbitrary test measures supported on suitable domains and quantify how the admissible domain depends on the underlying regularity of the function or operator. Our extrapolation guarantees are independent of the test distribution, depending only on its support. We also present a series of numerical experiments across a range of functions and operators that support the main theoretical findings.
\end{abstract}

\begin{keywords}
out-of-distribution generalization, stable extrapolation, high-dimensional approximation, operator learning, holomorphic functions and operators
\end{keywords}

\begin{MSCcodes}
41A63, 41A10, 46E50, 68Q32, 68T07
\end{MSCcodes}

\section{Introduction}

The ability of models to extrapolate beyond their training domain is a central challenge in modern Machine Learning (ML) and, in particular, Scientific ML (SciML). Understanding and enhancing the \textit{out-of-distribution (OOD)} performance of ML models -- namely, how well a model generalizes to examples drawn from a test distribution that differs from the training distribution -- is a substantial area of research. Yet, as we discuss in \S \ref{ss:OOD-ML}, existing theory of OOD generalization in ML largely addresses the case of small \textit{distributional shifts}, where the training and test distributions are close in some metric. Such guarantees are very general, placing minimal assumptions on the model or target object, but say little about larger shifts. However, a growing body of work in the SciML literature (described further in \S \ref{ss:related}) has shown that trained Deep Neural Network (DNN) and Deep Neural Operator (DNO) models can often \textit{extrapolate}: they generalize well much further from their training distribution than such pessimistic guarantees imply.

This raises a key question: \textit{under what assumptions on the target function or operator is stable extrapolation possible and how far can one extrapolate?} We address this question for classes of holomorphic functions and operators using polynomials, DNNs and DNOs. Our results establish explicit algebraic convergence rates for the OOD generalization error in the $L^\infty$-norm and quantify how the admissible extrapolation domain depends on the underlying regularity 
Our results are \textit{distribution-agnostic}, as they depend on the support of the test distribution only. 
They also reveal a certain \textit{blessing of high dimensionality}: increasing smoothness in higher-index coordinates permits extrapolation in increasingly large domains and, in appropriate regimes, preserves the underlying algebraic convergence order.
Our results are complemented with numerical experiments illustrating extrapolation in both function and operator learning and supporting our main conclusions.

This work combines classical numerical analysis of polynomial extrapolation, modern approximation theory in high and infinite dimensions and the recent advances in DNN/DNO approximation and generalization theory. 
We provide a coherent and comprehensive answer to this question of OOD generalization in the setting holomorphic regularity -- a common assumption in practice, especially in SciML -- thereby bridging a key gap between existing distributional shift theory in ML and practical performance of SciML models.

\subsection{OOD generalization in ML}\label{ss:OOD-ML}

Consider a standard regression problem in ML where a model $\hat{f}$ is trained to approximate an unknown $f : \bbR^d \rightarrow \bbR$. Let $\varrho$ and $\mu$ be two probability measures on $\bbR^d$, the \textit{in-distribution} and \textit{out-of-distribution} measures, respectively. Then typical OOD guarantees in ML (see \S \ref{ss:related} for references) take the form
\be{
\label{standard-ML-OOD-bound}
\nm{f - \hat{f}}^p_{L^p_{\mu}(\bbR^d)} \lesssim \nm{f - \hat{f}}^p_{L^p_{\varrho}(\bbR^d)}+ D(\mu , \varrho),
}
assuming some mild regularity of $f$, $\hat{f}$, such as Lipschitz continuity. Here and elsewhere, we write $A \lesssim B$ to mean there exists a numerical constant $c > 0$ such that $A \leq c B$, and likewise for $A \gtrsim B$. Such an estimate bounds the OOD error in terms of the in-distribution error plus a term $D(\mu,\varrho)$ measuring the distance between the two measures. Often, $D = W_p$ is the Wasserstein $p$-distance. For instance, when $p = 2$, a typical bound is
\be{
\label{Wp-bound}
\nm{f - \hat{f}}^2_{L^2_{\mu}(\bbR^d)} \leq \nm{f - \hat{f}}^2_{L^2_{\varrho}(\bbR^d)} + C(f,\hat{f},\varrho,\mu) W_2(\mu , \varrho).
}
Here $C(f,\hat{f},\varrho,\mu) = (L_{f} + L_{\hat{f}}) \sqrt{4(L_f + L_{\hat{f}})^2 (m_2(\mu) + m_2(\varrho)) + 16 (| f(0) |^2 + |\hat{f}(0)|^2)}$, $L_{\cdot}$ denotes the Lipschitz constant and $m_2(\cdot)$ denotes the second moment.
Unfortunately, these error bounds are only meaningful for small distribution shifts. For example, if $\varrho = \cU([0,1]^d)$ and $\mu = \cU([0,1+\delta]^d)$, then $W_p(\varrho,\mu) \asymp \delta \sqrt{d}$. Hence these bounds become meaningless unless $\delta$ is small. More generally, convergence of the in-distribution error to zero does not imply the same for the OOD error in the presence of a fixed discrepancy $D(\mu,\varrho)$.

\begin{figure}[t]
    \centering
    \begin{minipage}[c]{0.48\linewidth}
        \centering
        \includegraphics[width=\linewidth, trim = 0 0 150 20,clip]{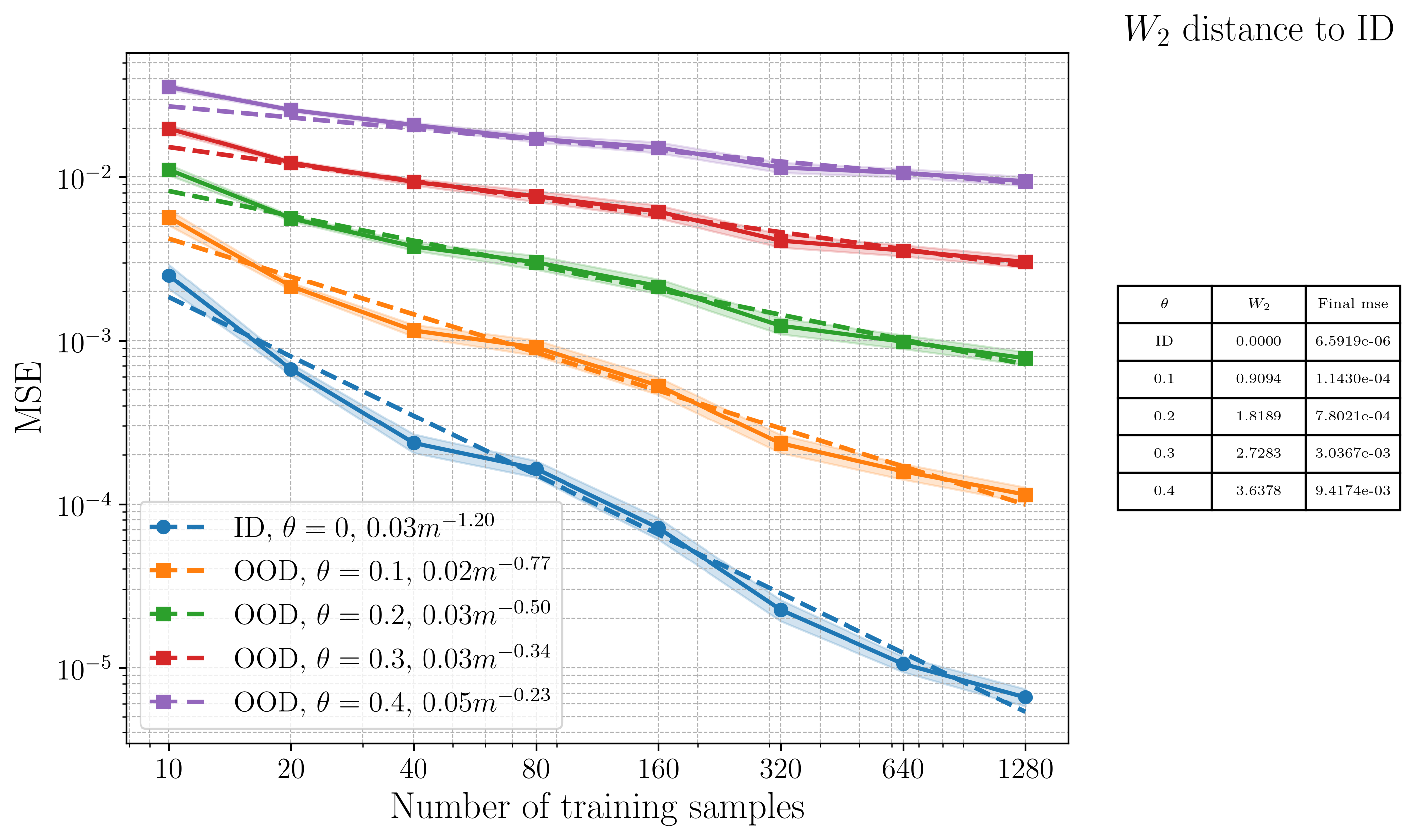}
        \\
        {\small
       \begin{tabular}{c|cccc}
            \hline
            $\theta$ & 0.1 & 0.2 & 0.3 & 0.4 \\
	   $W_2(\varrho,\mu)$ & 0.91 & 1.82 & 2.73 & 3.64 \\
            \hline
        \end{tabular} 
        }
    \end{minipage}
   \begin{minipage}[c]{0.48\linewidth}
        \centering
        \includegraphics[width=\linewidth]
        {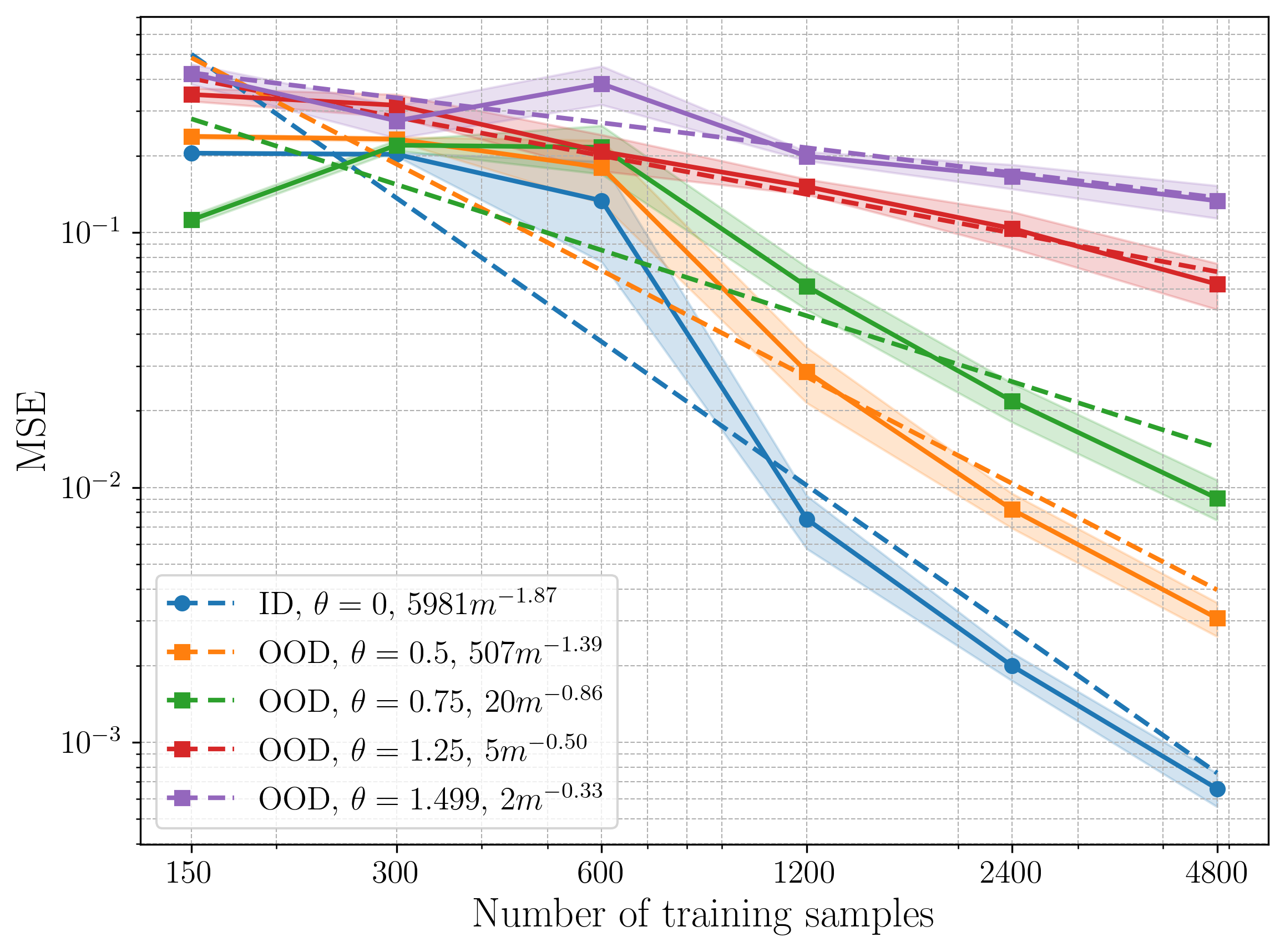}
        \\
        {\small
       \begin{tabular}{c|ccccc}
            \hline
            $\theta$ & 0.5 & 0.75 & 1.25 & 1.499\\ 
	   $W_2(\varrho,\mu)$ & 0.10 & 0.19 & 0.52 & 0.89 \\
            \hline
        \end{tabular} 
        }
    \end{minipage}
    \caption{For both function (left) and operator  (right) learning tasks, the OOD generalization error decays algebraically even though $W_{2}(\varrho,\mu)$ is large. Top row: the squared relative $L^2_{\mu}$-norm error versus number of training samples $m$ for various test distributions $\mu = \mu_{\theta}$. Here $\theta$ is a parameter that controls how far the test distribution is from the training distribution, which corresponds to the case $\theta = 0$. Bottom row: the distance $W_{2}(\varrho,\mu)$ between the training and distributions $\varrho$ and $\mu$. This experiment considers the \textit{wing weight} function (left), learned with a fully-connected DNN, and the solution operator of the Navier--Stokes PDE (right), learned with a so-called Fourier neural operator. See \S \ref{ss:ex-DNNs} and \S \ref{ss:ex-DNO}, respectively, for further details and discussion on these experiments.}
\label{fig:intro-examples}
\end{figure}

However, in practice, ML models often generalize much better than such bounds suggest. To illustrate, in Fig.\ \ref{fig:intro-examples} we show in-distribution and OOD generalization errors for function approximation and operator learning tasks. In the both cases, the terms $W_2(\varrho,\mu)$ are large. Yet, the OOD error decays as $m \rightarrow \infty$, albeit at a slower rate than the in-distribution error (yet still algebraic). Our aim in this paper is to explain this phenomenon.

\subsection{Main results}\label{ss:main-res-intro}

Our results identify classes of functions and operators for which the OOD error converges algebraically for \textit{arbitrary} measures supported in certain explicitly quantified domains. We now describe our setup for functions, before considering operators later.

\textbf{In- and out-of-distribution measures.} 
For reasons we discuss in \S \ref{s:holo-funs}, we consider functions $f : \bbR^{\bbN} \rightarrow \bbR$ of infinitely-many variables. However, our results apply seamlessly to functions of $d < \infty$ variables (Remark \ref{rem:fin-dim}). For the training distribution, we consider the uniform probability measure $\varrho = \cU([-1,1]^{\bbN})$.
Now let $\bm{\omega} = (\omega_i)_{i \in \bbN} \geq \bm{1}$  (understood componentwise). For the test distribution, we consider any probability measure $\mu$ with $\mathrm{supp}(\mu) \subseteq D_{\bm{w}}$, where
 \be{
\label{Dr-def}
D_{\bm{\omega}} = [-\omega_1 , \omega_1] \times [-\omega_2,\omega_2] \times \cdots .
}

\textbf{Training data and learning problem.} We consider the standard regression setting in ML, where $\bm{x}_1,\ldots,\bm{x}_m \sim_{\mathrm{i.i.d.}} \varrho$ and we are given noisy training data $(\bm{x}_i , y_i : = f(\bm{x}_i) + e_i)$, $i = 1,\ldots,m$. 
Our results consider a bounded noise model, which may be adversarial, whose contribution is controlled by $\nm{\bm{e}}_2$. Now consider a model class $\cM$  -- i.e., a set of functions $\bbR^{\bbN} \rightarrow \bbR$, taken in this work to be family of polynomials or DNNs -- and the estimator 
\be{
\label{l2-loss-intro}
\hat{f} \in \argmin{g \in \cM} \frac1m \sum^{m}_{i=1} |g(\bm{x}_i) - y_i |^2.
}

\textbf{Class of functions.} Given $\varepsilon > 0$ and $\bm{b} \in \ell^1(\bbN)$ with $\bm{b} \geq \bm{0}$, we consider the class $\cH(\bm{b},\varepsilon)$ of  $(\bm{b},\varepsilon)$-holomorphic functions (with norm at most one). As we discuss further in \S \ref{ss:contributions}-\ref{ss:related}, this is an important class in high-dimensional function and operator learning.

\textbf{Extrapolation condition.} As noted, our results employ a multiplicative OOD error bound which allows for large distribution shifts, with size determined by regularity of the functions being learned. To state this condition, we require several further concepts. First, we define the \textit{monotone} $\ell^p$-space, $0 < p \leq \infty$, denoted $\ell^p_{\mathsf{M}}(\bbN)$, as the space of all sequences $\bm{z} = (z_i)_{i \in \bbN}$ whose minimal monotone majorant $\tilde{\bm{z}} =  (\tilde{z}_i)_{i \in \bbN} \in \ell^p(\bbN)$, where $\tilde{z}_i : = \sup_{j \geq i} | z_j |$. Second, we define $\zeta : [1,\infty) \rightarrow [1,\infty)$ as $\zeta(\omega) = \omega + \sqrt{\omega^2-1}$. We now introduce the following condition relating $\bm{b}$, $\bm{\omega}$ and a further scalar $p \in (0,1)$: 
\be{
\label{b-r-cond-intro}
\bm{b} \odot \bm{r}_p \in \ell^p_{\mathsf{M}}(\bbN),\quad \nm{\bm{b} \odot (\bm{r}_p-\bm{1}) }_1 < \varepsilon,\qquad \text{where } \bm{r}_p =  (\zeta(\omega_i)^{2/p-1})_{i \in \bbN}.
}
Here $\odot$ is the Hadamard product. The first condition is a decay requirement on the sequence $\bm{b}$. The second limits the enlargement of the domain using the  holomorphy budget $\varepsilon$. As we see below, smaller values of $p$ give faster convergence rates, but increase $\bm{r}_p$ and therefore strengthen both requirements. We discuss concrete examples of this condition in \S \ref{ss:condition-examples}.

\thm{[Theorem \ref{thm:poly-optimized-least-squares}, informal version]
Let $\bm{\omega} \geq \bm{1}$ and $D_{\bm{\omega}}$ be as in \ef{Dr-def}. Then there exists a polynomial model class $\cM = \cP$ with the following property. Let $f \in \cH(\bm{b},\varepsilon)$, where $\bm{b}$ is such that \ef{b-r-cond-intro} holds for some $0 < p < 1$. Then any minimizer $\hat{f}$ of \ef{l2-loss-intro} satisfies
\bes{
\nm{f - \hat{f}}_{L^{\infty}_{\mu}(\bbR^{\bbN})} \lesssim C(\bm{b},\varepsilon,\bm{\omega},p)  ( m / \log^4(m)  )^{1-1/p} + \nm{\bm{e}}_2 / \log^2(m)
}
with high probability, for any measure $\mu$ supported in $D_{\bm{\omega}}$. Moreover, all elements of $\cP$ are polynomials with at most $m/\log^4(m)$ terms in the Legendre basis.
}

\thm{[Theorem \ref{thm:deep-learning-extrap}, informal version]
\label{thm:informal-2}
Consider the setup of the previous theorem, where $\omega_{\min} = \min \{ \omega_i \} > 1$. Then there exists a family of feedforward tanh DNNs $\cM = \cN$ achieving the same bound, up to an additional additive term proportional to $2^{-m}$. Moreover $\mathrm{width}(\cN) \lesssim m/\log(\omega_{\min})$ and $\mathrm{depth}(\cN) \lesssim \log ( \log(m)/\log(\omega_{\min}) )$.
}

Finally, we also consider \textit{operator learning} \cite{brugiapaglia2026short,boulle2024mathematical,kovachki2024operator,subedi2026operator}, where the target is an operator $F : \cX \rightarrow \cY$ between two separable Hilbert spaces. We consider an in-distribution probability measure $\nu$ on $\cX$, an OOD measure $\mu$ on $\cX$,
training data $(X_i, Y_i : =  F(X_i) + E_i )$, $i = 1,\ldots,m$, where $X_i \sim_{\mathrm{i.i.d.}} \nu$,
a class of holomorphic operators $\cH(\bm{b},\varepsilon ; \cX , \cY)$ (see Definition \ref{def:holo-op}), a family $\cN$ of \textit{neural operators} $\cX \rightarrow \cY$ and an estimator $\hat{F} \in \argmin{G \in \cN} \frac1m \sum^{m}_{i=1} \nm{G(X) - Y_i }^2$.

\thm{
[Theorem \ref{thm:operator-learning-extrap}; informal version]
\label{thm:main-opl-informal}
Let $\bm{\omega} > \bm{1}$, where $\omega_{\min} = \min \{ \omega_i \} > 1$, and suppose that Assumptions \ref{ass:opl}-\ref{ass:opl-ood} hold. Then there is a class of DNOs $\cN$ 
such that, for any $F \in \cH(\bm{b},\varepsilon;\cX , \cY)$ and $\bm{b}$ satisfying \ef{b-r-cond-intro}, the estimator $\hat{F}$ achieves a similar bound with high probability. Moreover, $\mathrm{width}(\cN) \lesssim m/\log(\omega_{\min})$ and $\mathrm{depth}(\cN) \lesssim \log ( \log(m)/\log(\omega_{\min}) )$.
}

See \S \ref{s:operator-learning} for full details on this case. This theorem asserts a similar OOD bound for holomorphic operators to that given in Theorem \ref{thm:informal-2} for  holomorphic functions.

\subsection{Contributions}\label{ss:contributions}

We now discuss several key features of this work.

\textbf{(i) Practical relevance, unknown anisotropy.} Holomorphic functions and operators arise in many parametric PDE and operator learning problems. See \S \ref{ss:related} for further discussion. Their anisotropic regularity permits algebraic approximation rates even in infinite dimensions. 
A notable feature of our results is that the training problem (i.e., the model class $\cM$ and loss function) are independent of $\bm{b}$. They are consequently applicable to the \textit{unknown anisotropy} setting, where the estimator has no \textit{a priori} knowledge of the smoothness of the target function or operator. As discussed in \cite[\S 3.3]{adcock2024learning}, this is particularly relevant in applications, where the target may be a black-box. 

\textbf{(ii) Distribution-agnostic bounds.} Our $L^\infty$-bounds hold uniformly over test distributions $\mu$ supported
in a prescribed extrapolation domain $D_{\bm{\omega}}$. The model class may depend on $\bm{\omega}$, but it requires no further knowledge of $\mu$. Notably, our bounds require neither a bounded density $\D \mu / \D \varrho$ nor proximity of $\mu$ to $\varrho$ in some metric. In particular, $\mu$ may be singular with respect to $\varrho$ and place all its mass outside the training support. This is a major conceptual departure of our work from standard theory (recall \S \ref{ss:OOD-ML}).

\textbf{(iii) Quantifiable extrapolation.} \ef{b-r-cond-intro} quantifies how far one can extrapolate
while retaining algebraic convergence. For example, suppose that $\omega_i = \omega$, $\forall i \in \bbN$. Then \ef{b-r-cond-intro} holds whenever $\bm{b} \in \ell^p_{\mathsf{M}}(\bbN)$ and $(\zeta(\omega)^{2/p-1}-1)\nm{\bm{b}}_1<\varepsilon$ and yields the noise-free convergence rate $(m/\log^4(m))^{1-1/p}$. Thus, one can extrapolate to arbitrary accuracy as $m \rightarrow \infty$ to hypercubes whose side lengths are a constant size larger than the in-distribution domain $[-1,1]^{\bbN}$. This is very different to bounds of \S \ref{ss:OOD-ML}, where the discrepancy $D(\mu,\varrho)$ limits the accuracy. 
But \ef{b-r-cond-intro} also allows for more. Assuming sufficient decay of the $b_i$, it permits domains where $\omega_i \rightarrow \infty$ as $i \rightarrow \infty$. For example, suppose that $b_i = c i^{-a}$ and $\omega_i = 1+i^b$ for $a,b,c > 0$ satisfying \eqref{b-r-cond-intro}. For $c>0$ sufficiently small, our results
give rates arbitrarily close to $(m/\log^4(m))^{1-(a+b)/(1+2b)}$. This example illustrates that one can extrapolate arbitrarily far in the $i$th variable as $i \rightarrow \infty$. Such a conclusion is, on the face of it, quite remarkable, and very unintuitive from the perspective of distribution shift theory. See \S \ref{ss:condition-examples} for further discussion. 

\textbf{(iv) Stability and a blessing of high dimensionality.} In low dimensions, polynomial approximations to holomorphic functions achieve in-distribution errors that decay exponentially. Exponential rates can be retained for extrapolation, but at the cost of exponential instability. Such instability can be avoided, but with a dramatic reduction in convergence order from exponential to algebraic. By contrast, natural in-distribution rates in infinite dimensions are algebraic. Our results in (iii) reveal a \textit{blessing of high dimensionality}. Stable extrapolation to admissible hypercubes preserves the same algebraic convergence order as the in-distribution rates. Moreover, algebraic rates, albeit slower, persist for certain hyperrectangles with unbounded side lengths. See Remark \ref{rem:blessing} for further details.
 
\textbf{(v) Sample efficiency and optimal rates.} Problems in SciML are often data-starved, making it critical to establish results on the sample complexity of learning tasks. We consider standard i.i.d.\ samples from the underlying measure $\varrho$. In particular, the samples, being \textit{independent} of $\bm{\omega}$, do not need to be adapted to the extrapolation domain. Our results show that holomorphic functions and operators can be learned in a sample-efficient manner with algebraic rates in terms of $m$, which are also optimal in various settings (see Remark \ref{rem:rate-optimality}).

\textbf{(vi) Generalization to rougher inputs.} A crucial consideration in operator learning is generalization to \textit{rougher} inputs (i.e., functions of lower regularity) than those seen in training \cite{jiang2026what, zhu2023reliable, subedi2026operator, dehoop2023convergence,benitez2024out}. Fig.\ \ref{fig:intro-examples} already shows an example of this effect, as the eigenvalues of covariance operator of the test distributions decay more slowly than those of the training distribution (see \S \ref{ss:ex-DNO} for further examples).
Theorem \ref{thm:main-opl-informal} shows this is possible: the learned model can generalize to test distributions $\mu$ that are \textit{arbitrarily rough}, and do so in a quantifiable manner, given sufficient smoothness of the target operator. See \S \ref{ss:examples-opl} for further discussion.

\subsection{Related work}\label{ss:related}

We now discuss various related areas of work.

\textbf{The challenge of OOD generalization.}
Practical demands mean that SciML models are often applied to inputs drawn from distributions differing significantly from their training distribution. Many recent works have discussed the importance and challenges of OOD generalization \cite{brivio2026ptpi,setinek2026simshift,wang2023scientific,muckley2023interpretable,li2024physics,bonnet2023airfrans,chen2024dataefficient,shikhman2026diagnosing,nguyen2026out,jiang2026what,chu2026do}. More broadly, OOD generalization is a critical issue in AI for Science \cite{zhang2025artificial}, and a key component of the \textit{robustness} of AI models. 

Understanding and enhancing the OOD performance is also a large topic in ML in general, with well-established concepts and techniques such as domain adaptation and generalization, distributionally robust optimization, invariant risk minimization, risk extrapolation  and others.
See, e.g., \cite{xu2021how,ben-david2006analysis,arjovsky2021out,liu2023towards,krueger2021OO,yuan2022towards,nguyen2026out}. Much of the literature focuses on classification, as opposed to regression problems, making it not directly applicable to SciML \cite{yuan2022towards}. Further, distribution shifts in SciML are typically highly structured \cite{shikhman2026diagnosing}, which renders concepts and approaches, including the distribution shift bounds discussed in \S \ref{ss:OOD-ML}, potentially less relevant.

\textbf{Existing OOD works.}
Distributional shift bounds of the form \ef{standard-ML-OOD-bound} can be found in many places in the ML literature (see \cite{hou2023instance,ben-david2010theory,benitez2024out,guerra2025learning} and references therein). Use of the $W_p$-distance is standard \cite{sinha2018certifying,lee2018minimax,blanchet2019quantifying,levine2020wasserstein,kuhn2025distributionally,benitez2024out,guerra2025learning,hou2023instance} and the specific estimate \ef{Wp-bound} is from \cite[Prop.\ 3.1]{guerra2025learning}, which is based on \cite{benitez2024out,hou2023instance}. See \cite{ben-david2010theory} for similar bounds using the TV distance and so-called $\cH$-divergences. Our results contrast with these general-purpose \textit{additive} bounds, as they rely on a multiplicative bound (Lemma \ref{lem:multiplicative-OOD}). They are also distribution-free, as they depend only on the support. A distribution-free bound was also developed in \cite[Cor.\ 5.2]{hong2024bridging} using the modulus of continuity. But it is also additive in nature, thus only relevant to small shifts.

OOD generalization in SciML has been investigated, primarily empirically, in many works \cite{brivio2026ptpi,muckley2023interpretable,li2025probing,shikhman2026diagnosing,yuan2022towards,zhu2023reliable,benitez2024out,setinek2026simshift,fesser2023understanding,chu2026do,shikhman2026diagnosing,bonfanti2024generalization}. Several of these demonstrate successful generalization under
substantial shifts \cite{zhu2023reliable,brivio2026ptpi,dehoop2023convergence,shikhman2026diagnosing, nguyen2026out}, a phenomenon our work strives to explain. 
Among theoretical contributions, \cite{benitez2024out} derives OOD bounds for neural operators applied to high-frequency PDE problems. Yet these employ the bounds of \S \ref{ss:OOD-ML}, and are therefore valid only for small shifts. More closely related to our work, \cite{dehoop2023convergence} establishes algebraic convergence rates for Bayesian learning of linear operators. Our work can be viewed a generalization from linear operators to holomorphic operators, albeit in a non-Bayesian context. The work \cite{xu2021how} also shares some similarity with ours, since it considers genuine extrapolation versus small shifts. But it does not provide generalization error bounds for concrete function classes.

\textbf{Assumptions and technical ingredients.} Infinite-dimensional homolorphy was developed extensively in the study of parametric DEs. It is now known that many families of parametric DEs have solution maps that are $(\bm{b},\varepsilon)$-holomorphic functions. The list includes parametric diffusion equations, parametric heat equations, PDEs over parametrized domains and various parametric initial value problems. See \cite{nobile2008anisotropic,nobile2008sparse,chkifa2015breaking,adcock2024optimal,schwab2019deep, cohen2015approximation,adcock2022sparse} and references therein. Holomorphic operators have also received substantial attention in recent years. See \cite{adcock2024optimalb,herrmann2024neural,reinhardt2024statistical,westermann2026performance,kovachki2024operator} and references therein. Besides their relevance to applications, holomorphic operators represent one of the few known classes of operators that are not subject to the \textit{curse of parametric complexity} or \textit{curse of sample complexity}, in contrast to finite-regularity operators \cite{lanthaler2024parametric,adcock2025sample}.
However, virtually all theory on learning holomorphic functions or operators considers in-distribution performance. Our work extends this to the OOD regime.

Besides holomorphy, our work has two further technical ingredients. First, modern polynomial approximation theory of infinite-dimensional functions. See, e.g., \cite{cohen2011analytic,chkifa2015breaking,adcock2022sparse,adcock2024learning,cohen2015approximation} and references therein. Second, the powerful technique of \textit{emulation} of polynomials via DNNs/DNOs \cite{guhring2021approximation,lu2021deep,yarotsky2017error,schwab2019deep,opschoor2022exponential} (see \cite[\S 7.1]{adcock2024learning} for a review). Note that emulation techniques have been at the core of the modern approximation theory of DNNs/DNOs which has developed in the last decade \cite{elbrachter2021deep,devore2020neural}. While emulation techniques are normally employed to show the DNNs/DNOs with certain approximation guarantees, in this work we follow ideas first developed in \cite{adcock2021deep} to employ them to establish full generalization bounds.

\textbf{Relation to classical analytic continuation.}
Extrapolation of holomorphic functions is synonymous with analytic continuation, which is a classical topic in numerical analysis (see \cite{trefethen2020quantifying} and references therein). The work \cite{demanet2019stable} considered analytic continuation of univariate functions via polynomials and linear least squares. Our results can be seen as a generalization of \cite{demanet2019stable} to multivariate functions.
In \cite{trefethen2020quantifying} it is shown that univariate analytic continuation is ill-posed without further constraints, such as boundedness. The class $\cH(\bm{b},\varepsilon)$ consists of unit norm, and therefore bounded, functions and hence avoids this problem.
But \cite{trefethen2020quantifying} shows that analytic continuation is always ill-conditioned with respect to $\ell^{\infty}$-norm perturbations of the data, with infinite condition number. Our results do not contradict this statement, as they only show $\ell^2$-norm stability via the term $\nm{\bm{e}}_2$. See Remark \ref{rem:condition-number} for further discussion.

\subsection{Limitations}\label{ss:limitations}
Holomorphy is a strong assumption. Although relevant to many SciML problems, other problems lack this regularity. Indeed, there is a growing discussion in operator learning as to what are the right structures beyond regularity that may govern learnability \cite{brugiapaglia2026short,boulle2024mathematical,kovachki2024operator,subedi2026operator}. Our work does not strive to address this question, which is largely open even in the in-distribution context (see also \S \ref{s:conclusion}). On a related note, the fact that our DNN/DNO results use emulation means that the resulting estimators perform no better than corresponding polynomial estimators. Our results also employ handcrafted DNN/DNO model classes, which are some ways away from those typically used in practice (see also \S \ref{s:conclusion}). Nonetheless, our work aligns with a recent trend in operator learning, in which more classical approaches such as polynomials \cite{westermann2026performance} and kernel methods \cite{batlle2024kernel} have been revisited and generalized to the operator setting, and shown to be competitive with DNOs in terms of in-distribution performance. Our theoretical results provide similar conclusions for OOD performance.

There are also many different types of OOD generalization in SciML. In operator learning, for instance, there is much interest in multi-operator learning and foundation models, i.e., DNOs that can simultaneously solve multiple PDEs \cite{choi2025defining,jiang2026what,chen2024dataefficient,yuan2022towards, herde2024poseidon,yang2026generalization,  shikhman2026diagnosing}. This is a type of OOD generalization where one wants to generalize to PDEs not seen in training. While not unrelated, it differs from our setting, as we consider learning a single operator only. Our setting is also different from extrapolation-in-time, as is common in neural PDE solvers \cite{fesser2023understanding,wang2023longtime, jiang2026what}. Another distinct type of extrapolation is superresolution, where one trains a DNO on a coarse grid, then uses finer grids at test time \cite{li2021fourier} (see also \cite{subedi2026zeroshot} and references therein). 

Finally, we stress that the objective of this work is to understand the OOD generalization of estimators, rather than methods of enhancing it. We discuss this important topic further in \S \ref{s:conclusion}.

\subsection{Outline}

In \S \ref{s:holo-funs} we introduce the class of functions considered and discuss their approximation via polynomials. In \S \ref{s:in-to-out} we begin our focus on OOD errors, by deriving a key multiplicative bound which underpins all our main results. Our main results are presented in \S \ref{s:sample-efficient-learning} (for functions) and \S \ref{s:operator-learning} (for operators). In \S \ref{s:experiments} we present numerical experiments and in \S \ref{s:conclusion} we conclude and discuss open problems. Proofs of all results in this paper can be found in \S \ref{s:proofs}, while \S \ref{app:experiments} contains additional details on the numerical experiments.

\section{Holomorphic functions and polynomial approximation}\label{s:holo-funs}

As mentioned in \S \ref{ss:main-res-intro}, we consider functions $f : \bbR^{\bbN} \rightarrow \bbR$ and the \textit{in-distribution} measure $\varrho = \cU([-1,1]^{\bbN})$. However, our theory applies seamlessly to functions of finitely-many variables (Remark \ref{rem:fin-dim}). Working in infinite dimensions is convenient theoretically and practically relevant, even for finite-dimensional functions, as finite-dimensional functions behave like infinite-dimensional functions in terms of their approximation theory as soon as the dimension is moderately large (Remark \ref{rem:fin-dim-funs}).


\subsection{Polynomial expansions and $s$-term approximation}

We first recap some standard notions from polynomial approximation. See, e.g., \cite{adcock2022sparse,cohen2015approximation}. Let $P_n$, $n \in \bbN_0$, denote the Legendre polynomials on $[-1,1]$, with normalization $P_n(1) = 1$. We define the normalized polynomials
\be{
\label{psi-n-def}
\psi_n(x) = \sqrt{2 n + 1} P_n(x),\quad \forall n \in \bbN_0
}
and note that these form an orthonormal basis with respect to the uniform probability measure on $[-1,1]$. To extend to $\bbR^{\bbN}$, we proceed by tensorization. Let $\bn = (n_1,n_2,\ldots) \in \bbN_0^{\bbN}$ and write $\mathrm{supp}(\bn) = \{ i : n_i \neq 0 \} \subseteq \bbN$ for its support. We define $\cF =  \{ \bn \in \bbN^{\infty}_0 : | \supp(\bn) | < \infty \}$
and consider the multivariate polynomials $\Psi_{\bn}(\bx) = \prod_{i \in \bbN} \psi_{n_i}(x_i)$.
The set $\{ \Psi_{\bn} \}_{\bn \in \cF}$ forms an orthonormal basis of $L^2_{\varrho}(\bbR^{\bbN})$ and therefore any $f \in L^2_{\varrho}(\bbR^{\bbN})$ has the convergent expansion
\be{
\label{f-exp}
f = \sum_{\bn \in \cF} c_{\bn} \Psi_{\bn},\qquad \text{where }
c_{\bn} = \ip{f}{\Psi_{\bn}}_{L^2_{\varrho}(\bbR^{\bbN})} = \int_{\bbR^{\bbN}} f(\bx) \overline{\Psi_{\bn}(\bx)} \D \varrho(\bx).
}
We are interested in $n$-term approximations to such functions. An \textit{$n$-term approximation} to $f$ based on an multi-index set $S \subset \cF$, $|S| = n$, has the form
\be{
\label{f-S-def}
f_{S} = \sum_{\bn \in S} c_{\bn} \Psi_{\bn} \in \cP_{S},\qquad \text{where } \cP_S = \spn \{ \Psi_{\bm{n}} : \bm{n} \in S \}.
}

\subsection{Infinite-dimensional holomorphic functions}

We now introduce the class of holomorphic functions considered. In the univariate setting, it is well known that the expansion of a function $f$ using the first $n$ Legendre polynomials converges with rate $\rho^{-n}$, where $\rho \geq 1$ is the parameter of the largest \textit{Bernstein ellipse} $\cE_{\rho}$ within which $f$ is holomorphic. Here
\bes{
\cE_{1} = [-1,1],\qquad \cE_{\rho} = \{ (z + z^{-1})/2 : z \in \bbC,\ 1 \leq |z| \leq \rho  \} \subset \bbC,\quad \rho > 1.
}
Generalizing via tensor products, we define the \textit{Bernstein polyellipse} of parameter $\bm{\rho} = (\rho_1,\rho_2,\ldots) \in [1,\infty)^{\bbN}$ as $\cE_{\bm{\rho}} = \cE_{\rho_1} \times \cE_{\rho_2} \times \cdots \subset \bbC^{\bbN}$.

\defn{[$(\bm{b},\varepsilon)$-holomorphy]
\label{d:holomorphy}
Let $\bm{b} \in [0,\infty)^{\bbN}$ and $\varepsilon > 0$. A function $f : \bbR^{\bbN} \rightarrow \bbC$ is \textit{$(\bm{b},\varepsilon)$-holomorphic} if it is holomorphic in every Bernstein polyellipse $\cE_{\bm{\rho}}$ with parameter $\bm{\rho} = (\rho_i)^{\infty}_{i=1} \in [1,\infty)^{\bbN}$ satisfying
\be{
\label{b-eps-holo}
\sum^{\infty}_{i=1} \left ( (\rho_i + \rho^{-1}_i) / 2 - 1 \right ) b_i \leq \varepsilon.
}
}
See \S \ref{ss:related} for relevant literature on this class.
For convenience we write
\bes{
\cR_{\bm{b},\varepsilon} = \bigcup \{ \cE_{\bm{\rho}} : \bm{\rho} \in [1,\infty)^{\bbN},\ \text{$\bm{\rho}$ satisfies \eqref{b-eps-holo}} \} \subseteq \bbC^{\bbN}
}
for the region and $\cH(\bm{b},\varepsilon) =  \{ \text{$f : D \rightarrow \bbC$ $(\bm{b},\varepsilon)$-holomorphic},\ \nm{f}_{L^{\infty}(\cR_{\bm{b},\epsilon})} \leq 1 \}$
for the set of functions that are holomorphic in $\cR_{\bm{b},\varepsilon}$ with uniform norm at most one. 

Note that the sequence $\bm{b}$ determines the type of anisotropic behaviour of functions in $\cH(\bm{b},\varepsilon)$. Indeed, if $b_j$ is large for some $j \in \bbN$, then \eqref{b-eps-holo} holds only for small values of $\rho_j$, meaning that $f$ is less smooth with respect to the variable $x_j$. Conversely, if $b_j$ is small (or even $b_j = 0$), then $f$ is more smooth (entire) in the variable $x_j$.

\rem{[Unknown anisotropy]
\label{rem:unknown-anisotropy}
In some cases, suitable parameters $(\bm{b}, \varepsilon)$ may be known, e.g., via \textit{a priori} analysis of a given problem. However, in practical cases, values for these parameters are typically unknown. We refer to this as the \textit{unknown anisotropy} setting. Our main focus is on this more practical case, with our goal being to design estimators with guaranteed convergence rates that are independent of $(\bm{b},\varepsilon)$. See \cite[\S 3.3]{adcock2024learning} for more discussion on known versus unknown anisotropy.
}

\rem{[Finite-dimensional functions]
\label{rem:fin-dim}
Definition \ref{d:holomorphy} is somewhat complicated, in that it requires the function to have a holomorphic extension to a union of Bernstein polyellipses. This is used in infinite dimensions to obtain algebraic rates of convergence of the best $n$-term approximation. In finite dimensions, it is enough for the function to be holomorphic in a single Bernstein polyellipse. Concretely, let $f : \bbR^d \rightarrow \bbR$ be holomorphic in the finite-dimensional Bernstein polyellipse $\cE_{\bar\rho_1} \times \cdots \times \cE_{\bar\rho_d} \subset \bbC^d$. Now let 
\eas{
b_i = \varepsilon \left ( (\bar\rho_i + \bar\rho^{-1}_{i})/2 - 1 \right )^{-1},\ i \in [d],\qquad b_i = 0,\ i \in \bbN \backslash [d].
}
Then the natural extension $\tilde{f}$ of $f$ to a function of infinitely-many variables, defined by $\tilde{f}(\bm{x}) = f(x_1,\ldots,x_d)$ for all $\bm{x} = (x_i)_{i \in \bbN} \in \bbR^{\bbN}$, is $(\bm{b},\varepsilon)$-holomorphic. Hence, all results that follow also apply seamlessly to finite-dimensional holomorphic functions.
}

\subsection{Best $n$-term approximation of holomorphic functions}

Best $n$-term approximation is a type of nonlinear approximation \cite{devore1998nonlinear} in which one chooses the subset $S \subset \cF$, $|S| = n$, that minimizes the error $f - f_S$ in a certain norm, where $f_S$ is the corresponding approximation \ef{f-S-def}. We shall consider this with respect to the $L^2_{\varrho}$-norm (although other $L^p_\varrho$-norms, $p \neq 2$,
could equally be considered). We define an $L^2_{\varrho}$-norm best $n$-term approximation as
\bes{
f_n = f_{S^*},\qquad \text{where }S^* \in \argmin{} \{ \nm{f - f_S}_{L^2_{\varrho}(\bbR^{\bbN})} : S \subset {\mathcal{F}},\ |S| = n \}.
}
Due to Parseval's identity, $S^*$ has an explicit characterization as the index set corresponding to $n$ largest coefficients of $f$ in absolute value (note that $S^*$, and therefore $f_n$, need not be unique, due to ties, but this makes no difference in what follows). The following result is well known (see, e.g., \cite[Thm.\ 3.28]{adcock2022sparse} or \cite[\S 3.2]{cohen2015approximation}). It demonstrates that the best $n$-term approximation of any function in $\cH(\bm{b},\varepsilon)$ converges with algebraic rate whenever $\bm{b}$ is $\ell^p$-summable.

\thm{
[Algebraic convergence of the best $s$-term approximation]
\label{t:best-s-term}
Let $\varepsilon > 0$ and {$\bm{b} \in [0,\infty)^{\bbN}$ be such that $\bm{b} \in \ell^p(\bbN)$ for some $0 < p < 1$}. Then, for any $2 \leq q \leq \infty$,
\be{
\label{n-term-poly-holo-err}
\nm{f - f_n}_{L^{q}_{\varrho}(\bbR^{\bbN})} \leq C(\bm{b},\varepsilon,p) \cdot n^{1-\frac1q-\frac1p},\qquad \forall f \in \cH(\bm{b},\varepsilon),\ n \in \bbN.
}
}

This remarkable result states that there are classes of functions whose approximation by polynomials is free from the
\textit{curse of dimensionality}, despite depending on infinitely-many variables. Recall that $\bm{b}$ controls the anisotropy of functions in $\cH(\bm{b},\varepsilon)$. The faster $\bm{b}$ decays (i.e., the smaller $p$), the more anisotropic the functions and the faster the rate in \ef{n-term-poly-holo-err}.

\rem{
[Finite-dimensional functions and exponential rates]
\label{rem:fin-dim-funs}
Standard polynomial approximation theory asserts that the $n$-term approximation of a function of $d < \infty$ variables that is holomorphic in $\cE_{\rho_1} \times \cdots \times \cE_{\rho_d}$ converges exponentially fast in $n^{1/d}$. While exponential, this rate exhibits the curse of dimensionality, and raises the question: what happens for large $d$? The bound \ef{n-term-poly-holo-err} provides an answer. By assuming holomorphy \textit{and} anisotropy, one achieves $d$-independent algebraic rates (here we recall that \ef{n-term-poly-holo-err} also holds for any holomorphic function of $d$ variables, due to Remark \ref{rem:fin-dim}). Typically, the error follows the algebraic rate for small $n$, before transitioning to the exponential rate as $n$ grows larger. However, when $d \approx 10$ or larger, this transition point typically occurs at a large value of $n$, meaning that the exponential rates are often not witnessed in practical computations \cite{adcock2024learning}. Thus, algebraic rates are the relevant even in finite dimensions, except when $d$ is small. See \S \ref{ss:num-exp-poly} for several examples of this phenomenon.
}

\section{From in-distribution to out-of-distribution}\label{s:in-to-out}

The bound \ef{n-term-poly-holo-err} considers only the in-distribution error in the $L^q_{\varrho}$-norm. Further, it only asserts the existence of a polynomial achieving the specified algebraic rate, but no insight into how to compute from data. In this section, we focus on the OOD error, for an arbitrary measure $\mu$. We introduce an OOD generalization guarantee which, in contrast to those discussed in \S \ref{ss:OOD-ML}, is multiplicative rather than additive. We then apply this to the polynomial setting in Theorem \ref{t:OOD-s-term}, demonstrating that it guarantees good OOD generalization up to the size of the in-distribution error.

\subsection{A multiplicative OOD generalization guarantee}

Let $\cM$ be an arbitrary \textit{model class} or \textit{hypothesis set} for the regression problem, i.e., a set of functions $\bbR^{\bbN} \rightarrow \bbR$ from which we aim to compute an estimator $\hat{f} \in \cM$ to an unknown function $f$. Given two measures $\varrho,\mu$ on $\bbR^{\bbN}$, we define the \textit{distribution shift} constant of $\cM$ as  
\be{
\label{dist-shift-const}
\Delta(\cM ; \varrho,\mu) = \sup \left \{ \frac{\nm{g_1 - g_2}_{L^{\infty}_{\mu}(\bbR^{\bbN})}}{\nm{g_1 - g_2 }_{L^{2}_{\varrho}(\bbR^{\bbN})}} : g_1,g_2 \in \cM,\ g_1 \neq g_2 \right \}.
}
This constant measures the worst-case OOD behaviour of the difference of two elements of $\cM$ relative to their in-distribution behaviour. Note that it is possible to formulate $\Delta(\cM ; \varrho , \mu)$ in terms of other norms. We use the $L^2_{\varrho}$-norm as we generally consider estimators that are defined through empirical least-squares fits. We use the $L^{\infty}_{\mu}$-norm in order to have error bounds that depend on the support of the measure $\mu$ only (since, in our work, all functions are continuous).

\lem{
[Multiplicative OOD generalization bound]
\label{lem:multiplicative-OOD}
Let $\cM$ and $\Delta(\cM ; \varrho,\mu)$ be as above and $f \in L^{\infty}_{\mu}(\bbR^{\bbN})$. Then, for any $\hat{f} \in \cM$,
\be{
\label{OOD-ID-bound}
\begin{split}
\nm{f - \hat{f}}_{L^{\infty}_{\mu}(\bbR^{\bbN})} \leq& ~ \inf_{g \in \cM} \left \{ \Delta(\cM ; \varrho,\mu) \nm{f - g}_{L^{2}_{\varrho}(\bbR^{\bbN})} +  \nm{f - g}_{L^{\infty}_{\mu}(\bbR^{\bbN})} \right \} 
\\
& + \Delta(\cM ; \varrho,\mu) \nm{f - \hat{f}}_{L^{2}_{\varrho}(\bbR^{\bbN})} .
\end{split}
}
}

\prf{
By the triangle inequality and the definition of $\Delta(\cM ; \varrho,\mu)$, we have
\eas{
\nm{f - \hat{f}}_{L^{\infty}_{\mu}(\bbR^{\bbN})} & \leq \nm{f - g}_{L^{\infty}_{\mu}(\bbR^{\bbN})} + \nm{g - \hat{f}}_{L^{\infty}_{\mu}(\bbR^{\bbN})}
\\
& \leq \nm{f - g}_{L^{\infty}_{\mu}(\bbR^{\bbN})}+ \Delta(\cM ; \varrho,\mu) \nm{g - \hat{f}}_{L^2_{\varrho}(\bbR^{\bbN})} 
}
for any $g \in \cM$. We now use the triangle inequality once more to get the result.
}

This result implies that accurate OOD generalization is possible provided (i) the constant $\Delta(\cM ; \varrho,\mu)$ is not too large, (ii) the model class $\cM$ contains an element $g$ that provides simultaneously good approximations to $f$ in both the $L^2_{\varrho}$- and $L^{\infty}_{\mu}$-norms and (iii) the in-distribution generalization error is small. Further, \ef{OOD-ID-bound} also implies stability: the OOD generalization behaviour of $\hat{f}$ will be robust to errors in the training data, as long as $\Delta(\cM ; \varrho,\mu)$ is not too large and the in-distribution generalization performance of $\hat{f}$ is robust to such errors. 

\subsection{Application to polynomial model classes}

This begs the question: for holomorphic functions, when is it possible to achieve (i)--(iii)? Our next result focuses on (i)-(ii), and shows the existence of a polynomial model class that achieves these objectives. The third objective, (iii), i.e., the construction of a suitable estimator, is the focus of the next section.

\thm{
\label{t:OOD-s-term}
Let $\bm{\omega} \geq \bm{1}$, $\varrho = \cU([-1,1]^{\infty})$, $\bm{\omega} \geq \bm{1}$ and  $D_{\bm{\omega}}$ be as in \ef{Dr-def}. Suppose that $f \in \cH(\bm{b},\varepsilon)$, where $\bm{b} \geq \bm{0}$ and $\varepsilon > 0$ satisfy
\be{
\label{b-r-cond-main}
\bm{b} \odot \bm{r}_p \in \ell^p(\bbN),\quad \nm{\bm{b} \odot (\bm{r}_p-\bm{1}) }_1 < \varepsilon,\qquad \text{where } \bm{r}_p = \zeta(\bm{\omega})^{2/p-1}
}
for some $0 < p <1$. Then, for every $k > 0$, there exists a set $S \subset \cF$ depending on $\bm{b}$, $\varepsilon$, $\bm{\omega}$ and $k$ with $|S| \leq k$ such that
\be{
\label{dist-shift-const-polys}
\Delta(\cP_{S} ; \varrho,\mu) \leq \sqrt{k},
}
where $\cP_S$ is as in \ef{f-S-def}, and a $g \in \cP_S$ such that
\be{
\label{alg-best-k-OOD}
\Delta(\cP_S ; \varrho,\mu)\nm{f - g}_{L^2_{\varrho}(\bbR^{\bbN})} + \nm{f - g}_{L^{\infty}_{\mu}(\bbR^{\bbN})} \leq C(\bm{b},\varepsilon,\bm{\omega},p) k^{1-1/p}
}
for any probability measure $\mu$ supported in $D_{\bm{\omega}}$.
}

This is a  key result of the paper. It provides a precise criterion \ef{b-r-cond-main} under which one achieves the algebraic rate \ef{alg-best-k-OOD} with a distribution shift constant \ef{dist-shift-const-polys} growing mildly in $k$. Notice that the model class depends on the anisotropy parameters $\bm{b}$ and $\varepsilon$, which makes it unsuitable for the unknown anisotropy setting (recall (i) in \S \ref{ss:contributions}). In the next section, we shall remove this requirement through a mild strengthening of \ef{b-r-cond-main}, which is equivalent to \ef{b-r-cond-intro}.

\rem{[The blessing of high dimensionality]
\label{rem:blessing}
Let $f : \bbR \rightarrow \bbR$ be holomorphic in $\cE_{\rho}$, $n \in \bbN$, $\cM =  \bbP_{n}$ be the space of polynomials of degree at most $n$ and suppose that $\mathrm{supp}(\mu) \subseteq [-\omega,\omega]$. This case was considered in depth in \cite{demanet2019stable}. As shown therein, the constant $\Delta(\cP_{S} ; \varrho,\mu) \lesssim n \zeta(\omega)^n $ and there exists a $g \in \cM$ such that
\bes{
\Delta(\cP_{S} ; \varrho,\mu) \nm{f - g}_{L^2_{\varrho}(\bbR)} + \nm{f - g}_{L^{\infty}_{\mu}} \lesssim_{\rho,\omega} n (\zeta(\omega)/\rho)^n  .
}
Thus extrapolation is possible for any $\omega$ satisfying $\zeta(\omega) < \rho$, with an exponential rate of convergence in $n$. But this comes at the expense of a potentially exponentially-large stability constant $\Delta(\cP_{S} ; \varrho,\mu) \lesssim n \zeta(\omega)^n$. To compare with Theorem \ref{t:OOD-s-term}, let $k > 0$ be arbitrary and set $n = \lfloor \frac{\log(k)}{2 \log \zeta(\omega)} \rfloor$ so that $\Delta(\cP_{S} ; \varrho,\mu) = \ord{\sqrt{k}}$. Here, the symbol $\ord{\cdot}$ suppresses a log factor in $k$ and $\omega$-dependent constant (which are irrelevant to the present discussion). With this choice of $n$, one ensures a similar stability level. However, this reduces the convergence rate 
\bes{
\Delta(\cP_{S} ; \varrho,\mu) \nm{f - g}_{L^2_{\varrho}(\bbR)} + \nm{f - g}_{L^{\infty}_{\mu}} = \ord{k^{\frac12(1-\log(\rho) / \log(\zeta(\omega)))}}
}
to algebraic in $k$.
See \cite{demanet2019stable} for a similar derivation, as well as lower bounds showing that such convergence is optimal. Here lies the blessing of high dimensionality. In low (in particular, one) dimension, the in-distribution error decays exponentially fast. However, stable extrapolation can only be achieved by reducing the rate from exponential to  algebraic. Conversely, in infinite dimensions, the in-distribution error decays algebraically, and stable extrapolation can be achieved with algebraic rates, albeit with a potentially smaller order.
}

\section{Sample-efficient OOD learning of high-dimensional functions} \label{s:sample-efficient-learning}

The previous results show that there are (polynomial) model classes that control the OOD error of an estimator $\hat{f}$ up to its in-distribution error. However, they do not describe how to compute such an estimator from data and, critically, how much data is required to do so. This is the focus of the present section, which contains the main results in this paper on function approximation. In \S \ref{s:operator-learning}, we establish analogous results for operator learning.
We now draw $\bm{x}_1,\ldots,\bm{x}_m \sim_{\mathrm{i.i.d.}} \varrho$ and consider the training data
\be{
\label{training-data-intro}
(\bm{x}_i, y_i ),\ i = 1,\ldots,m,\qquad \text{where }y_i = f(\bm{x}_i) + e_i,
}
and $e_i \in \bbR$. We then consider the $\ell^2$-loss over some model class $\cM$, yielding the estimator $\hat{f}$ in \ef{l2-loss-intro}. We first consider polynomial model classes. Next, we use emulation techniques to establish similar results for DNN model classes.

\subsection{Polynomial model classes}\label{ss:polynomial-models}

It is possible to obtain an estimator by directly combining Theorem \ref{t:OOD-s-term}, which identifies a suitable polynomial space $\cP_{S}$, with rather standard tools from the analysis of linear least-squares fitting with i.i.d.\ random samples. However, such a result is 
unrealistic, since the model class (and consequently the estimator) requires knowledge of the anisotropy parameters $(\bm{b},\varepsilon)$. Our main result establishes the existence of an estimator that is independent of $(\bm{b},\varepsilon)$, and therefore suitable for the unknown anisotropy setting described in Remark \ref{rem:unknown-anisotropy}. For this, we need the following assumption, which slightly strengthens \ef{b-r-cond-main}:
\be{
\label{b-r-cond-2-main}
\bm{b} \odot \bm{r}_p \in \ell^p_{\mathsf{M}}(\bbN),\quad \nm{\bm{b} \odot  (\bm{r}_p-\bm{1}) }_1 < \varepsilon,\qquad \text{where } \bm{r}_p = \zeta(\bm{\omega})^{2/p-1}.
}
Recall from \S \ref{ss:main-res-intro} that $\ell^p_{\mathsf{M}}(\bbN)$ is the monotone $\ell^p$-space.

\thm{
\label{thm:poly-optimized-least-squares}
Let $0 < \epsilon < 1$, $m \geq \bar{m}$, where $\bar{m} = \bar{m}(\epsilon) \in \bbN$ depends on $\epsilon$ only, $\varrho = \cU([-1,1]^{\bbN})$, $\bm{\omega} \geq \bm{1}$ and $D_{\br}$ be as in \ef{Dr-def}. Then there exists a polynomial model class $\cP$ depending on $m$, $\bm{\omega}$ and $\epsilon$ only with the following property. Let $f \in \cH(\bm{b},\varepsilon)$, where $\bm{b} \geq \bm{0}$ and $\varepsilon > 0$ satisfy \ef{b-r-cond-2-main} for some $0 < p <1$, draw $\bm{x}_1,\ldots,\bm{x}_m \sim_{\mathrm{i.i.d}} \varrho$ and consider \ef{l2-loss-intro} with noisy training data \ef{training-data-intro}. Then, with probability at least $1-\epsilon$, any minimizer $\hat{f}$ satisfies
\bes{
\nm{f - \hat{f}}_{L^{\infty}_{\mu}(\bbR^{\bbN})} \lesssim  \frac{C(\bm{b},\varepsilon,\bm{\omega},p)}{\sqrt{\epsilon}} \left ( \frac{m}{L} \right )^{1-1/p}+ \frac{\nm{\bm{e}}_2}{\sqrt{L}}
} 
for any probability measure $\mu$ supported in $D_{\br}$, where $L = L(m,\epsilon) : = \log^4(m)+\log(1/\epsilon)$ and $\bm{e} = (e_i)^{m}_{i=1}$.
Moreover, $\hat{f} \in \cP_{S}$ for some $S$ satisfying $|S| \leq m / L$. 
}

Here and elsewhere, a \textit{polynomial model class $\cP$} is a subset of all algebraic polynomials, i.e., $\cP \subseteq \spn \left \{ \bm{x} \mapsto \bm{x}^{\bm{n}} : \bm{n} \in  \cF \right \}$, where $\bm{x}^{\bm{n}} = \prod_{i \in \bbN} x^{n_i}_{i}$ for $\bm{x} = (x_i)_{i \in \bbN}$ and $\bm{n} = (n_i)_{i \in \bbN}$.
This result shows that the algebraic OOD error rate suggested by Theorem \ref{t:OOD-s-term} can be obtained from $m$ samples via a least-squares estimator, at the loss of the log term $L(m,\epsilon)$. Further, the OOD error of the estimator is robust to $\ell^2$-bounded noise, due to the term $\nm{\bm{e}}_2 / \sqrt{L}$.

\rem{
[Condition \ef{b-r-cond-main} versus \ef{b-r-cond-2-main}]
These conditions are identical whenever the sequence $\bm{b} \odot  \bm{r}_p$ is monotonically nonincreasing. Observe that in the in-distribution setting (i.e., $\bm{\omega} = \bm{1}$), they reduce to $\bm{b} \in \ell^p(\bbN)$ and $\bm{b} \in \ell^p_{\mathsf{M}}(\bbN)$, respectively. It has been shown \cite{adcock2024optimal} that it is impossible to construct estimators that converge as $m \rightarrow \infty$ for all $\bm{b} \in \ell^p(\bbN)$. However, it is possible to construct an estimator for all $\bm{b}$ satisfying $\bm{b} \in \ell^p_{\mathsf{M}}(\bbN)$. Theorem \ref{thm:poly-optimized-least-squares} shows that it is also possible to construct an estimator in the OOD setting that converges for all $\bm{b}$ satisfying \ef{b-r-cond-2-main}. It is currently unknown whether \ef{b-r-cond-2-main} can be weakened. However, as we discuss in Remark \ref{rem:rate-optimality} later, under this condition the estimator developed in Theorem \ref{thm:poly-optimized-least-squares} yields rates that are optimal up to the log term $L$ for specific choices of $\bm{\omega}$.
}

\rem{
[Practicality of the estimator]
\label{rem:poly-est-impractical}
Theorem \ref{thm:poly-optimized-least-squares} is based on a sparse polynomial approximation procedure, where $\cP =\cup_{S \in \cS} \cP_{S}$ is a union of subspaces of dimension at most $m / L(m,\epsilon)$. Roughly speaking, the collection $\cS$ is constructed by selecting all possible index sets $S$ asserted by Theorem \ref{t:OOD-s-term} for any $(\bm{b}$, $\varepsilon)$ satisfying \ef{b-r-cond-2-main}. Being a union of polynomial subspaces, it shares some similarities with the estimator of \cite{cohen2017discrete}, although the choice of sets $S$ is different. The strengthened condition \ef{b-r-cond-2-main} is required over \ef{b-r-cond-main} to ensure that $|\cS| < \infty$. This is a necessary component of our analysis which, in general, uses compressed sensing techniques. Unfortunately, the estimator in Theorem \ref{thm:poly-optimized-least-squares} is generally impractical to compute, since, in principle, it would involve solving $|\cS|$ linear least-squares problems and then choosing the one with the smallest value of the loss function. This is, in some senses, analogous to the classical problem of a sparse approximate solution to a linear system $A z = b$ by minimizing $\nm{A z - b}^2_2$ over the set of $s$-sparse vectors (since this set is a union of $s$-dimensional subspaces). In standard compressed sensing, one obtains practical algorithms by via convexification -- for instance, by minimizing the LASSO penalty $\lambda \nm{z}_1 + \nm{A z - b}^2_2$. A similar approach could be considered here, by using a weighted $\ell^1$-norm and following ideas from \cite{adcock2022sparse}. However, this is not the main focus of this work.
}

\subsection{Deep learning and DNN model classes}\label{s:deep-learning}

We now focus on DNN model classes. Let $\sigma : \bbR \rightarrow \bbR$ be an activation function. In this work, we consider feedforward DNNs of the form 
\be{
\label{Phi_NN_layers}
N : \bbR^{n_0} \rightarrow \bbR^{n_{M+1}},\ \bm{z} \mapsto N(\bm{z}) = \cA_{M} ( \sigma ( \cA_{M-1} ( \sigma ( \cdots \sigma ( \cA_0 (\bm{z}) ) \cdots ) ) ) ),
}
where $\cA_l : \bbR^{n_{l}} \rightarrow \bbR^{n_{l+1}}$ are affine maps, $n_1,\ldots,n_{M}$ are the widths of the hidden layers, and $n_0$ and $n_{M+1}$ are the input and output dimensions, respectively. We define $\mathrm{width}(N) = \max \{ n_1,\ldots, n_{M} \}$ and $\mathrm{depth}(N) = M$.
We denote a DNN model class of the form \eqref{Phi_NN_layers} with a fixed architecture (i.e., fixed activation function, depth and widths) as $\cN$, and write $\mathrm{width}(\cN) = \max \{ n_1,\ldots, n_{M} \}$ and $\mathrm{depth}(\cN) = M$.

Let $\cN$ be a DNN model class as above, with input dimension $n_0 = n \in \bbN$ and output dimension $n_{M+1} = 1$. Given training data \ef{training-data-intro} of an unknown function $f$, we consider estimator \ef{l2-loss-intro} with $\cM = \cN$. Notice that the variable $\bm{x} \in \bbR^{\bbN}$ while the input dimension of any $N \in \cN$ is $n_0 < \infty$. We assume that $N$ acts on the first $n_0$ entries of $\bm{x}$ only, but to avoid unnecessary notation, we will simply write $N(\bm{x})$ instead of $N(\bm{x}_{[n_0]})$.

\thm{
\label{thm:deep-learning-extrap}
Let $0 < \epsilon < 1$, $m \geq \bar{m}$, where $\bar{m} = \bar{m}(\epsilon) \in \bbN$ depends on $\epsilon$ only, $\varrho = \cU([-1,1]^{\bbN})$, $\bm{\omega} \geq \bm{1}$ satisfy $\omega_{\min} : = \min \{ \omega_{i} \} > 1$ and $D_{\bm{\omega}}$ be as in \ef{Dr-def}. Then there exists a DNN model class $\cN$ with tanh activation functions depending on $m$, $\bm{\omega}$ and $\epsilon$ only with the following property. Let $f \in \cH(\bm{b},\varepsilon)$, where $\bm{b} \geq \bm{0}$ and $\varepsilon > 0$ satisfy \ef{b-r-cond-2-main} for some $0 < p <1$, draw $\bm{x}_1,\ldots,\bm{x}_m \sim_{\mathrm{i.i.d}} \varrho$ and consider the problem \ef{l2-loss-intro} with noisy training data \ef{training-data-intro}. Then, with probability at least $1-\epsilon$, any minimizer $\hat{f}$ satisfies
\bes{
\nm{f - \hat{f}}_{L^{\infty}_{\mu}(\bbR^{\bbN})} \lesssim  \frac{C(\bm{b},\varepsilon,\bm{\omega},p)}{\sqrt{\epsilon}} \left ( \frac{m}{L} \right )^{1-1/p}+ \frac{\nm{\bm{e}}_2}{\sqrt{L}} + 2^{-m},
} 
for any $\mu$ probability measure supported on $D_{\bm{\omega}}$, where $L = L(m,\epsilon) : = \log^4(m)+\log(1/\epsilon)$ and $\bm{e} = (e_i)^{m}_{i=1}$.
Further, the class $\cN$ satisfies
\bes{
n_0 \lesssim \frac{m}{L},\qquad \mathrm{width}(\cN) \lesssim \frac{m}{\log(\omega_{\min})}, \qquad \mathrm{depth}(\cN) \lesssim \log \left ( \frac{\log(m)}{\log(\omega_{\min})} \right ).
}
}
This result shows one can achieve the same error bounds as in Theorem \ref{thm:poly-optimized-least-squares}, up to an extra $\ord{2^{-m}}$ term, via deep learning with DNNs that are not excessively wide or deep. This class uses tanh activations, yet other activations could readily be considered (such as ReLU or general activations that are definable in an \emph{o-minimal} sense \cite{kratsios2026algorithmic}, for which the emulation results employed in the proof of Theorem~\ref{thm:deep-learning-extrap} hold). It also assumes that $\omega_{\min} > 1$, although we believe this could be removed by using somewhat deeper DNNs. 

Theorem \ref{thm:deep-learning-extrap} employs a specific `handcrafted' family of DNNs. In the final layer, the weights can take arbitrary real values, but in the hidden layers, only a certain discrete set of values is allowed for the weights and biases (see the proof for the full construction). This is because its proof is based on combining Theorem \ref{thm:poly-optimized-least-squares} with tools from DNN emulation. Specifically, we approximately emulate the elements of the polynomial model class $\cP$ as tanh DNNs, then show that minimizers of the DNN training problem yield approximate minimizers of the polynomial training problem. Finally, the result follows through a careful perturbation analysis to show the overall effect of the approximate emulation is at most $\ord{2^{-m}}$.

\subsection{Examples and further discussion}\label{ss:condition-examples}

We now discuss the condition \ef{b-r-cond-2-main} and the various rates described by the above theorems. We do this by considering the following two examples.

\examp{[Constant side lengths]
\label{ex:const-side lengths}
Suppose first that $\omega_1 = \omega_2 = \cdots = \omega > 1$, so that $D_{\bm{\omega}}$ is a hypercube with side lengths  $2 \omega$. Then \ef{b-r-cond-2-main} holds if and only if
\be{
\label{b-eps-cond-const}
\bm{b} \in \ell^p_{\mathsf{M}}(\bbN)\quad \text{and}\quad \omega < \omega_*(\bm{b},\varepsilon,p) : = \zeta^{-1} \left ( \left ( 1 + \varepsilon/\nm{\bm{b}}_1 \right )^{\frac{p}{2-p}} \right ).
}
Note here that $\zeta^{-1} : [1,\infty) \rightarrow [1,\infty)$ is the \textit{Joukowsky map} $\zeta^{-1}(r) = \frac12\left (r+r^{-1} \right)$. Observe that the first condition $\bm{b} \in \ell^p_{\mathsf{M}}(\bbN)$ 
implies that the in-distribution error satisfies
\bes{
\nm{f - \hat{f}}_{L^{\infty}_{\mu}(\bbR^{\bbN})} \lesssim \frac{C(\bm{b},\varepsilon,p)}{\sqrt{\epsilon}} \left ( \frac{m}{L} \right )^{1-1/p} + \frac{\nm{\bm{e}}_2}{\sqrt{L}}
}
(this follows from Theorem \ref{thm:poly-optimized-least-squares} with $\omega = 1$). Hence, Theorems \ref{thm:poly-optimized-least-squares} and \ref{thm:deep-learning-extrap} show that one can achieve \textit{precisely the same rate} for the OOD error, for any measure $\mu$ supported in $[-\omega,\omega]^{\bbN}$. Notice that this goes significantly beyond the realm of the usual distribution shift theory described in \S \ref{ss:OOD-ML}, since such a $\mu$ need not be a small perturbation of $\varrho$.
}

\examp{[Unbounded side lengths]
\label{ex:unbounded-side lengths}
Now suppose that $b_i = c i^{-a}$ and $\omega_i = i^b$ for $a,b , c > 0$. Then \ef{b-r-cond-2-main} holds for $p \in (0,1)$ if and only if $a$, $b$ and $\varepsilon$ satisfy
\be{
\label{b-eps-cond-unbounded}
\frac1p < \frac{a+b}{1+2b}\quad \text{and} \quad c \sum^{\infty}_{i=1} i^{-a} \left ( \zeta(i^b)^{2/p-1} - 1 \right ) < \varepsilon.
}
If this holds, the OOD error rate is $\ord{(m/L)^{1 - \frac{a+b}{1+2b}+\delta}}$ 
for any $\delta > 0$. Conversely, the in-distribution rate is $\ord{(m/L)^{1-a+\delta}}$. 
Hence we achieve algebraic rates of decay of the OOD error for \emph{any measure supported in a hypercube whose side lengths can become arbitrarily large}, albeit at a slower algebraic rate than the in-distribution error. This setting is a stark contrast to theory of small distribution shifts discussed in \S \ref{ss:OOD-ML}.
}


\rem{
[Optimality of the rates]
\label{rem:rate-optimality}
It is known that the in-distribution error cannot decay any faster than $m^{1-1/p} / \log(m)$ in the worst case for $\bm{b} \in \ell^p_{\mathsf{M}}(\bbN)$, $\nm{\bm{b}}_{p,\mathsf{M}} = 1$, \cite[Thm.\ 4.2]{adcock2024optimalb}. This is a type of minimax lower bound, which considers the worst-case error over $\cH(\bm{b},\varepsilon)$ given arbitrary i.i.d.\ random samples and an arbitrary (nonlinear) estimator. As a result of Example \ref{ex:const-side lengths}, the estimators of Theorems \ref{thm:poly-optimized-least-squares} and \ref{thm:deep-learning-extrap} yield optimal OOD error rates (up to the log term $L$) over measures $\mu$ supported in a domain with constant side lengths. Currently, it is unknown whether the algebraic rates in the setting of Example \ref{ex:unbounded-side lengths} are also optimal.
}

\rem{
[Ill-conditioning and best guaranteed accuracy]
\label{rem:condition-number}
Our results show that the estimator is robust to perturbations, yet, as noted in \S \ref{ss:related}, analytic continuation in one dimension is ill-conditioned with infinite condition number. We now explain this seeming contradiction. The key point is that the former statements pertain to $\ell^{\infty}$-norm perturbations in the data. Suppose now that the errors $e_i$ in the data satisfy a uniform bound $\nm{\bm{e}}_{\infty} \leq \chi$, implying that $\nm{\bm{e}}_2 \leq \sqrt{m} \chi$. Then, applying Theorem \ref{thm:poly-optimized-least-squares} with an optimal choice of $m$ to balance the two error terms, we deduce that one can compute an $\hat{f}$ that satisfies (for a possibly different constant)
\be{
\label{best-accuracy}
\nm{f - \hat{f}}_{L^{\infty}_{\mu}(\bbR^{\bbN})} \lesssim C(\bm{b},\varepsilon,\bm{\omega},p) \chi^{\frac{1/p-1}{1/p-1/2}}.
}
For this argument to imply that analytic continuation had a finite condition number, this bound would have had to scale like $\chi$. But the power in \ef{best-accuracy} is less than one, as indeed it must be so as not to contradict the one-dimensional result.

Viewing $\chi = \epsilon_{\mathsf{mach}}$ as machine epsilon, the bound \ef{best-accuracy} can also be interpreted as the best guaranteed accuracy in finite-precision computations. It is interesting to note that this is some fractional power of machine epsilon that depends on the smoothness, and setting $\bm{\omega} = \bm{1}$, this is even the case for the in-distribution error. We believe this is an artefact of the estimator, which is based on the $\ell^2$-loss, or its analysis. Determining whether \ef{best-accuracy} represents the best achievable accuracy for \textit{any} estimator is an interesting topic for future work.
}

\section{Operator learning}\label{s:operator-learning}

Operator learning is an extension of the classical function approximation problem, where the target object is an operator between two, typically infinite-dimensional, spaces. As is customary, we consider operators $F : \cX \rightarrow \cY$ mapping between two separable Hilbert spaces $\cX,\cY$. We endow $\cX$ with a probability measure $\nu$ and, given $X_1,\ldots,X_m \sim_{\mathrm{i.i.d.}} \nu$, consider training data
\be{
\label{training-data-opl}
(X_i,Y_i),\ i = 1,\ldots,m,\quad \text{where } Y_i = F(X_i) + E_i \in \cY
}
and $E_i \in \cY$ represents noise. The task is to approximate $F$ from the data \ef{training-data-opl}.

Operator learning is often carried out using DNOs. These are extensions of DNNs that are designed to handle the infinite-dimensional nature of the inputs and outputs. There are many different operator learning strategies, such as DeepONets \cite{lu2021learning}, Fourier Neural Operators (FNOs) \cite{li2021fourier}, PCA-Net \cite{lanthaler2023operator} and many others. We refer to \cite{brugiapaglia2026short,boulle2024mathematical,kovachki2024operator,subedi2026operator} for reviews.

A standard operator learning paradigm, termed \textit{encoder-decoder nets}, involves approximating $F$ with a DNO of the form
\be{
\label{encoder-decoder}
F \approx \hat{F} : = \hat{\cD}_{\cY} \circ \hat{N} \circ \hat{\cE}_{\cX},
}
where $\hat{\cE}_{\cX} : \cX \rightarrow \bbR^{d_{\cX}}$ in an \textit{encoder} for $\cX$, $\hat{\cD}_{\cY} : \bbR^{d_{\cY}} \rightarrow \cY$ is a \textit{decoder} for $\cY$ and $\hat{N} : \bbR^{d_{\cX}} \rightarrow \bbR^{d_{\cY}}$ is a neural network. The encoder and decoder are either specified by the problem, learned spearately from data or learned concurrently with $\hat{N}$. Note that DeepONets and PCA-Net both fall into this category, as do various other approaches. Our main theoretical results consider certain types of encoder-decoder DNOs.

\subsection{Assumptions and definitions}\label{ss:opl-ass}

We now introduce the key assumptions needed for our theoretical results.

\begin{assumption}
\label{ass:opl}
The probability measure $\nu$ has mean zero and finite second moments, i.e., $\int_{\cX} \nm{X}^2_{\cX} \D \mu(X) < \infty$. Further, let $\lambda_1 \geq \lambda_2 \geq \cdots \geq 0$ and $\{ \phi_i \}_{i \in \bbN} \subset \cX$ be the eigenvalues and orthonormal basis of eigenvectors of the covariance operator of $\nu$. Then the real-valued random variables $\xi_i = \ip{X}{\phi_i}_{\cX} / \sqrt{3\lambda_i}$ satisfy $\xi_i \sim_{\mathrm{i.i.d.}} \cU([-1,1])$.
\end{assumption}

Here $\sqrt{3}$ is a normalization factor, as $\bbE \xi^2_i = 1/3$.
Assumption \ref{ass:opl} means that $\nu$ is the law of the $\cX$-valued random variable $X = \sum^{\infty}_{i=1} \sqrt{3\lambda_i} \xi_i  \phi_i$, where $\xi_i \sim_{\mathrm{i.i.d.}} \cU([-1,1])$.
Such an assumption is quite common in practice. Indeed, it is standard in operator learning literature to generate training data by simulating inputs $X$ via such an expansion \cite{boulle2024mathematical,kovachki2024operator,subedi2026operator}. We slightly deviate from standard practice by assuming the $\xi_i$ are uniformly-distributed random variables, rather than standard normal random variable. In the latter case, $\nu$ is a Gaussian process on $\cX$. Dealing with normal random variables presents some technical challenges, which are beyond the scope of this paper. See \S \ref{s:conclusion} for further discussion.

\begin{assumption}
\label{ass:opl-2}
Let the $\phi_i$ and $\lambda_i$ be as in Assumption \ref{ass:opl}. The encoder $\hat{\cE}_{\cX}$ is given by $\hat{\cE}_{\cX} : X \mapsto (\ip{X}{\phi_i}_{\cX} / \sqrt{3\lambda_i})^{d_{\cX}}_{i=1}$ and the decoder $\hat{\cD}_{\cY}$ is linear.
\end{assumption}

This assumption means that our DNO architecture is a type of PCA-Net \cite{lanthaler2023operator}, at least in terms of its encoder, since $\hat{\cE}_{\cX}$ involves analyzing $X \in \cX$ using the first $d_{\cX}$ PCA basis functions. In PCA-Net, one normally first learns the PCA basis empirically from data before learning the DNN $\hat{N}$.
To avoid complications, we assume the exact PCA basis is known. However, we believe it is possible to extend our analysis to the case where an empirical PCA basis is employed by tracking the influence of this error on the overall error bound. This is a topic for future work. 
In PCA-Net, the decoder $\hat{\cD}_{\cY}$ is also constructed using the empirical PCA basis on $\cY$ specified by the pushforward measure $F \sharp \nu$. However, we do not require this, since we just assume linearity of the map. 

\begin{assumption}
\label{ass:opl-ood}
The probability measure $\mu$ on $\cX$ satisfies the following. Let $\upsilon$ denote the law of the random variable $\bm{\xi}$ on $\bbR^{\bbN}$, defined as $\bm{\xi} = \left ( \ip{X}{\phi_i}_{\cX} / \sqrt{3\lambda_i} \right )_{i \in \bbN}$, where $X \sim \mu$. Then $\mathrm{supp}(\upsilon) \subseteq D_{\bm{\omega}}$ for some $\bm{\omega} \geq \bm{1}$, where $D_{\bm{\omega}}$ is as in \ef{Dr-def}.
\end{assumption}

This states that the projection of samples $X \sim \mu$ via the PCA basis $\{ \phi_i \}_{i \in \bbN}$, scaled by the PCA eigenvalues $\{ \lambda_ i \}_{i \in \bbN}$, are bounded, and the $i$th such projection $\ip{X}{\phi_i}_{\cX} / \sqrt{3\lambda_i}$ belongs to $[-\omega_i,\omega_i]$ almost surely. Note that we do not require the random variables $\ip{X}{\phi_i}_{\cX} / \sqrt{3\lambda_i}$ to be uncorrelated, hence this assumption is quite general. We discuss several examples that satisfy it in \S \ref{ss:examples-opl} below. In particular, we show that this assumption allows for functions sampled from the test distribution $\mu$ to be less smooth than functions sampled from the training distribution $\nu$ -- an important practical scenario.

Finally, we now define the class of holomorphic operators we consider in this work.

\defn{
\label{def:holo-op}
Suppose that Assumption \ref{ass:opl} holds. We say that $F : \cX \rightarrow \cY$ is a \textit{$(\bm{b},\varepsilon)$-holomorphic operator} if $F(X) = f (  \ip{X}{\phi_i}_{\cX} )^{\infty}_{i=1} / \sqrt{3\lambda_i} )$, $\forall X \in \cX$, and $f : \bbR^{\bbN} \rightarrow \cY$ is a $(\bm{b},\varepsilon)$-holomorphic $\cY$-valued function. We write $\cH(\bm{b},\varepsilon ; \cY)$ for the class of $\cY$-valued $(\bm{b},\varepsilon)$-holomorphic functions and $\cH(\bm{b},\varepsilon ; \cX , \cY)$ for the class of $(\bm{b},\varepsilon)$-holomorphic operators.
}

\rem{
[Infinite encoders and decoders]
\label{rem:inf-enc-dec}
It is convenient at this point to define the `infinite' encoder and decoder
\bes{
\cE_{\cX} :  X \mapsto (\ip{X}{\phi_i}_{\cX} / \sqrt{3\lambda_i} )_{i \in \bbN},\quad \cD_{\cX} :  \bm{x} \mapsto \sum_{i \in \bbN} \sqrt{3\lambda_i} x_i \phi_i.
}
Since $\{ \phi_i \}_{i \in \bbN}$ is an orthonormal basis, we have that $\cD_{\cX} \circ \cE_{\cX} = \cI_{\cX}$. Notice that Assumption \ref{ass:opl} is equivalent to $\cE_{\cX} \sharp \nu = \varrho$,
where we recall that $\varrho = \cU([-1,1])^{\bbN}$, and Assumption \ref{ass:opl-ood} can be written equivalently as $\mathrm{supp}(\cE_{\cX} \sharp \mu) \subseteq  D_{\bm{\omega}}$. Further, the operator $F$ and function $f$ in Definition \ref{def:holo-op} satisfy $F = f \circ \cE_{\cX}$ and $f = F \circ \cD_{\cX}$.
}

\subsection{Main result}

Our main results consider operator learning using encoder-decoder nets of the form \ef{encoder-decoder}. Given an encoder $\hat{\cE}_{\cX}$ and decoder $\hat{\cD}_{\cY}$ we let $\cN$ be a family of DNNs $N : \bbR^{d_{\cX}} \rightarrow \bbR^{d_{\cY}}$ and consider the estimator
\be{
\label{DNN-training-prob-operators}
\hat{F} =  \hat{\cD}_{\cY} \circ \hat{N} \circ \hat{\cE}_{\cX},\qquad \text{where }\hat{N} \in \argmin{N \in \cN} \frac1m \sum^{m}_{i=1} \nm{Y_i - \hat{\cD}_{\cY} \circ N \circ \hat{\cE}_{\cX} (X_i)}^2_{\cY}.
}
Before stating our result we need one further piece of notation. Recall from Assumption \ref{ass:opl-2} that $\hat{\cD}_{\cY} : \bbR^{d_{\cY}} \rightarrow \cY$ is linear, therefore its range $\hat{\cY} : = \hat{\cD}_{\cY}(\bbR^{d_{\cY}})$ denotes an at most $d_{\cY}$-dimensional subspace. We now write $\hat{\cB}_{\cY} : \cY \rightarrow \hat{\cY}$ for the orthogonal projection onto $\hat{\cY}$.

\thm{
\label{thm:operator-learning-extrap}
Let $0 < \epsilon < 1$ and $m \geq \bar{m}$, where $\bar{m} = \bar{m}(\epsilon)$ depends on $\epsilon$ only. Suppose that Assumptions \ref{ass:opl}-\ref{ass:opl-ood} hold with $\bm{\omega} > \bm{1}$ satisfying $\omega_{\min} : = \min \{ \omega_{i} \} > 1$.  Then there exists a DNN model class $\cN$ with $n_0 = d_{\cX}$ and $n_{M+1} = d_{\cY}$ and depending on $m$, $\bm{\omega}$ and $\epsilon$ only with the following property. Let $F \in \cH(\bm{b},\varepsilon ; \cX , \cY)$, where $\bm{b} \geq \bm{0}$ and $\varepsilon > 0$ satisfy \ef{b-r-cond-2-main} for some $0 < p < 1$, draw $X_1,\ldots,X_m \sim_{\mathrm{i.i.d.}} \nu$ and let $\hat{F}$ be as in \ef{DNN-training-prob-operators} with noisy training data \ef{training-data-opl}. Then, with probability at least $1-\epsilon$, $\hat{F}$ satisfies 
\eas{
\nm{F - \hat{F}}_{L^{\infty}_{\mu}(\cX ; \cY)} \lesssim & ~ \frac{C(\bm{b},\varepsilon,\bm{\omega},p)}{\sqrt{\epsilon}} \left ( \frac{m}{L} \right )^{1-1/p} + \sqrt{\frac{m}{\epsilon L}} \nm{F - \hat{\cB}_{\cY} \circ F }_{L^{2}_{\nu}(\cX;\cY)} 
 \\
& +\frac{\nm{\bm{E}}_{2;\cY}}{ \sqrt{L} } + 2^{-m} ,
}
provided $d_{\cX} \geq c m / L$ for some universal constant $c$, where $L = L(m,\epsilon) : = \log^4(m) + \log(1/\epsilon)$ and $\bm{E} = ( E_i)^{m}_{i=1}$. Further, the class $\cN$ satisfies
\bes{
\mathrm{width}(\cN) \lesssim \frac{m}{\log(\omega_{\min})}, \qquad \mathrm{depth}(\cN) \lesssim \log \left ( \frac{\log(m)}{\log(\omega_{\min})} \right ).
}
}
Here and elsewhere $\nm{\cdot}_{L^{p}_{\mu}(\cX ; \cY)}$ denotes the $L^p$-Bochner norm (see \S \ref{ss:additional-notation}). This result provides an OOD generalization guarantee for operator learning. It shows that the generalization error for learning holomorphic operator decays with algebraic rates of convergence in terms of $m$ depending on the support of the measure $\mu$, as described by Assumption \ref{ass:opl-ood}. We discuss examples of this condition in the next subsection.

As with previous results, the DNNs are not excessively wide or deep. Similar to Theorem \ref{thm:deep-learning-extrap}, while this result considers tanh DNNs, other activations could readily be considered. Distinct from Theorem \ref{thm:deep-learning-extrap}, however, this result describes learning classes of operators taking values in an arbitrary separable Hilbert space $\cY$. For this reason, the error bound contains an additional term that acounts for the decoding error on $\cY$: namely,
\bes{
\nm{F - \hat{\cB}_{\cY} \circ F }_{L^{2}_{\nu}(\cX;\cY)} = \sqrt{\int_{\cX} \nm{F(X) - \hat{\cB}_{\cY} \circ F(X)}^2_{\cY} \D \nu(X) } \equiv \nm{\cI_{\cY} - \hat{\cB}_{\cY} }_{L^2_{F \sharp \nu}(\cY ; \cY)}.
}
Here $\cI_{\cY} : \cY \rightarrow \cY$ is the identity map. Such an expression is quite standard in the operator learning literature. Overall, this term implies that the decoder error depends on how well $\hat{\cB}_{\cY}$ approximates the identity map. Concrete estimates for this term have been established in various settings, such as PCA-Nets and DeepONets \cite{lanthaler2023operator,lanthaler2022error}.




\subsection{Examples and further discussion}\label{ss:examples-opl}

We now discuss Assumption \ref{ass:opl-ood} in more detail.

\examp{[Rougher test distributions]
\label{ex:opl-rougher-test}
Suppose that $\mu$ has the same PCA eigenbasis as $\nu$, with PCA eigenvalues $\{ \lambda'_i \}_{i \in \bbN}$ satisfying $\lambda'_i \geq \lambda_i$, $\forall i \in \bbN$. In this case, $\mu$ is the law of the random variable $X = \sum^{\infty}_{i=1} \xi_i \sqrt{3\lambda'_i} \phi_i$, where $\xi_i \sim_{\mathrm{i.i.d.}} \cU([-1,1])$.
Note that, if $\{ \phi_i \}_{i \in \bbN}$ is an orthonormal basis of smooth functions, such as the Fourier basis, this means that functions sampled from the test distribution $\mu$ can be less smooth than functions sampled from the training distribution $\nu$.
In this case, Assumption \ref{ass:opl-ood} holds if and only if
\bes{
\lambda'_i \leq \omega_i^2 \lambda_i,\quad \forall i \in \bbN.
}
meaning that the eigenvalues $\lambda'_i$ may decay as $i \rightarrow \infty$ at a slower rate than the $\lambda_i$, depending on the growth rate of the $\omega_i$. Indeed, after choosing $\bm{\omega}$ minimally so that the above inequality holds, we obtain the OOD rate $\ord{(m/L)^{1-1/p}}$ provided $F \in \cH(\bm{b},\varepsilon ; \cX , \cY)$ for $\bm{b}$, $\varepsilon$ satisfying
\be{
\label{opl-eig-cond}
\bm{b} \odot  \bm{r}_p \in \ell^p_{\mathsf{M}}(\bbN),\quad \nm{\bm{b} \odot  (\bm{r}_p-\bm{1}) }_1 < \varepsilon,\qquad \text{where } \bm{r}_p = \left ( \zeta\left(\sqrt{\lambda'_i / \lambda_i}\right)^{2/p-1} \right)_{i \in \bbN}.
}
Concretely, let $\lambda_i = i^{-v}$ and  $\lambda'_i = i^{-u}$ for $v \geq u > 1$ and suppose, as in Example \ref{ex:unbounded-side lengths}, that $b_i = c i^{-a}$ for some $c > 0$. Then, it is a short argument to show that algebraic convergence occurs -- i.e., \ef{opl-eig-cond} holds for some $p \in (0,1)$ -- whenever $a > 1 + \frac{v-u}{2}$ and $c$ is sufficiently small. Thus, the OOD error decays algebraically even when $\mu$ is much rougher than $\nu$ (i.e., $u \ll v$), provided $F$ has sufficient regularity.
}

In the previous example, the entries $\xi_i$ of the random variable $\bm{\xi}$ defined in Assumption \ref{ass:opl-ood} were uncorrelated, since $\nu$ and $\mu$ differed only in their PCA eigenvalues. However, Assumption \ref{ass:opl-ood} allows for correlated distributions, thus allowing $\nu$ and $\mu$ to be very different in nature. We now explore several such examples

\examp{
[Arbitrarily rough test distributions]
\label{ex:opl-arbitrarily-rough}
Suppose that $\nm{X}_{\cX} \leq c$ for $X \sim \mu$ a.s.. Note that such a distribution can be arbitrarily rough. Assumption \ref{ass:opl-ood} holds, provided $\omega_i \geq c / \sqrt{3\lambda_i}$, $\forall i \in \bbN$.
In this case, Theorem \ref{thm:operator-learning-extrap} shows that one can generalize to arbitrarily rough test distributions given sufficient smoothness: namely, $(\bm{b},\varepsilon)$-holomorphy, where $\bm{b}$ satisfies \ef{b-r-cond-2-main} for some $\bm{\omega}$ whose entries $\omega_i$ grow proportional to $1/\sqrt{\lambda_i}$. Thus, and perhaps unsurprisingly, faster decay of the eigenvalues $\lambda_i$ (i.e., a smoother training distribution $\nu$) translates into higher smoothness of the operator in order to ensure convergence for an arbitrarily rough test distribution $\mu$.

For a concrete example, let $\cX = L^2([0,1])$ and $\phi_i(x) = \sqrt{2}\sin(i \pi x)$, $\forall i \in \bbN$, be the sine basis. Let $\nu$ be the law of the random variable $X = \sqrt{2}\sum^{\infty}_{i=1} \xi_i i^{-v/2} \sin(i \pi \cdot)$, where $\xi_i \sim_{\mathrm{i.i.d.}} \cU([-1,1])$,
for some $v > 1$. Now let $d > 1$ and $\mu$ be the law of the random variable $X = \sum^{d-1}_{i=0} (2+\xi_i)  \bbI_{[i/d,(i+1)/d)}(\cdot)$, where  $\xi_i \sim_{\mathrm{i.i.d.}} \cU([-1,1])$. Functions drawn from $\nu$ belong almost surely to $H^s([0,1])$ for $s < (v-1)/2$, while functions drawn from $\mu$ are pieceiwse constant and only belong to $H^s([0,1])$ almost surely for $s < 1/2$. Since $\nm{X}_{\cX} \leq 3$ for $X \sim \mu$ a.s.\ we may choose $\omega_i = 3 i^{v/2}$ in order for Assumption \ref{ass:opl-ood} to hold. 
Such choice is independent of $d$. However, a direct calculation gives that $| \ip{X}{\phi_i}_{\cX} | \leq 6\sqrt{2} d / i$,  $\forall i \in \bbN$,
almost surely for $X \sim \mu$, meaning that one can also choose $\omega_i = \max \{ 1 , 6 \sqrt{2} d i^{v/2-1} \}$. Therefore, if $d$ is not large, a more slowly-growing choice of $\omega_i$ -- which equates to a less stringent smoothness requirement on $F$ -- is also possible.
}

\section{Numerical experiments}\label{s:experiments}

We now present a series of numerical experiments. These experiments consider polynomial, DNN and DNO estimators, with full implementation details given in Appendix \ref{app:experiments}. The intention of these experiments is to illustrate the main messages of the paper: namely, OOD generalization errors decay at algebraic rates depending on the support of the test distribution. But they also elucidate a number of gaps between the theory developed in this work and practice. This motivates future work, as we discuss further in \S \ref{s:conclusion}.

\rem{
[Unknown supports]
\label{rem:unknown-supports}
The estimators we consider in our experiments are all independent of $\bm{\omega}$, the parameter defining the hyperrectangle $D_{\bm{\omega}}$ in which the support of the test distribution $\mu$ is contained. This is a key difference with our theory, since the estimators of Theorems \ref{thm:poly-optimized-least-squares}, \ref{thm:deep-learning-extrap} and \ref{thm:operator-learning-extrap} all depend on $\bm{\omega}$. An important takeaway from our experiments is that polynomial and DNN estimators possess even stronger OOD performance than the theory suggests, as they achieve algebraic rates of convergence \textit{without any knowledge} of the test distribution $\mu$. As we discuss further in \S \ref{s:conclusion}, understanding why this is the case is an open problem. 
}

\subsection{Polynomial estimators for multivariate function approximation}\label{ss:num-exp-poly}

We first consider multivariate function approximation. Our aim is to conduct a precise study of the effect of changing the support of the test distribution on the convergence rate and to compare the empirically-observed rates to the theoretical rates. Therefore, in order to avoid effects implicit to training DNNs -- such as architecture design, initializations and optimization solvers -- in this subsection we consider polynomial estimators. For the reasons discussed in Remark \ref{rem:poly-thy-practce}, we do not implement exact the estimator of Theorem \ref{thm:poly-optimized-least-squares}. Instead, we consider an efficient Adaptive Least Squares (ALS) estimator based on \cite{migliorati2019adaptive}. This estimator is described in \S \ref{app:ALS}.

\rem{[Differences between Theorem \ref{thm:poly-optimized-least-squares} and the experiments]
\label{rem:poly-thy-practce}
There are several reasons why we do not implement the estimator of Theorem \ref{thm:poly-optimized-least-squares}. First, as discussed in Remark \ref{rem:poly-est-impractical} it is generally impracticable. While it is possible to design an implementable estimator using weighted $\ell^1$-minimization, such estimators can be memory-intensive and slow to compute \cite{adcock2022sparse}, thus limiting ones ability to perform convergence studies involving large sample sizes. In comparison, the ALS estimator is relatively inexpensive to compute, yet it lacks theoretical guarantees, even in the in-distribution setting \cite{migliorati2019adaptive}. However, the experiments in this section suggest it possesses excellent in-distribution and OOD generalization properties.
}

To conduct a precise study, we consider the family of functions 
\be{
\label{delta-fun}
f : [-1,1]^d \rightarrow \bbR,\ \bm{x} = (x_i)^{d}_{i=1} \mapsto \prod^{d}_{i=1} \frac{\sqrt{2 \delta_i + \delta^2_i}}{x_i+1+\delta_i},
}
for positive parameters $\delta_i > 0$ (the factor in the numerator ensures this function has unit norm, and therefore avoids scale effects for large $d$ due to the $d$-fold product) \cite{adcock2024learning}. These functions are holomorphic in $[-1,1]^d$ with singularities at any $\bm{x}$ for which $x_i = -1-\delta_i$ for some $i$. Therefore they are holomorphic in $\cE_{\rho_1} \times \cdots \times \cE_{\rho_d}$ for any $\rho_i \geq 1$ satisfying
\be{
\label{f-delta-holo-region}
\rho_i < \zeta(1+\delta_i)\quad \forall i \in [d],
}
where we recall $\zeta$ from \S \ref{ss:main-res-intro}.
Because of Remark \ref{rem:fin-dim}, the function \ef{delta-fun} -- or, more precisely, its extension to $[-1,1]^{\bbN}$ -- is $(\bm{b},\varepsilon)$-holomorphic for any $\bm{b} = (b_i)^{\infty}_{i=1}$ satisfying
\be{
\label{b-delta-relation}
b_i > \varepsilon / \delta_i,\ \forall i \in [d],\qquad b_i = 0,\ \forall i > d.
}
In our experiments, we consider $\delta_i = 2 i^2$.
For training and and test distributions, we consider $\varrho = \cU([-1,1]^d)$ and $\mu = \cU\left ( \otimes^{d}_{i=1} [-\omega_i,\omega_i] \right )$, respectively.

\begin{figure}[t]
  \begin{minipage}[c]{0.35\linewidth}
        \centering
\includegraphics[width=\linewidth]{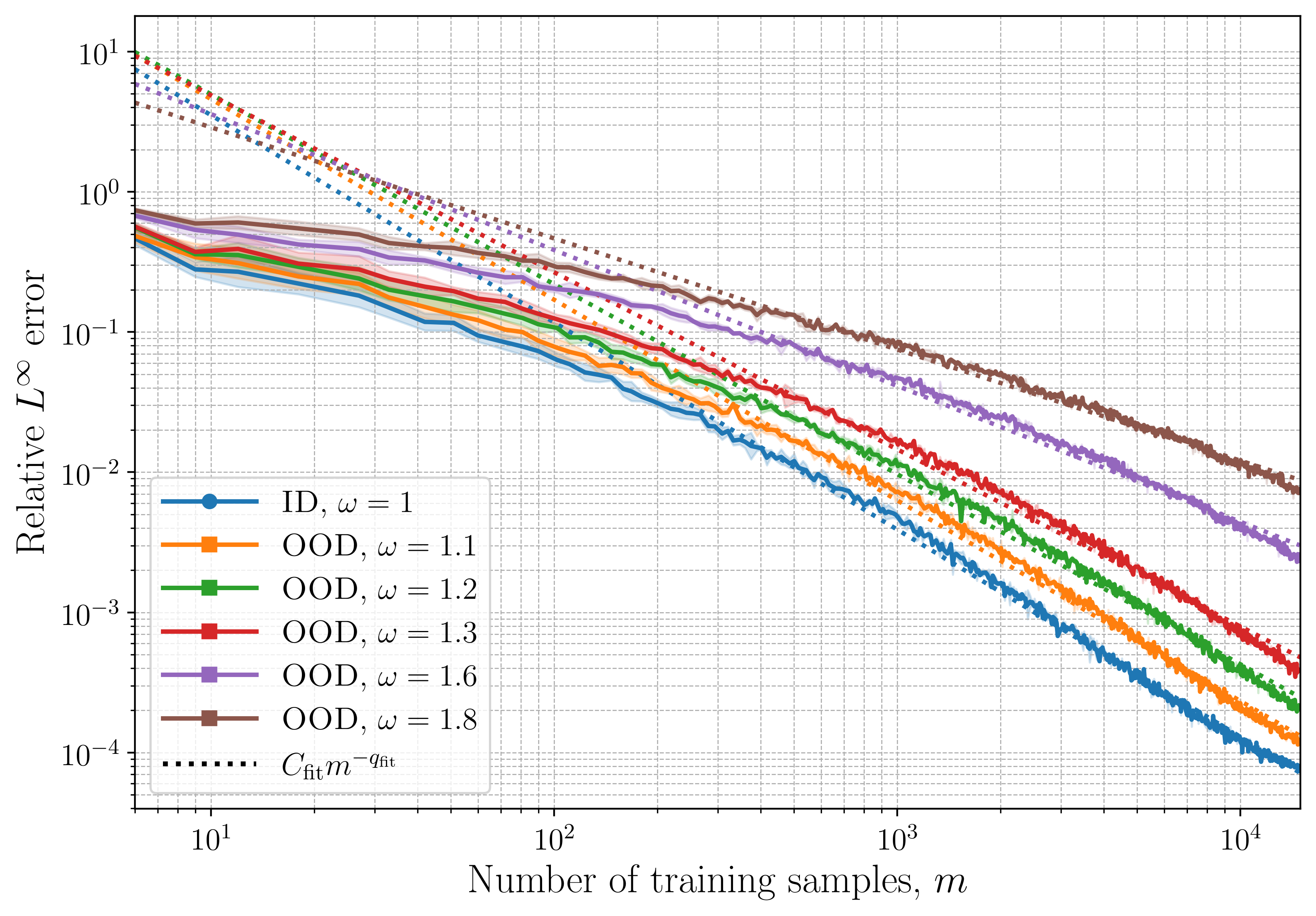}
    \end{minipage}
  \begin{minipage}[c]{0.35\linewidth}
        \centering
\includegraphics[width=\linewidth]{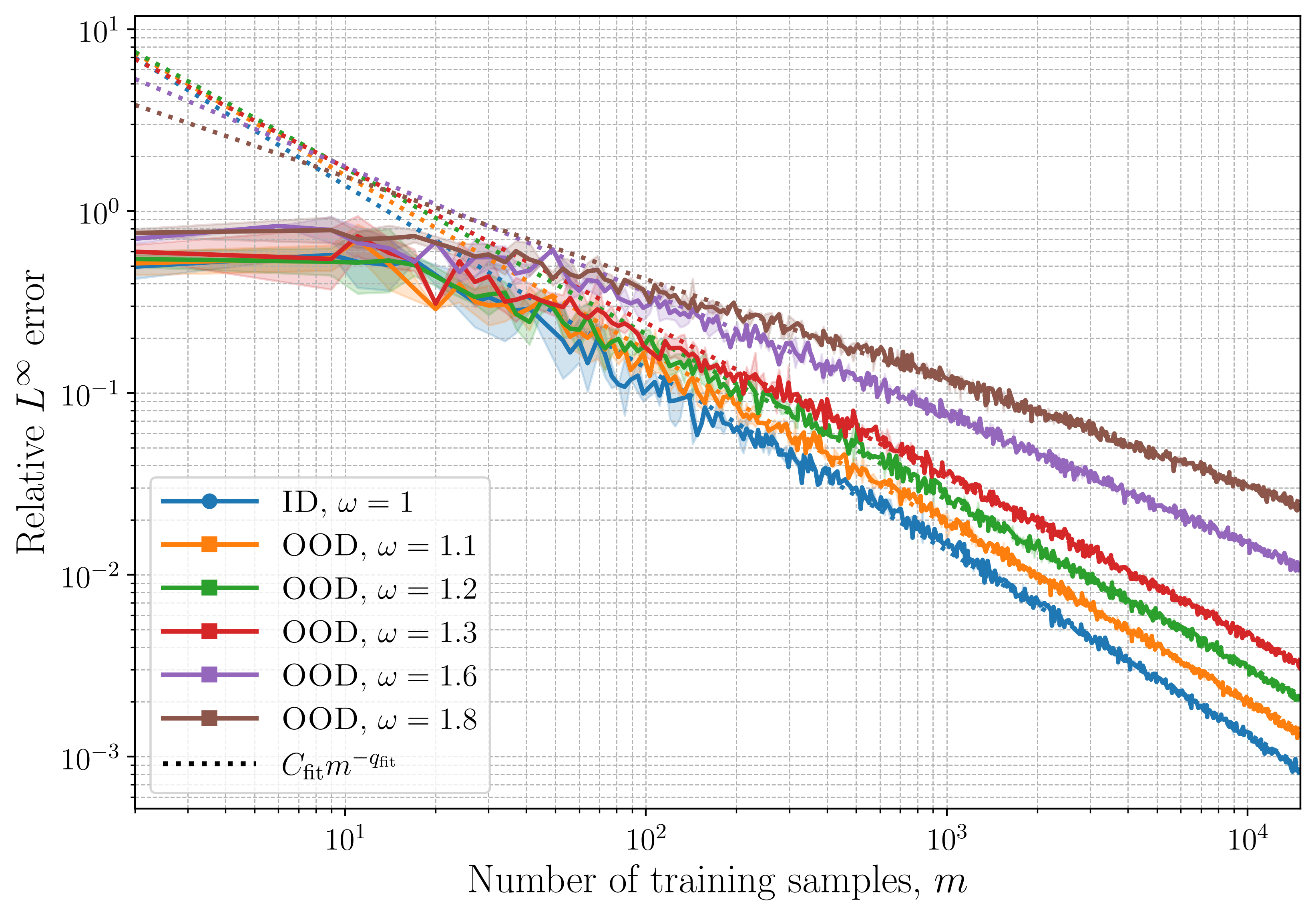}
    \end{minipage}
        \begin{minipage}[c]{0.28\linewidth}
        \centering
        {\small
\begin{tabular}{rrr}
\hline
$\omega$ & $q_{\mathrm{fit}}$ & $q_{\mathrm{thy}}$ \\
\hline
1.0 & 1.004 & 1.000 \\
1.1 & 0.960 & 0.409 \\
1.2 & 0.912 & 0.148 \\
1.3 & 0.854 & 0.033 \\
1.6 & 0.689 & --- \\
1.8 & 0.564 & --- \\
\hline
\end{tabular} 
}   \end{minipage}
    \caption{Error versus number of training samples $m$ for function approximation using polynomial estimators computed by the ALS algorithm. The function is given by \eqref{delta-fun} with $\delta_i = 2i^2$, $\forall i \in [d]$ and $d = 8 $ (left) or $d = 32$ (middle). The training and test distributions are $\varrho = \cU([-1,1]^d)$ and $\mu = \cU(\otimes^{d}_{i=1}[-\omega_i,\omega_i])$, where $\omega_1 = \cdots = \omega_d = \omega$. The value $q_{\mathrm{thy}}$ is given by \ef{q-thy-const-side lengths} in the case $d = 32$. A dashed line means that $q_{\mathrm{thy}} \leq 0$, i.e., no convergence. The values $q_{\mathrm{fit}}$ were computed by empirical fitting in the case $d=32$.
%
    }
    \label{fig:poly-constant-sides}
\end{figure}

Following Example \ref{ex:const-side lengths}, in Fig.\ \ref{fig:poly-constant-sides} we first consider constant side lengths $\omega_1 = \cdots = \omega_d = \omega \geq 1$. Since $\min_{i} \delta_i = \delta_1 = 2$,  \ef{f-delta-holo-region} implies that $f$ is holomorphic in any hypercube $[-\omega,\omega]^d$ with $\omega < 3$.
Fig.\ \ref{fig:poly-constant-sides} shows that the OOD error decays for all values of $\omega$ considered. Recalling Remark \ref{rem:fin-dim-funs}, we notice that in lower dimension $d = 8$, the rate of convergence appears slightly faster than algebraic, while when $d$ increases to $d = 32$ we witness algebraic convergence. 

In Fig.\ \ref{fig:poly-constant-sides}, we also compare the empirically-obtained algebraic rate of convergence $q_{\mathrm{fit}}$ with the rate predicted by our theory $q_{\mathrm{thy}}$. Recall that our theory predicts a rate $(m/L)^{1-1/p}$ whenever \ef{b-r-cond-2-main} holds for some $0 < p < 1$. For a given value of $\omega$, recalling \ef{b-delta-relation} and the fact that $\delta_i = 2 i^2$, we determine that
\be{
\label{q-thy-const-side lengths}
q_{\mathrm{thy}} = \frac12 \left [ \log \left ( 1 + 2 \Big /\sum^{d}_{i=1}  i^{-2}  \right ) \Big / \log(\zeta(\omega)) - 1 \right ],
}
meaning that our theory predicts a rate that is $m^{-q_{\mathrm{thy}} + \delta}$ for any $\delta > 0$.

As Fig.\ \ref{fig:poly-constant-sides} shows, the estimator consistently outperforms the theoretical convergence rate, except in the case $\omega = 1$, which corresponds to the in-distribution error. Further, after manipulating \ef{q-thy-const-side lengths}, we see that $q_{\mathrm{thy}} \leq 0$, meaning our theory predicts no convergence, whenever $\omega \geq 1.343$. Yet, the estimator still converges algebraically for $\omega = 1.6$ and $\omega = 1.8$. This suggests the main condition \ef{b-r-cond-2-main} may not be sharp -- improving it is an open problem.



Next, we consider varying side lengths. Following Example \ref{ex:unbounded-side lengths}, we set $\omega_i = i^b$ for several different values of $b$. Similar to the previous example, we witness algebraic convergence, or faster in the $d = 8$ case. We also compare the rate against our theory. In this case, the theoretical rate is given by
\be{
\label{q-thy-unbounded-side lengths}
q_{\mathrm{thy}} = 1/p_{\mathrm{thy}}-1,\quad \text{where $p_{\mathrm{thy}}$ solves } \sum^{d}_{i=1} \left ( \zeta(i^b)^{2/p-1} -1 \right ) i^{-2} = 2.
}
Once more, the empirical convergence rate is better than our theory predicts. In particular, the estimator converges for $b = 0.6$ and $b = 0.8$. This is outside the realm of our theory, which only asserts convergence for $b < 0.562$. Indeed, \ef{q-thy-unbounded-side lengths} implies that $q_{\mathrm{thy}} \leq 0$ whenever $b \geq 0.562$.
This once again points towards a potential lack of sharpness in the condition \ef{b-r-cond-2-main} and an open problem for future work.


\begin{figure}[t]
    \centering
    \begin{minipage}[c]{0.35\linewidth}
        \centering
        \includegraphics[width=\linewidth]{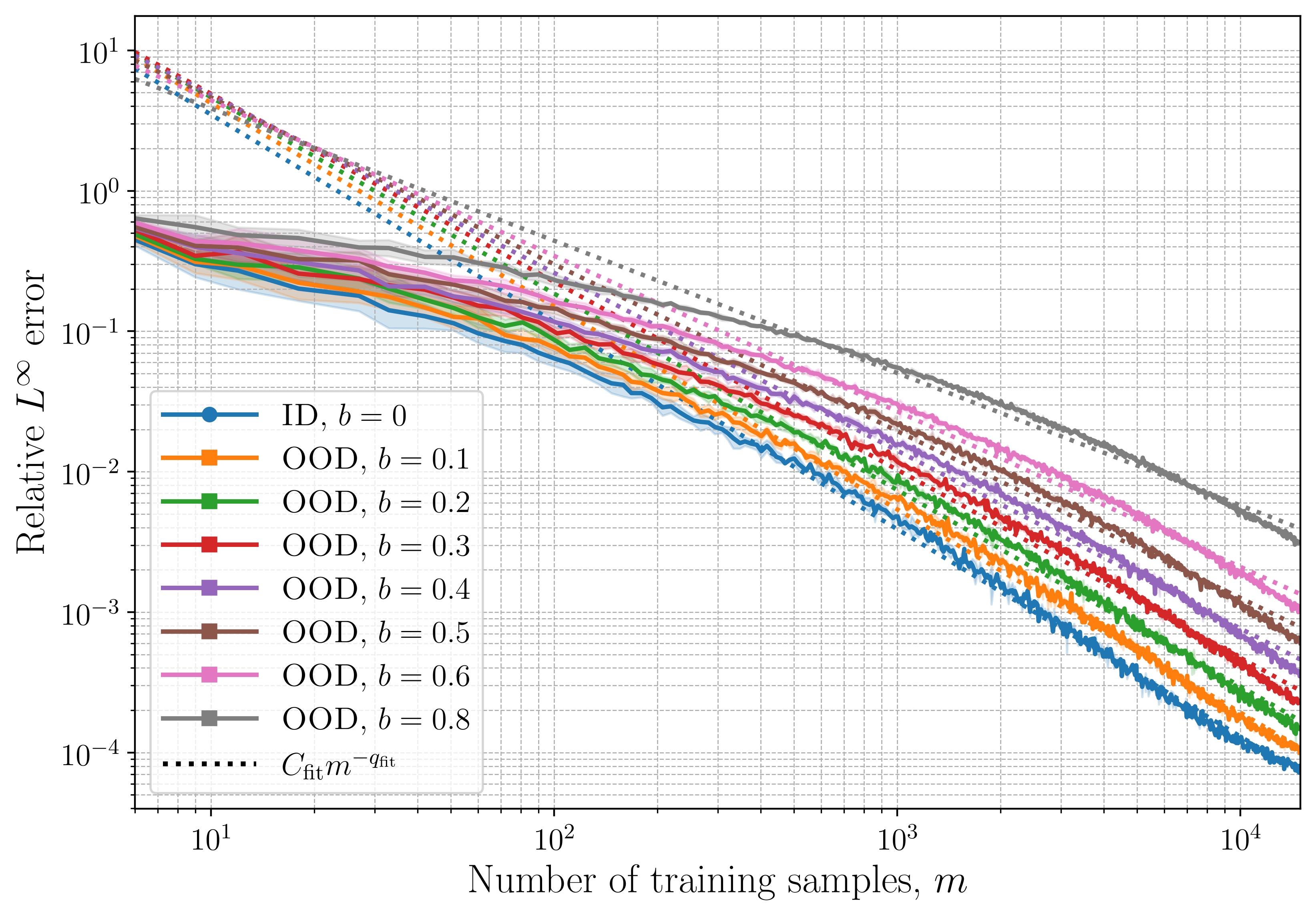}
    \end{minipage}
        \begin{minipage}[c]{0.35\linewidth}
        \centering
        \includegraphics[width=\linewidth]{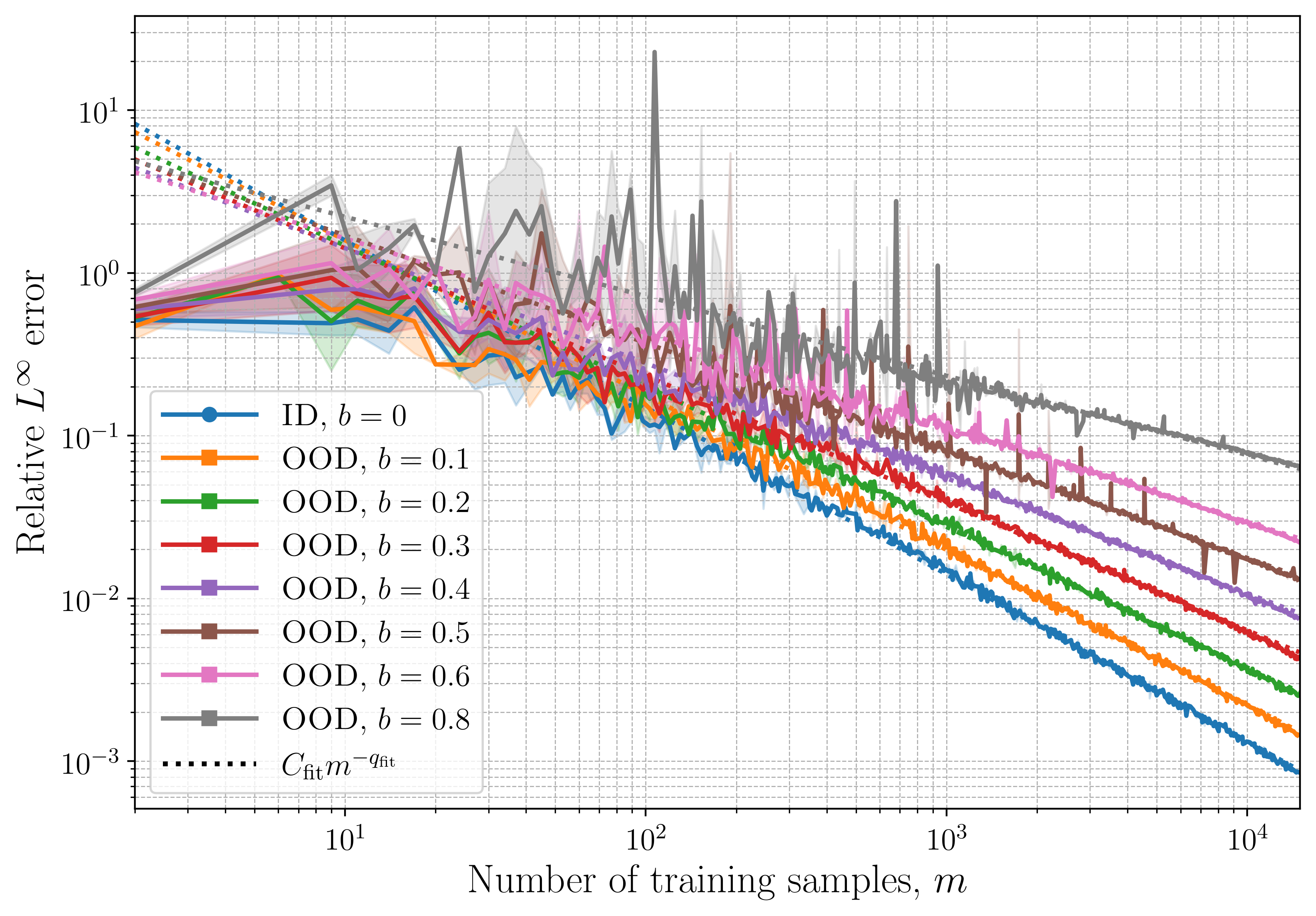}
    \end{minipage}
    \hfill
    \begin{minipage}[c]{0.28\linewidth}
        \centering
        {\small
\begin{tabular}{rrr}
\hline
$b$ & $q_{\mathrm{fit}}$ & $q_{\mathrm{thy}}$ \\
\hline
0.0 & 1.026 & 1.000 \\
0.1 & 0.951 & 0.855 \\
0.2 & 0.864 & 0.430 \\
0.3 & 0.782 & 0.237 \\
0.4 & 0.707 & 0.120 \\
0.5 & 0.660 & 0.039 \\
0.6 & 0.581 & --- \\
0.8 & 0.488 & --- \\
\hline
\end{tabular} 
}   \end{minipage}
\caption{The same as Fig.\ \ref{fig:poly-constant-sides}, except with varying side lengths $\omega_i = i^b$ and the value $q_{\mathrm{thy}}$ given by \ef{q-thy-unbounded-side lengths} in the case $d = 32$ once more.}
    \label{fig:poly-unbounded-sides}
\end{figure}

%

\subsection{DNN estimators for multivariate function recovery}\label{ss:ex-DNNs}

We now consider DNN estimators. Rather than examine precise algebraic rates, our goal is to demonstrate the OOD performance for a series of different functions motivated by applications, in scenarios where the support of the test distribution varies. As in the previous section, we do not implement the estimator described in our theoretical result (Theorem \ref{thm:deep-learning-extrap}). Instead, we opt for a fully-connected, feedforward architecture that is independent of $\bm{\omega}$ (recall Remark \ref{rem:unknown-supports}). Details of the architecture and training procedure are given in \S \ref{app:DNN}.

We consider several standard test functions from the \textit{Virtual Library of Simulation Experiments} \cite{surjanovic2013virtual}, specifically, the \textit{Physical Models} section of the \textit{Emulation/Prediction Test Problems} test set. This is a standard library for high-dimensional approximation and integration tasks, where 
each function represents a model for a certain physical process. See \cite{surjanovic2013virtual} for details. To align with the rest of the paper, we rescale the input variables to $[-1,1]^d$ and consider $\varrho = \cU([-1,1]^d)$ as the training distribution.

In Fig.\ \ref{fig:intro-examples}, we consider the \textit{wing weight} function. This is a $d = 10$ variate function, which, after rescaling to $[-1,1]^{10}$ is holomorphic inside the hyperrectangle $D_{\bm{\omega}} = \otimes^{10}_{i=1}[-\omega_i,\omega_i]$ for any $\bm{\omega} = (\omega_i)^{10}_{i=1}$ satisfying
$\bm{\omega}  < \bm{\omega}_{\max}$, where $\bm{\omega}_{\max}  = (7,6.5,4,9,2.10,3,2.6,2.43,5.25,+\infty).
$
In Fig.\ \ref{fig:intro-examples} we consider $\mu = \cU(D_{\bm{\omega}})$, where $\omega_i = 1 + \theta (\omega_{i,\max}-1)$, $\forall i \in [10]$, with the value $\omega_{\max,10} = 10$. In other words, we study extrapolation to a domain whose side lengths are a constant fraction of the maximum possible extrapolation domain. As expected, in all cases we see algebraic convergence with the rate depending on the size of $\theta$. Note that we plot the squared $L^2$-error in this experiment, in order to make a comparison with the distributional shift bound \ef{Wp-bound}. As we see from the table, the Wasserstein distance $W_2(\varrho,\mu)$ is significantly larger than this error, thus demonstrating the pessimistic nature of such bounds.

\begin{figure}[t]
    \centering
    \begin{subfigure}[b]{0.4\textwidth}
        \centering
     \includegraphics[width=\textwidth, trim= 0 0 150 0,clip] {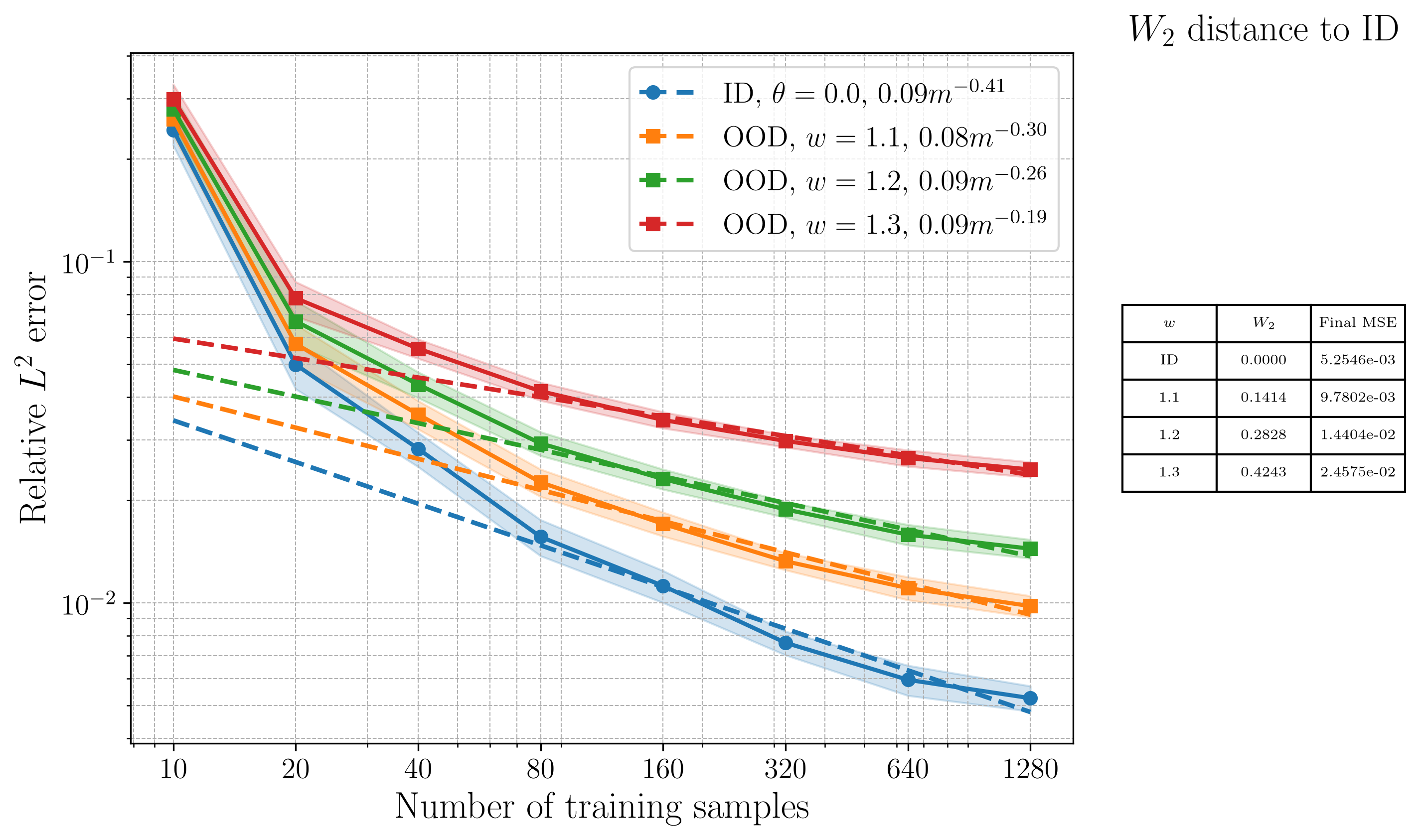}\\
        {\small
       \begin{tabular}{c|cccc}
            \hline
            $\omega$ & 1.1 & 1.2 & 1.3\\
	   $W_2(\varrho,\mu)$ & 0.141 & 0.283 & 0.424 \\
            \hline
        \end{tabular} 
        }
    \end{subfigure}
    \begin{subfigure}[b]{0.4\textwidth}
        \centering
        \includegraphics[width=\textwidth, trim= 0 0 150 0,clip]{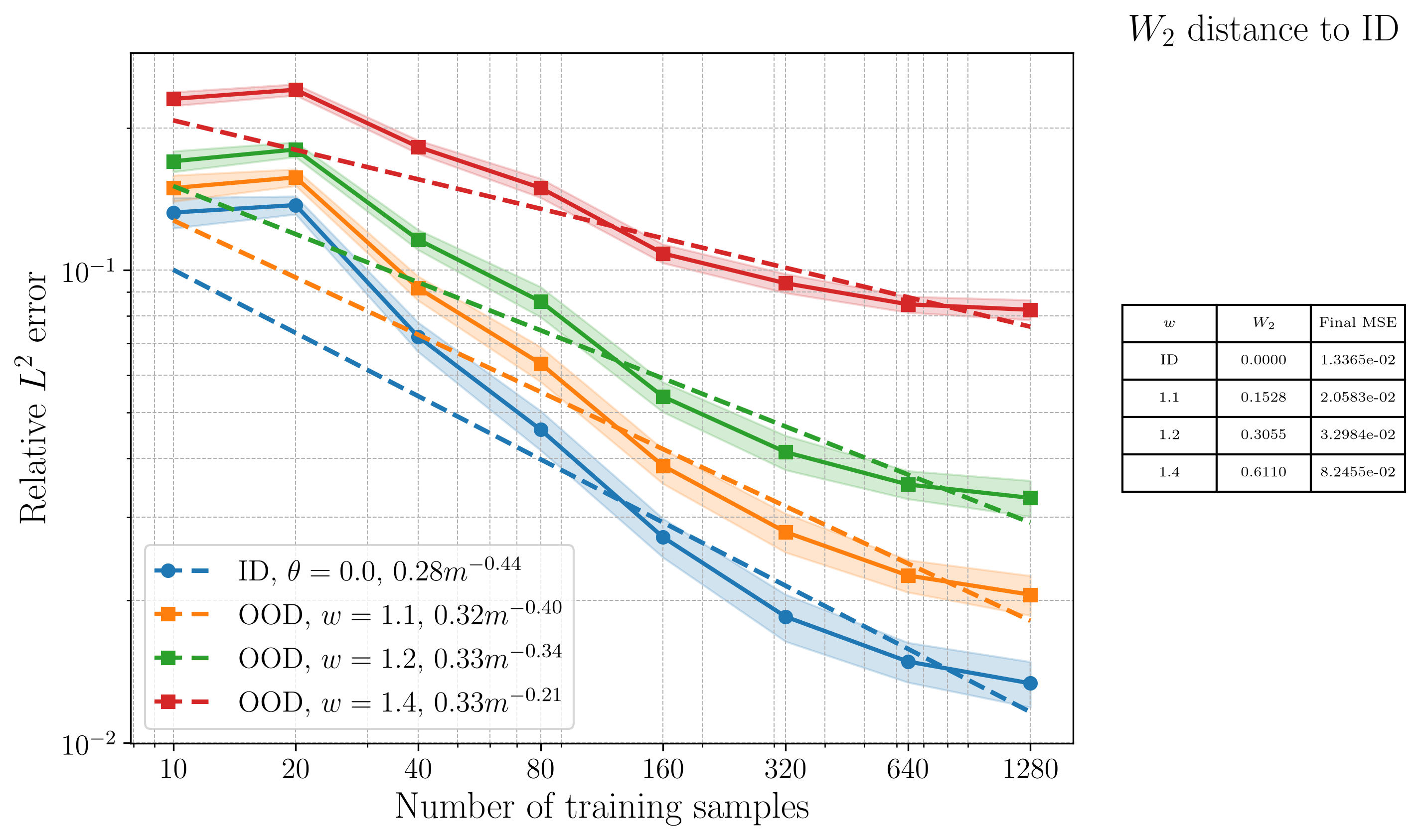}
        \\
        {\small
       \begin{tabular}{c|cccc}
            \hline
            $\omega$ & 1.1 & 1.2 & 1.4  \\
	   $W_2(\varrho,\mu)$ & 0.153 & 0.306 & 0.811 \\
            \hline
        \end{tabular} 
        }
    \end{subfigure}
    \caption{
Error versus number of training samples $m$ for approximating the circuit (left) and piston (right) functions via tanh DNNs with $L = 15$ layers and width $W = 150$.}
    \label{fig:circuit-piston}
\end{figure}

In Fig.\ \ref{fig:circuit-piston} we consider two other functions from the same library, the \textit{circuit} and \textit{piston} functions. These are $d = 6$ and $d=7$ variate functions, respectively, that are holomorphic inside the hyperrectangle
$D_{\bm{\omega}} = \otimes^{6}_{i=1}[-\omega_i,\omega_i]$ for any $\bm{\omega} = (\omega_i)^{6}_{i=1}$ satisfying $\bm{\omega} = (\omega_1,\ldots,\omega_{6}) < \bm{\omega}_{\max} : = (2,2.11,1.40,2.85,20.427,1.40)$ and $\bm{\omega} = (\omega_1,\ldots,\omega_{7}) < \bm{\omega}_{\max} : = (3,1.67,1.5,1.5,10.97,67,35)$, respectively.
We consider $\varrho = \cU([-1,1]^d)$ and $\mu = \cU([-\omega,\omega]^d)$, where $\omega \geq 1$ is a parameter we vary. For both functions, we witness algebraic convergence for all values of $\omega$. It is noticeable that the algebraic rate decreases slightly as $\omega$ increases. For the circuit function, the algebraic power decreases from $0.4$ in the $\omega = 1$ (i.e., in distribution) case to $0.19$ in the $\omega  = 1.3$ case. Similarly, for the piston function it decreases from 0.44 when $\omega  = 1$ to 0.21 when $\omega  = 1.4$. 


\rem{[Differences between Theorem \ref{thm:deep-learning-extrap} and numerical experiments]
\label{rem:dnn-differences}
The primary difference is the choice of DNN family $\cN$. As discussed in \ref{s:deep-learning}, Theorem \ref{thm:deep-learning-extrap} uses a specific handcrafted family of DNNs, where only certain values for the internal weights and biases are allowed. Conversely, in our numerical experiments, we consider fully-connected, feedforward DNNs where all weights and biases are trained and can take arbitrary real values. In particular, the resulting estimator is independent of $\bm{\omega}$, whereas the DNN family, and consequently the DNN estimator, of Theorem \ref{thm:deep-learning-extrap} depends on $\bm{\omega}$ (recall Remark \ref{rem:unknown-supports}).
}

\subsection{DNO estimators for PDE operators}\label{ss:ex-DNO}

Finally, we consider operator learning. As in the previous subsections, we do not strive to implement the DNO estimators described in Theorem \ref{thm:operator-learning-extrap}. Instead, due to their wide popularity and robust performance, we consider FNOs. Details of the architecture and training procedure are given in \S \ref{app:DNO}.

In our first example, we consider the Darcy flow PDE
\begin{equation}
\begin{split}
    -\nabla\cdot(a(x)\nabla u(x)) =f(x), \ x \in D,\qquad 
    u(x) = 0,\ x \in \partial D.
    \end{split}
    \label{eq: darcy}
\end{equation}
Here $D=(0,1)^2$, $u(x)$ is the pressure field, $f(x)$ is the forcing function and $a(x)$ is the input permeability field. We take $f(x) = 1$ and consider the operator $F$ mapping from the diffusion coefficient $a(x)$ to the solution $u(x)$, i.e., $F : a \mapsto u$.

\begin{figure}[t]
    \centering
\begin{subfigure}[b]{0.32\textwidth}
\centering
 \includegraphics[width=\textwidth] {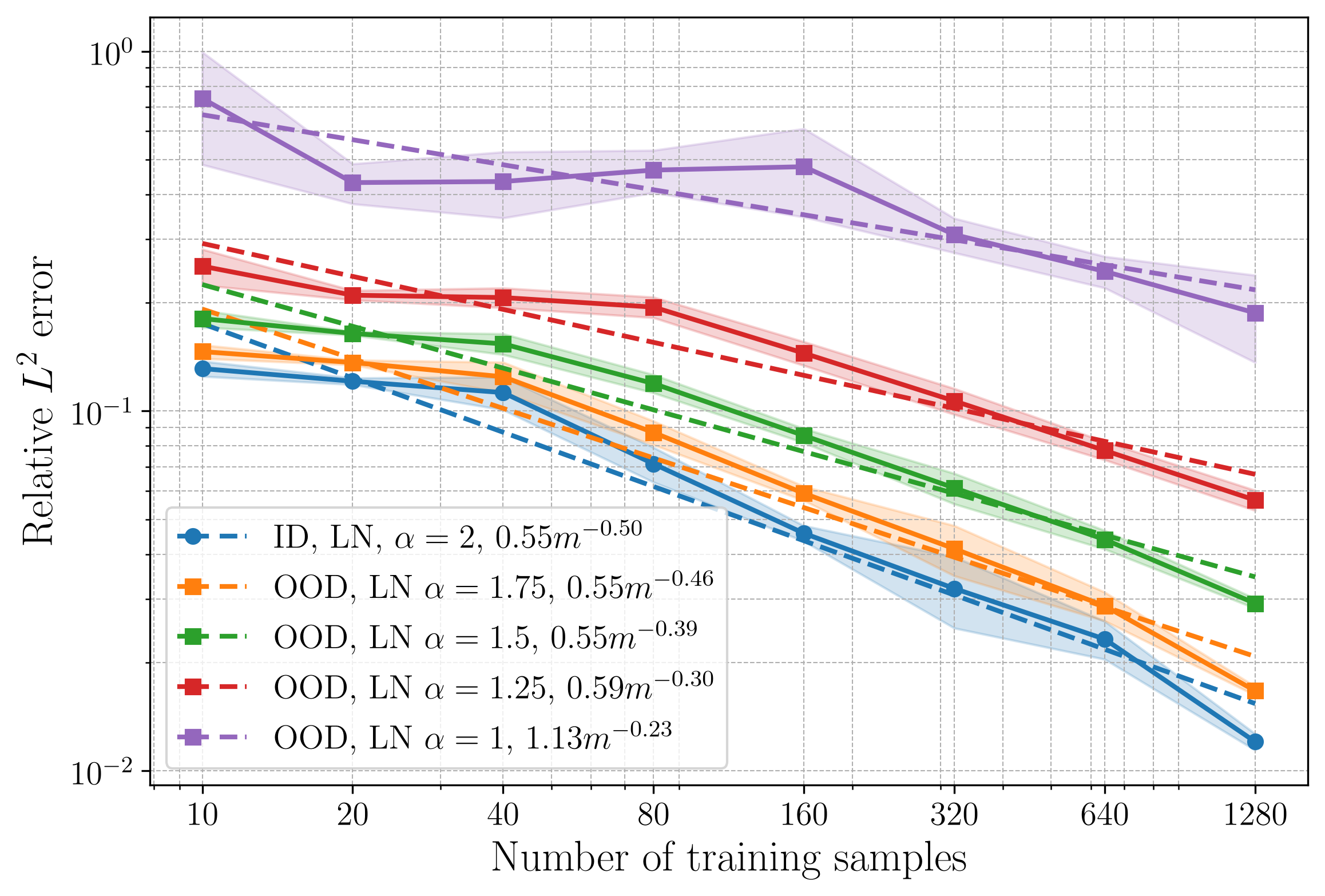}
 
    \end{subfigure}
\begin{subfigure}[b]{0.32\textwidth}
\centering
 \includegraphics[width=\textwidth] {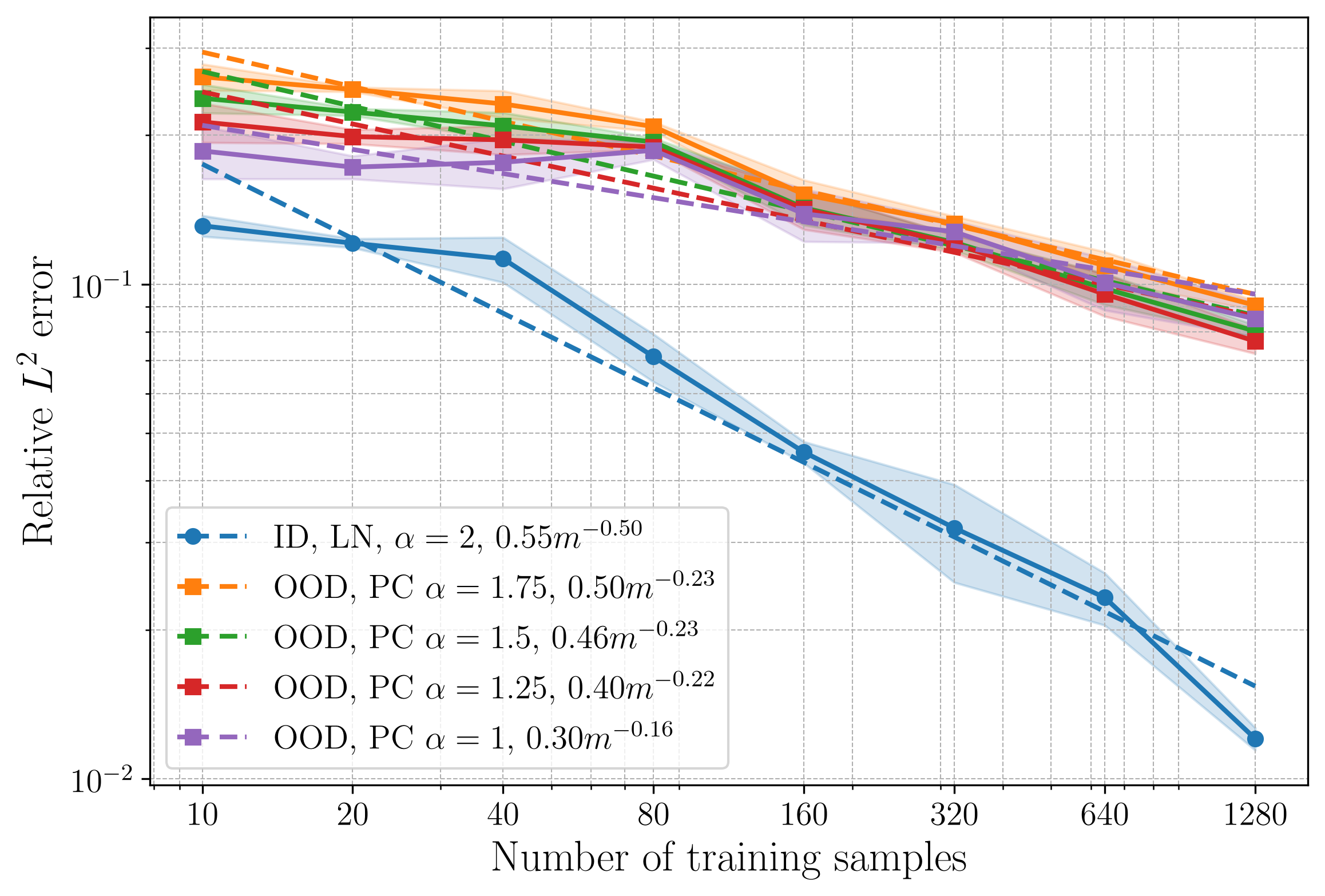}
    \end{subfigure}
\begin{subfigure}[b]{0.32\textwidth}
\centering
 \includegraphics[width=\textwidth] {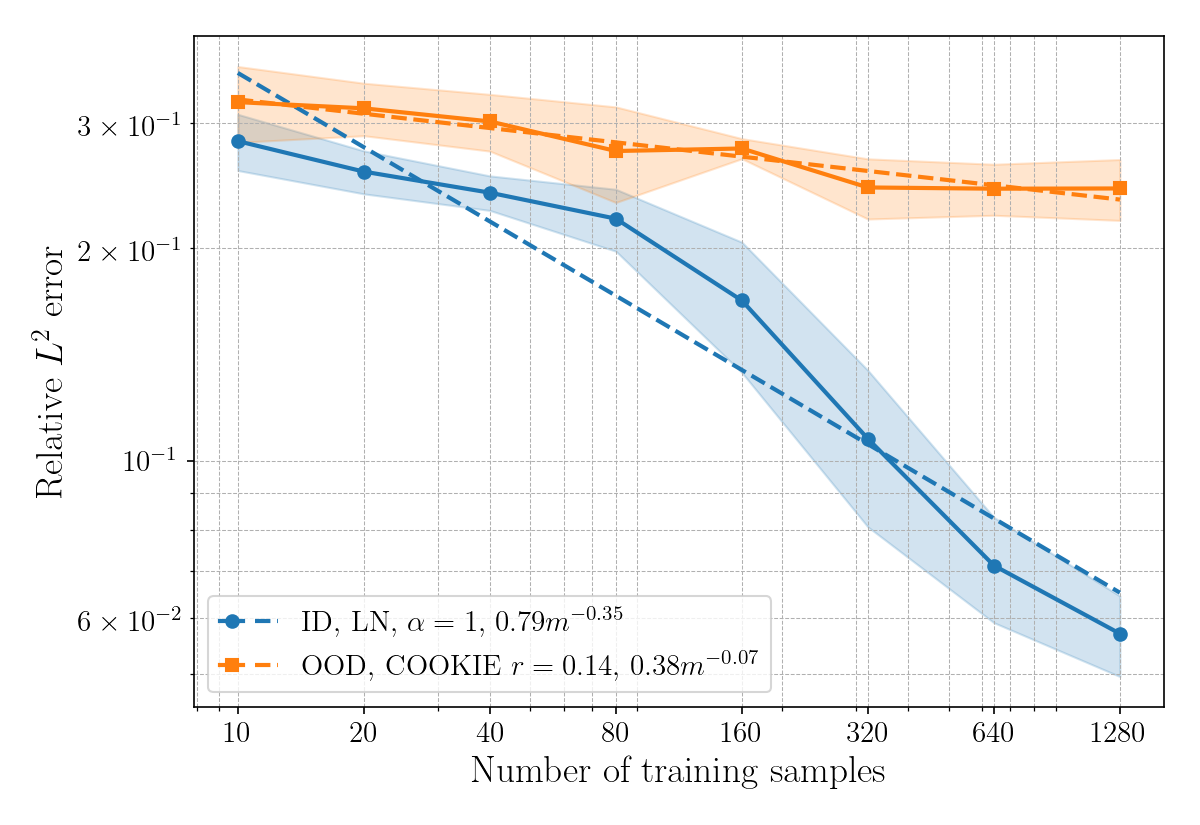}
    \end{subfigure}
    \caption{Error versus number of training samples $m$ for the Darcy flow problem
 with  $\nu =  \mu^{\mathrm{LN}}_{3,2} $ and $\mu  =  \mu^{\mathrm{LN}}_{3,\alpha} $ (left), $\mu = \mu^{\mathrm{PC}}_{3,\alpha} $ (middle); $\nu = \mu_{3,1}^{LN}$ and $\mu = \mu^{\mathrm{cookies}}$ (right).}
    \label{fig:Darcy-exp}
\end{figure}


In Fig.\ \ref{fig:Darcy-exp}, we follow standard practice and consider log-normal random fields, with corresponding probability measures
\be{
\label{mu-lognormal}
 \mu^{\mathrm{LN}}_{\tau,\alpha} =  \psi \#
    \mathcal{N}\left(
        0,(-\Delta+\tau^2I)^{-\alpha}
    \right),\qquad \text{where } \psi(x) = \exp(x) .
}
Here $\Delta$ is the Laplacian, which we equip with homogeneous Neumann boundary conditions. In Fig.\ \ref{fig:Darcy-exp} (left) we consider $\nu =  \mu^{\mathrm{LN}}_{3,2}$ for training and for testing we consider $\mu =  \mu^{\mathrm{LN}}_{3,\alpha}$ for various $\alpha < 2$. This corresponds to the setting of Example \ref{ex:opl-rougher-test}, since the eigenvalues of the covariance operator of $\mathcal{N}\left(
        0,(-\Delta+\tau^2I)^{-\alpha}
    \right)$ are $(\pi^2 (k^2_1 + k^2_2) + \tau^2)^{-\alpha}$ for $k_1,k_2 \in \bbN_0$, which means that the test functions are rougher for smaller $\alpha$. Fig.\ \ref{fig:Darcy-exp} illustrates the same qualitative behaviour predicted by our theory: the OOD generalization error decays at an algebraically in $m$, at a rate that decreases as $\alpha \rightarrow 1^{+}$.
Next we test the conclusions of Example \ref{ex:opl-arbitrarily-rough}. We once more consider $\nu = \mu^{\mathrm{LN}}_{\tau,\alpha}$, while for testing we use distributions of piecewise constant functions. We consider two examples. The first is a standard distribution in operator learning \cite{li2021fourier}, denoted $\mu^{\mathrm{PC}}_{\tau,\alpha}$, which takes the same form as \ef{mu-lognormal},
 except with $\psi(x)$ replaced by a piecewise constant function $\phi$ taking value $\phi(x) = 2$ when $x >0$ and $\phi(x) = 0.5$ otherwise. The second is the `cookies' distribution, denoted $\mu = \mu^{\mathrm{cookies}}$, which is common in parametric PDEs \cite{back2011stochastic,chkifa2015discrete}. See \ef{cookies}. In both cases, we see the behaviour expected: namely, algebraic convergence (albeit at a slow rate) even though realizations from the test distribution are very different from those seen in training. 

In our second example, we consider the incompressible Navier--Stokes equations in vorticity form on the torus $\bbT^2$:
\begin{equation}
\begin{aligned}
    \partial_t w(x,t)
    +u(x,t)\cdot\nabla w(x,t)
    &=
    \nu\Delta w(x,t)+f(x),
    && x\in \bbT^2,\quad t\in(0,T],\\
    \nabla\cdot u(x,t)&=0,
    && x\in \bbT^2,\quad t\in[0,T],\\
    w(x,0)&=w_0(x),
    && x\in \bbT^2.
\end{aligned}
\label{eq:ns-vorticity}
\end{equation}
Here, $u$ is the velocity field, $w=\nabla\times u$ is the vorticity,
$w_0$ is the initial vorticity, $\nu>0$ is the viscosity and $f$ is a
forcing function. We learn the solution operator
\(
F:
    w_0
    \longmapsto
    w|_{\mathbb{T}^2\times(0,T]}.
\)
Our setup is a standard one in the operator learning literature \cite{li2021fourier}. For training and test distributions, we use $\mu_{\tau,\alpha} =  \mathcal{N}(
0,\sigma(-\Delta+\tau^2I)^{-\alpha})$, where $\Delta$ is equipped with periodic boundary conditions.
Throughout, we use parameters $\sigma=7^{3/2}$, $\tau=7.5$, $\nu=10^{-3}$ (equivalently $\mathrm{Re}=10^3$), $T=50$ and the forcing $f(x)=0.1(\sin(2\pi(x_1+x_2))+\cos(2\pi(x_1+x_2)))$. For training, we use $\alpha=2.5$ while in testing we use various values $\alpha < 2.5$, thus producing rougher test distributions. 

Fig.\ \ref{fig:intro-examples} gives our results for the parameter $\theta : = 2.5-\alpha$. Once more, we see algebraic convergence for all values of $\alpha$, with the rate decreasing as $\alpha \rightarrow 1^{+}$. In this figure, we plot the squared $L^2$-error, to make a direct comparison with the bounds \ef{Wp-bound}. Notably, the errors are significantly smaller than the corresponding Wasserstein distances $W_2(\nu,\mu)$, thus demonstrating the inadequacy of distributional shift theory in explaining the OOD generalization performance of operator learning models, which is the main thesis of this work.

\rem{[Differences between Theorem \ref{thm:operator-learning-extrap} and the experiments]
\label{rem:dno-differences}
As in Remark \ref{rem:dnn-differences}, the architectures used in our experiments (i.e., FNOs) differ from those used in Theorem \ref{thm:operator-learning-extrap}, which are handcrafted PCA-Net-type architectures. Also, we use Gaussian or lognormal random fields for training in our numerical experiments, while Assumption \ref{ass:opl} means that in our theory the training distribution has a Karhunen--Lo\`eve expansion defined by i.i.d.\ $\cU([-1,1])$ random variables, rather than $\cN(0,1)$ random variables. As we discuss below, extending our theory to deal with the latter -- and consequently Gaussian random fields -- is an open problem.
}

\section{Conclusion and future work}\label{s:conclusion}

OOD generalization is one of the biggest challenges in SciML. 
While standard theory shows robustness of estimators to small distributional shifts, it has often been observed that estimators generalize much further beyond their training distribution. This paper provides a theoretical analysis of this phenomenon. Working with classes of holomorphic functions and operators, we derived convergence rates for OOD generalization that depend on the specific regularity and the support of the test distribution. These show successful generalization well beyond small distributional shifts and genuine extrapolation behaviour. Our results are quantifiable, distribution agnostic (i.e., depending only on the support) and, in certain cases, the rates are optimal. We also show that stable extrapolation is possible while retaining algebraic rates, which we refer to as a blessing of high dimensionality. Finally, in the context of operator learning, our results show it is possible to generalize from smooth training distributions to arbitrarily rough test distributions. Numerical results on function and operator learning tasks support our main conclusions.
There are various avenues for future work, several of which we now describe.

\textbf{From adaptation to generalization.} As noted in Remark \ref{rem:unknown-supports}, our theoretical estimators require knowledge of $\bm{\omega}$, i.e., a support estimate for the test distribution $\mu$. Thus, our results are of \textit{domain adaptation} form. In domain adaptation \cite{ben-david2006analysis,ben-david2010theory} one uses information about $\mu$ to construct ML models that generalize beyond their training set. Interestingly, our results require only very weak knowledge of $\mu$: namely, a support estimate. However, the estimators used in our numerical results are independent of $\bm{\omega}$ (see Remark \ref{rem:unknown-supports}), thus they pertain to the \textit{domain generalization} setting.
An important question is whether one can construct estimators independent of $\bm{\omega}$ that yield the same generalization bounds.. This would provide a theoretical explanation for the domain generalization behaviour observed in practice.

\textbf{Weaker conditions and sharper rates.} Our theoretical results provide a sufficient condition, \ef{b-r-cond-2-main}, linking holomorphic regularity to the admissible extrapolation domain and the convergence rate.
The experiments in \S \ref{ss:num-exp-poly} motivate investigating whether this condition can be weakened. They show convergence for domains larger than those ensured by \ef{b-r-cond-2-main}, as well as faster rates.
Note, however, that these observations concern particular functions, and do not determine worst-case convergence rates of the specific holomorphic function classes.
As discussed in Remark \ref{rem:rate-optimality}, the theoretical rates are optimal up to logarithmic factors for certain domains with constant side lengths. Establishing sharp rates for growing side lengths, as in Example \ref{ex:unbounded-side lengths}, remains open. Resolving these questions requires lower bounds that account for both the extrapolation domain and the required stability to perturbations.


\textbf{Gaussian distributions.} Our function approximation results consider the training distribution $\varrho = \cU([-1,1]^{\bbN})$  and arbitrary test distributions supported in the hyperrectangles $D_{\bm{\omega}}$. A natural extension is to tensor-product Jacobi distributions on $[-1,1]^{\bbN}$, using orthonormal Jacobi polynomials in place of Legendre polynomials. However, distributions with unbounded support present significantly more difficulties. This impacts the operator learning case as well, as in Assumption \ref{ass:opl} we impose that the random variables $\xi_i \sim \cU([-1,1])$. In particular, our analysis does not apply to the case where $\nu$ is a Gaussian measure on $\cX$, as in that case one would have $\xi_i \sim \cN(0,1)$. Gaussian measures are ubiquitous in operator learning and, for these reasons, we also used them in our numerical experiments. An important topic for future work is to generalize our analysis to this case. We note, however, that this limitation is not unique to our work -- existing statistical learning theory guarantees for operators typically assume boundedness of the training and test distributions \cite{liu2024deep,reinhardt2024statistical}.

\textbf{Practical architectures and statistical learning theory.} As mentioned in Remarks \ref{rem:dnn-differences} and \ref{rem:dno-differences}, our DNN/DNO architectures differ from those typically employed in practice. A challenge of future work is to extend our analysis to practical architectures. In particular, for operator learning, we consider a type of PCA-Net, while FNOs and DeepONets are substantially more popular.
The typical machinery to analyze generalization error with standard (as opposed to carefully handcrafted) architectures is statistical learning theory, which has been employed in \cite{liu2024deep,reinhardt2024statistical} and elsewhere to derive in-distribution generalization bounds for operator learning. Statistical learning theory considers the case of statistical noise whose variance may not be small, i.e., $e_i \sim \cN(0,\sigma^2)$ in \ef{training-data-intro}, and can be used to construct DNN/DNO estimators that achieve \textit{minimax optimal rates} for nonparametric regression. Minimax rates concern the best achievable accuracy form $m$ noisy samples uniformly over given class of functions or operators. Our work focuses on the more typical setting in scientific computing, of bounded, potentially adversarial noise that is small in norm. Thus, our results are more closely associated with \textit{optimal recovery rates}, i.e., the best achievable accuracy over a given class from $m$ noiseless samples.  See \cite{devore2025optimal} for discussion on minimax versus optimal recovery rates. A major distinction is that minimax rates typically decay no faster than the Monte Carlo rate, i.e., $\ord{m^{-1/2}}$, while optimal recovery rates can be arbitrarily fast depending on the function or operator class. This is also reflected in our results, where the algebraic order can be arbitrarily large with sufficient smoothness.
Nonetheless, it is an interesting problem for future work to use statistical learning theory to derive OOD generalization guarantees for holomorphic functions and operators in the presence of statistical noise.

\textbf{Structure beyond regularity.} As mentioned in \S \ref{ss:limitations}, holomorphy is a strong assumption which may not hold in practice. This has resulted in increasing discussion about what underlying structure may allow high-dimensional functions and operators to be efficiently learned \cite{brugiapaglia2026short,boulle2024mathematical,kovachki2024operator,subedi2026operator}. Related to this is the question of how to enforce such structure into SciML models. For instance, regularization functionals have been used to encode PDE priors such as boundary conditions \cite{li2024physics} and variational formulations \cite{goswami2023physicsinformed}. Other approaches include PDE-informed architectures \cite{boulle2022datadriven, gin2021deepgreen}. See also \cite{serrano2026test} and references therein. Much of this discussion centres on in-distribution generalization, but the same considerations apply OOD as well. Understanding what types of structures aid OOD generalization and, critically, how to enforce such structure in SciML models is a significant open problem.

\textbf{Enhancing OOD performance.} Similarly, there are many methods, largely empirical, that strive to improve OOD generalization in SciML. These include: \cite{zhu2023reliable}, which studies incorporating physics, fine-tuning with new observations and multi-fidelity data; \cite{chen2024dataefficient}, which considers in-context learning; \cite{goswami2022deep}, which considers transfer learning; and \cite{brivio2026ptpi,setinek2026simshift}, which consider domain adaptation. See also \cite{serrano2026test} and references therein. It would be interesting to understand if any of these approaches conveyed a theoretical advantage over the generalization guarantees derived in this work.

\textbf{Further types of OOD generalization.} Finally, as mentioned in \S \ref{ss:limitations}, there are different types of OOD generalization in SciML, especially in the realm of multi-operator learning and foundation models. Our work considers OOD generalization in the sense of differing training and test distributions. Extending our work to other types of OOD generalization is an interesting topic for future research.

\section*{Acknowledgments}

The authors would like to thank Chris Budd, Stefania Fresca, Anastasis Kratsios and Jakob Zech for helpful discussions.

\bibliographystyle{siamplain}
\bibliography{OODrefs}

\end{document}


\maketitle

\section{Proofs of the main results}\label{s:proofs}

In this appendix, we present the proofs of the main results from the paper. 

Overall, these proofs use, adapt and extend tools from the approximation theory of infinite-dimensional holomorphic functions \cite{cohen2010convergence,cohen2011analytic,chkifa2015breaking,adcock2022sparse,cohen2015approximation}, in combination with emulation techniques \cite{adcock2025near,de-ryck2021approximation,guhring2021approximation,li2020better,lu2021deep,yarotsky2017error,schwab2019deep,mhaskar1996neural,opschoor2022exponential,schwab2023deepa} to establish results for learning DNNs and DNOs. 
We commence in \S \ref{ss:additional-notation} with some additional notation and concepts.  \S \ref{ss:summability-bound} presents what is arguably the key technical innovation in this work, namely, new weighted summability bounds (Lemma \ref{lem:abstract-summability-1} and \ref{lem:abstract-summability-2}) and weighted $k$-term approximation error bounds (Theorems \ref{thm:weighted-summability-coeffics-1} and \ref{thm:weighted-summability-coeffics-2}) for holomorphic functions for \textit{exponentially growing} weights. As we explain, summability with respect to exponentially growing weights is directly related to extrapolation with polynomials. 

In \S \ref{ss:proofs-OOD-s-term} we use these results to prove our main result on the existence of good polynomial estimators for extrapolation, Theorem \ref{t:OOD-s-term}. Next, in \S  \ref{ss:sample-efficient-learning-proofs-II} we prove the main result on sample-efficient out-of-distribution learning via polynomials, Theorem \ref{thm:poly-optimized-least-squares}. Finally, in \S \ref{ss:deep-learning-proofs} we establish the main results for DNNs and DNNOs, Theorems \ref{thm:deep-learning-extrap} and \ref{thm:operator-learning-extrap}, respectively.

\subsection{Additional notation and concepts}\label{ss:additional-notation}

While the results in the main paper consider real-valued functions $f : D \rightarrow \bbR$, in our proofs we will establish results for Hilbert-valued functions $f : D \rightarrow \cY$, where $(\cY , \ip{\cdot}{\cdot}_{\cY})$ is a Hilbert space. This not only generalizes the results in the main paper, but is also needed in order to establish the results on operator learning considered in \S \ref{s:operator-learning}.

This necessitates some additional concepts (see, e.g., \cite{adcock2022sparse,adcock2025optimal} and references therein). First, given a set $D$ with a measure $\varrho$ and $1 \leq p \leq \infty$, we define the Lebesgue--Bochner space $L^p_{\varrho}(D;\cY)$ as the space consisting of (equivalence classes of) strongly $\varrho$-measurable functions $f: D \rightarrow \cY$ for which $\|f\|_{L^p_{\varrho}(D;\cY)}<\infty$, where
\bes{
\| f \|_{L^p_{\varrho}(D;\cY)} : = 
\begin{cases} 
\left( \int_{D} \nm{f( \bm{y} )}_{\cY}^p \D \varrho (\bm{y}) \right)^{1/p} & 1 \leq p < \infty ,
\\
\mathrm{ess} \sup_{\bm{y} \in D} \nm{f(\bm{y})}_{\cY}  & p = \infty.
\end{cases}
}
For simplicity we write $L^{p}_{\varrho}(D)$ when $\cY=\bbR$. 

Second, we require various notions relating to sequences whose entries take values in $\cY$. Let $d \in \bbN \cup \{ \infty \}$, $\Lambda \subseteq \cF$ be a multi-index set and $\bm{c} = ( c_{\bn})_{\bn \in \Lambda}$ be a sequence with $\cY$-valued entries. For $0 < p \leq \infty$, we write $\ell^p(\Lambda ; \cY)$ for the space of $\cY$-valued sequences $\bm{c} = (c_{\bn})_{\bn \in \Lambda}$ for which $\nm{\bm{c}}_{p;\cY} < \infty$, where
\eas{
\nm{\bm{c}}_{p;\cY} = \begin{cases} \left ( \sum_{\bn \in \Lambda} \nm{c_{\bn}}^p_{\cY} \right )^{\frac1p} & 0 < p < \infty ,
\\
\sup_{\bn \in \Lambda} \nm{c_{\bn}}_{\cY} & p = \infty. 
\end{cases}
}
When $\cY = \bbR$, we just write $\ell^p(\Lambda)$ and $\nm{\cdot}_p$. Note that $\nms{\cdot}_{p;\cY}$ defines a norm when $1 \leq p \leq \infty$ and quasi-norm when $0<p<1$. Now let $\bm{c} \in \ell^p(\Lambda ; \cY)$ and $s \in \bbN$ with $1 \leq s \leq | \Lambda |$. We define the $\ell^p$-norm \textit{best $s$-term approximation error} of $\bm{c}$  as
\be{
\label{sigma-s-def}
\sigma_s(\bm{c})_{p;\cY} = \inf \left \{ \nm{\bm{c} - \bm{z}}_{p;\cY} : \bm{z} \in \ell^p(\Lambda ; \cY),\ | \mathrm{supp}(\bm{z}) | \leq s \right \}.
}
Here, 
\bes{
\mathrm{supp}(\bm{z}) = \{\bn \in \Lambda : z_{\bn} \neq 0\}
}
is the \textit{support} of $\bm{z}$. Given a $\cY$-valued sequence $\bm{c}$ and $S \subseteq \Lambda$, we write $\bm{c}_{S}$ for the $\cY$-valued vector with $\bn$th entry equal to $c_{\bn}$ if $\bn \in S$ and zero otherwise. Observe that
\bes{
\sigma_{s}(\bm{c})_{p;\cY} = \inf \left \{ \nm{\bm{c} - \bm{c}_S}_{p;\cY} : S \subseteq \Lambda,\ |S| \leq s \right \}.
}
Now let $\bm{w} = (w_{\bn})_{\bn {\in \Lambda}} > \bm{0}$ be a vector of positive weights. For $1 \leq p \leq 2$, we define the weighted $\ell^p_{\bm{w}}({\Lambda}; \cY)$ {space} as the space of $\cY$-valued sequences $\bm{c} = (c_{\bn})_{\bn \in {\Lambda}}$ for which $\nm{\bm{c}}_{p,\bm{w};\cY} < \infty$, where 
\bes{
\nm{\bm{c}}_{p,\bm{w};\cY} =  \left ( \sum_{\bn \in {\Lambda}} w^{2-p}_{\bn} \nm{c_{\bn}}^p_{\cY} \right )^{\frac1p} ,\quad 1 \leq p \leq 2 .
}
Notice that $\nms{\cdot}_{p,\bm{w} ; \cY}$ coincides with the unweighted norm $\nms{\cdot}_{p;\cY}$ for any $\bm{w}$ when $p=2$, or for any $p$ when $\bm{w}=\bm{1}$.  Next, we define the weighted cardinality of an index set  $S \subseteq {\Lambda}$  as $| S |_{\bm{w}} = \sum_{\bn \in S} w^2_{\bn}$ and, for $0 < k \leq |\Lambda|_{\bm{w}}$, we define the  $\ell^p_{\bm{w}}$-norm \textit{weighted best $(k,\bm{w})$-term
 approximation error} of a $\cY$-valued vector or sequence $\bm{c} \in \ell^p_{\bm{w}}(\Lambda ; \cY)$ as
\begin{equation}
\label{weighted-k-w-term}
\sigma_{{k}}(\bm{c})_{p,\bm{w};\cY} = \inf \left \{ \nm{\bm{c} - \bm{z}}_{p,\bm{w};\cY} : \bm{z} \in \ell^p_{\bm{w}}(\Lambda ; \cY),\ | \mathrm{supp}(\bm{z}) |_{\bm{w}} \leq {k} \right \}.
\end{equation}
For later use, we also require the following result, which is a weighted version of Stechkin's inequality. See \cite[Lem.\ 3.12]{adcock2022sparse} (note the result therein is only stated for $\cY = \bbR$, but it readily extends to Hilbert-valued sequences).

\lem{[Weighted Stechkin's inequality]
\label{lem:weighted-stechkin}
Let $d \in \bbN \cup \{ \infty \}$, $0 < p \leq 2 \leq 2$, $k > 0$, $\bm{c} \in \ell^p_{\bm{w}}(\Lambda ; \cY)$ for some $\Lambda \subseteq \bbN^d_0$ and $\bm{w} = (w_{\bn})_{\bn \in \Lambda} > \bm{0}$. Then
\bes{
\sigma_k(\bm{c})_{q,\bm{w} ; \cY} \leq \nm{\bm{c}}_{p,\bm{w} ; \cY} k^{\frac1q-\frac1p}.
}
}
Finally, we also need several concepts involving lower and anchored sets. For $i \in [d]$, let $\bm{e}_i$ be the multi-index with $i$th entry equal to $1$ and all other entries equal to zero.
A set $\Lambda \subseteq \bbN^d_0$, where $d \in \bbN \cup \{ \infty \}$, is \textit{lower} if, whenever $\bn \in \Lambda$  and $\bm{n'} \leq \bn $, it also holds that $\bm{n'} \in \Lambda$. Moreover, $\Lambda$ is \textit{anchored} if it is lower and if, whenever $\bm{e}_j \in \Lambda$ for some $j \in [d]$, it also holds that $\{\bm{e}_1,\bm{e}_2, \ldots, \bm{e}_{j} \}\subseteq \Lambda$.
A scalar-valued sequence $\bm{d} = (d_{\bn})_{\bn \in \Lambda}$ is \textit{monotonically nonincreasing} if $d_{\bn} \geq d_{\bm{n'}}$ whenever 
$\bn \leq \bm{n'}$. It is \textit{anchored} if it is monotonically nonincreasing and $d_{\bm{e}_j} \leq d_{\bm{e}_i}$ whenever $i,j \in [d]$ with $i \leq j$. Now let $\bm{c} \in \ell^{\infty}({\Lambda};\cY)$. An \textit{anchored majorant} of $\bm{c}$ is any scalar-valued sequence $\bm{d} = (d_{\bn})_{\bn \in \Lambda}$ that is anchored and satisfies $d_{\bn} \geq \nm{c_{\bn}}_{\cY}$, $\forall \bn \in \Lambda$. The \textit{minimal anchored majorant} of $\bm{c}$ is the (unique) scalar-valued sequence $\tilde{\bm{c}}$ that is an anchored majorant and for which $\tilde{\bm{c}} \leq \bm{d}$ for any other anchored majorant $\bm{d}$. It has the explicit expression
\begin{equation}\label{def:min_anch}
\tilde{c}_{\bn} = 
\begin{cases}
\sup \{ \|c_{\bm{\mu}}\|_{\cY}: \bm{\mu} \geq \bn \} & \text{ if }  \bn \neq \bm{e}_j \text{ for any } j \in[d],  \\ 
\sup \{ \|c_{\bm{\mu}}\|_{\cY}: \bm{\mu} \geq \bm{e}_i \text{ for some } i \geq j \} & \text{ if } \bn =\bm{e}_j \text{ for some } j \in [d].
\end{cases}
\end{equation}
Given $0 < p \leq \infty$, we define the \textit{anchored $\ell^p$ space} $\ell^p_{\mathsf{A}}(\Lambda;\cY)$ as the space of sequences $\bm{c} \in \ell^{\infty}(\Lambda;\cY)$ for which $\tilde{\bm{c}} \in \ell^p(\Lambda)$, and define the (quasi-)norm $\|\bm{c}\|_{p,\mathsf{A};\cY} = \|\tilde{\bm{c}}\|_{p}$.
Finally, for $1 \leq s \leq |\Lambda|$ we also define the $\ell^p$-norm \textit{best $s$-term approximation error in anchored sets} by
\begin{equation}\label{def:best_anch}
\sigma_{s,\mathsf{A}}(\bm{c})_{p;\cY} = \inf \left \{ \nm{\bm{c} - \bm{z}}_{p;\cY} : \bm{z} \in \ell^p(\Lambda ; \cY),\ | \mathrm{supp}(\bm{z}) | \leq {s}, \, \mathrm{supp}(\bm{z}) \text{ anchored}\right \}.
\end{equation}

\subsection{Summability and weighted $k$-term approximation of polynomial coefficients}\label{ss:summability-bound}

In this section, we provide several results on weighted summability and weighted $k$-term approximation of polynomial coefficients of $(\bm{b},\varepsilon ; \cY)$-holomorphic functions. They are crucial steps towards establishing our main results, as they can be used to estimate both the in-distribution and out-of-distribution polynomial approximation errors.

Unweighted summability of polynomial coefficients of $(\bm{b},\varepsilon ; \cY)$-holomorphic functions is a topic at the core of infinite-dimensional polynomial approximation theory \cite{adcock2022sparse,cohen2015approximation} since, in view of Stechkin's inequality, it implies algebraic convergence of the best $s$-term approximation. The extension from unweighted to weighted summability and weighted $k$-term approximation was a crucial advance, since sets with low weighted cardinality are critical to both guaranteeing convergence in the $L^{\infty}$-norm and to ensuring that sample-efficient estimators can be constructed via least-squares or compressed sensing techniques. See \cite{rauhut2016interpolation,chkifa2018polynomial,rauhut2017compressive}, as well as \cite{adcock2022sparse,adcock2024efficient,adcock2025optimal}. However, these works consider only the in-distribution error and weights that grow algebraically fast. As we explain next, to address out-of-distribution errors, we require summability of polynomial coefficients with respect to exponentially growing weights.

Recall that $\varrho$ is the uniform probability measure on $[-1,1]^{\bbN}$.
To commence, observe first that
\be{
\label{Psi-ID-bound}
\nm{\Psi_{\bn}}_{L^{\infty}_{\varrho}(\bbR^{\bbN})} = \prod_{j \in \mathrm{supp}(\bn)} \sqrt{2 n_j + 1} =: u_{\bn},\quad \forall \bn \in \cF,
}
which follows directly from the definition of the $\Psi_{\bn}$. Next, we consider the behaviour of the Legendre polynomials on $D_{\bm{\omega}}$, where $D_{\bm{\omega}}$ is as in \ef{Dr-def}. For this, we use the integral representation
\bes{
P_n(x) = \frac{1}{\pi} \int^{\pi}_{0} (x+\sqrt{x^2-1} \cos(\theta))^n \D \theta , \quad \forall |x| > 1,\ n \in \bbN_0.
}
which follows from the contour integral representation of $P_n$ \cite[Tab.\ 18.10.1]{1964handbook}. Hence
\bes{
| P_n(x) | \leq \zeta(\omega)^{n},\quad \forall |x| \leq \omega,\ \omega > 1,\ n \in \bbN_0. 
}
We deduce that, for any probability measure $\mu$ supported in $D_{\bm{\omega}}$,
\be{
\label{Psi-OOD-bound}
\nm{\Psi_{\bn}}_{L^{\infty}_{\mu}(\bbR^{\bbN})} \leq \prod_{j \in \mathrm{supp}(\bn)} \sqrt{2 n_j + 1} \zeta(\omega_j)^{n_j} = : w_{\bn},\quad \forall \bn \in \cF.
}
The key consequence of these observations is the following. Let $f : D \rightarrow \cY$ with $\cY$-valued coefficients $\bm{c} = (c_{\bn})_{\bn \in \cF}$, $S \subset \cF$ and $f_S = \sum_{\bn \in S} c_{\bn} \Psi_{\bn}$. Then observe that the in-distribution $L^{\infty}$-norm error can be bounded by
\bes{
\nm{f - f_S}_{L^{\infty}_{\varrho}(D ; \cY)} \leq \nm{\bm{c} - \bm{c}_S}_{1,\bm{u} ; \cY}
}
where $\bm{u} = (u_{\bm{n}})_{\bm{n} \in \cF}$,
and the out-of-distribution 
$L^{\infty}$-norm error can be bounded by
\bes{
\nm{f - f_S}_{L^{\infty}_{\mu}(\bbR^{\bbN} ; \cY)} \leq \nm{\bm{c} - \bm{c}_S}_{1,\bm{w} ; \cY},
}
where $\bm{w} = (w_{\bm{n}})_{\bm{n} \in \cF}$.
Hence, the weighted Stechkin inequality (Lemma \ref{lem:weighted-stechkin}) implies algebraic decay of these errors, whenever $\bm{c} \in \ell^p_{\bm{u}}(\cF ; \cY)$ or $\bm{c} \in \ell^p_{\bm{w}}(\cF ; \cY)$, respectively. We now establish conditions under which these properties hold (note that the former is just a special case of the latter corresponding to $\bm{\omega} = \bm{1}$). The key innovation here is we show weighted summability of the polynomial coefficients of a $(\bm{b},\varepsilon;\cY)$-holomorphic function with respect to the exponentially-growing weights $\bm{w}$. Our first main result in this section is the following.

\thm{
[Weighted summability and weighted $k$-term approximation with exponentially growing weights]
\label{thm:weighted-summability-coeffics-1}
Let $\bm{\omega} \geq \bm{1}$ and suppose that $f \in \cH(\bm{b},\varepsilon ; \cY)$, where $\bm{b} \in \ell^1(\bbN)$ satisfies
\be{
\label{b-r-cond-1}
\bm{b} \odot \bm{r}_p \in \ell^p(\bbN),\quad \nm{\bm{b} \odot (\bm{r}_p-\bm{1}) }_1 < \varepsilon,\qquad \text{where } \bm{r}_p = \zeta(\bm{\omega})^{2/p-1}
}
for some $0 < p <1$. Then its coefficients $\bm{c} \in \ell^p_{\bm{w}}(\cF ; \cY)$, where $\bm{w}$ is given by \ef{Psi-OOD-bound}. Moreover, for any $k > 0$, there exists a set $S \subset \cF$ depending on $k$, $\bm{b}$, $\varepsilon$ and $\bm{\omega}$ only with $|S|_{\bm{w}} \leq k$ such that
\bes{
\nm{\bm{c} - \bm{c}_S}_{q,\bm{w} ; \cY} \leq C(\bm{b},\varepsilon,\bm{\omega},p) k^{1/q-1/p},
}
for any $q \in (p,2]$.
}

Secondly, we also establish a result on approximation of such an $f$ via anchored sets, which will be needed in order to establish the results in \S \ref{s:sample-efficient-learning} and \S \ref{s:operator-learning} of the paper. Here, we recall the definition of monotone $\ell^p$-space, $\ell^p_{\mathsf{M}}(\bbN)$ (see \S \ref{ss:main-res-intro}).

\thm{
[Weighted summability and anchored $s$-term approximation with exponentially growing weights]
\label{thm:weighted-summability-coeffics-2}
Let $\bm{\omega} \geq \bm{1}$ and suppose that $f \in \cH(\bm{b},\varepsilon ; \cY)$, where $\bm{b} \in \ell^1(\bbN)$ satisfies
\be{
\label{b-r-cond-2}
\bm{b} \odot \bm{r}_p \in \ell^p_{\mathsf{M}}(\bbN),\quad \nm{\bm{b} \odot (\bm{r}_p-\bm{1}) }_1 < \varepsilon,\qquad \text{where } \bm{r}_p = \zeta(\bm{\omega})^{2/p-1}
}
for some $0 < p <1$. Then its coefficients $\bm{c} \in \ell^p_{\bm{w}}(\cF ; \cY)$, where $\bm{w}$ is given by \ef{Psi-OOD-bound}. Moreover, for any $s \in \bbN$, there exists an anchored set $S \subset \cF$ depending on $s$, $\bm{b}$, $\varepsilon$ and $\bm{\omega}$ only with $|S| \leq s$ such that
\bes{
\nm{\bm{c} - \bm{c}_S}_{q,\bm{w} ; \cY} \leq C(\bm{b},\varepsilon,\bm{\omega},p) s^{1/q-1/p},
}
for any $q \in (p,2]$.
}

To prove these results, we require the following two abstract summability lemmas. Results of this type have a long history in the literature, first appearing in \cite{chkifa2015breaking,cohen2015approximation,cohen2010convergence} in the context of Chebyshev and Legendre polynomial coefficients of, firstly, parameter-to-solution maps of parametric PDEs and, latterly, $(\bm{b},\varepsilon ; \cY)$-holomorphic functions \cite[Thms.\ 3.28 \& 3.33]{adcock2022sparse}. These two lemmas are based directly on \cite[Lem.\ 5.5]{adcock2025optimal}, which considered abstract sequences of coefficients, but only algebraically growing weights.

\lem{
[Abstract summability]
\label{lem:abstract-summability-1}
Let $0 < p < 1$, $\varepsilon > 0$, $\xi : (1,\infty) \rightarrow [0,\infty)$ be continuous with $\xi(t)$ bounded as $t \rightarrow \infty$, $\bm{r} \geq \bm{1}$ and $\bm{b} \in \ell^1(\bbN)$ be such that
\be{
\label{b-r-cond}
\bm{b} \odot \bm{r} \in \ell^p(\bbN),\qquad \nm{\bm{b} \odot (\bm{r}-\bm{1}) }_1 < \varepsilon.
}
Let $c,\gamma > 0$ be constant and $\bm{d} = (d_{\bn})_{\bn \in \cF}$ be such that $|d_{\bm{0}}| \leq c$ and, for every $\bn \in \cF \backslash \{ \bm{0} \}$,
\bes{
|d_{\bn} | \leq \prod_{k \in \mathrm{supp}(\bn)} \xi(\rho_k) \rho^{-n_k}_{k} r^{n_k}_{k} (c n_k + 1)^{\gamma}
}
for all $\bm{\rho} = (\rho_j)_{j\in \bbN} \geq \bm{1}$ satisfiying $\{ k : \rho_k > 1 \} \supseteq \mathrm{supp}(\bn)$ and
\be{
\label{rho-admissible}
\sum^{\infty}_{j=1} \left ( \frac{\rho_j + \rho^{-1}_j}{2} - 1 \right ) b_j \leq \varepsilon.
}
Then $\bm{d} \in \ell^p(\cF)$ and $\nm{\bm{d}}_p \leq C(\bm{b},\varepsilon,\bm{r},p,c,\gamma,\xi)$. 
}

\prf{
We follow the setup of \cite[Lem.\ 5.5]{adcock2025optimal}. We split $\bbN = E \cup F$, where $E = [d]$ and $F = \bbN \backslash [d]$, where $d \in \bbN$ will be chosen later. By assumption, there exists an $0 < s < 1$ such that $\nm{\bm{b} \odot (\bm{r}-\bm{1})}_1 \leq s^2 \varepsilon$. Now let $\bm{a} = (a_j)_{j \in \bbN} > \bm{1}$ be a sequence that will be specified later that satisfies
\be{
\label{a-cond-1}
\sum^{\infty}_{j=1} (a_j - 1) b_j \leq s \varepsilon
}
Given $\bn \in \cF$, let $\tilde{\bm{\rho}} = \tilde{\bm{\rho}}(\bm{b},\varepsilon,\bn) \geq \bm{1}$ be defined by
\bes{
\frac{\tilde{\rho}_j + \tilde{\rho}^{-1}_j}{2} = \begin{cases} a_j & j \in \mathrm{supp}(\bn_E) \\ a_j + \frac{(1-s) \varepsilon n_j}{b_j \nm{\bn_F}_1 }& j \in  \mathrm{supp}(\bn_F) \\ 1 & j \notin \mathrm{supp}(\bn) \end{cases}
}
Observe that $\{ k : \rho_k > 1 \} = \mathrm{supp}(\bn)$ and $\tilde{\bm{\rho}}$ satisfies \ef{rho-admissible}. Now suppose that 
\be{
\label{a-cond-2}
a_{\min} = \inf_j \{ a_j \} > 1.
}
 and notice that
\bes{
\xi(\tilde{\rho}_j) \leq C(a_{\min},\xi),\quad \forall j \in \mathrm{supp}(\bn),
}
where the constant $C(a_{\min},\xi)$ is increasing as $a_{\min}$ decreases. This follows from the assumptions on $\xi$, which imply that it is uniformly bounded on the interval $[a_{\min},\infty)$, combined with the fact that $\tilde{\rho}_j \geq a_j > a_{\min}$ for $j \in \mathrm{supp}(\bm{n})$.
This implies that
\bes{
\xi(\tilde{\rho}_j) \leq C(a_{\min},\xi) n_j + 1,\quad \forall j \in \mathrm{supp}(\bn),
}
and therefore
\be{
\label{dn-Bn-bd}
|d_{\bn}| \leq B(\bn) : = \prod_{j \in \mathrm{supp}(\bn)}  \tilde{\rho}^{-n_j}_{j} r^{n_j}_{j} (\tilde{c} n_j + 1)^{\tilde{\gamma}},
}
where $\tilde{c} = \tilde{c}(a_{\min},\xi,c) = \max \{ C(a_{\min}),\xi ) , c \}$ is likewise increasing as $a_{\min}$ decreases and $\tilde{\gamma} = \gamma+1$. We now split the product and write $B(\bn) = B_E(\bn) B_F(\bn)$ and then get
\bes{
\nm{\bm{d}}^p_p \leq \sum_{\bn \in \cF_E} B_E(\bn)^p \cdot \sum_{\bn \in \cF_F} B_F(\bn)^p = : \Sigma_E \cdot \Sigma_F,
}
where $\cF_S = \{ \bn \in \cF : \mathrm{supp}(\bn) \subseteq S \}$ for $S = E,F$. To complete the proof, it suffices to show that $\Sigma_E , \Sigma_F < \infty$.

Consider $\Sigma_E$. We have
\bes{
\Sigma_E \leq \prod^{d}_{j=1} \left ( \sum^{\infty}_{n_j = 1}  (r_{j} /a_j)^{p n_j} (\tilde{c} n_j + 1)^{p\tilde{\gamma}} \right ).
}
This is finite, provided $\bm{a}$ satisfies
\be{
\label{a-cond-3}
a_j > r_j,\quad \forall j \in \bbN.
}
Now consider $\Sigma_F$. Using the definition of $\tilde{\rho}_j$, we have
\eas{
B_F(\bn) & \leq \prod_{j \in \mathrm{supp}(\bn_F)} \left ( \frac{b_j \nm{\bn_F}_1 r_j }{(1-s) \varepsilon n_j } \right )^{n_j} (\tilde{c} n_j + 1)^{\tilde \gamma} 
\\
& \leq \frac{\nm{\bn_F}_1 !}{\bn_F !} \left ( \frac{\E (\bm{b}\odot \bm{r})_F }{(1-s) \varepsilon} \right )^{\bn_F} \prod_{j \in \mathrm{supp}(\bn_F)} \left( ((\tilde{c}+1) n_j)^{\tilde \gamma} + 1 \right )
}
(here $\bm{n}_F ! = \prod_{j \in \mathrm{supp}(\bm{n}_F) } n_j !$ is the multi-index factorial).
The argument is now identical to the corresponding part of the proof of \cite[Lem.\ 5.5]{adcock2025optimal}, except with $\bm{b}$ replaced by $\bm{b} \odot \bm{r}$. In particular, we see that $\Sigma_F < \infty$ after taking $d$ sufficiently large, provided $\bm{b} \odot \bm{r} \in \ell^p(\bbN)$.

To summarize, we have shown $\bm{d} \in \ell^p(\bbN)$, provided there exists a sequence $\bm{a} > \bm{1}$ that satisfies \ef{a-cond-1}, \ef{a-cond-2} and \ef{a-cond-3}. Set
\be{
\label{aj-def}
a_j = 1 + \max \left \{ \frac{1+s}{2s}(r_j-1) , \frac{(1-s)s \varepsilon}{2 \nm{\bm{b}}_1} \right \}.
}
Since $0 < s < 1$, we see that \ef{a-cond-2} and \ef{a-cond-3} hold.  Also, letting $\cJ = \{ j \in \bbN : a_j = 1 + \frac{1+s}{2s}(r_j-1) \}$, we see that
\eas{
\sum^{\infty}_{j=1} (a_j-1)b_j & = \frac{1+s}{2s} \sum_{j \in \cJ} (r_j-1) b_j + \sum_{j \notin \cJ} \frac{(1-s) s \varepsilon}{2 \nm{\bm{b}}_1} b_j
\\
& \leq\frac{1+s}{2s} s^2 \varepsilon + \frac{(1-s)s}{2} \varepsilon \varepsilon 
\\
& = s \varepsilon,
}
where in the second step we used the fact that $\nm{\bm{b} \odot (\bm{r}- \bm{1}) }_1 \leq s^2 \varepsilon$. Therefore \ef{a-cond-3} holds as well. This completes the proof.
}

\lem{
[Abstract anchored summability]
\label{lem:abstract-summability-2}
Consider the setup of the previous lemma, except where \ef{b-r-cond} is replaced by 
\bes{
 \bm{b} \odot \bm{r} \in \ell^p_{\mathsf{M}}(\bbN),\qquad \nm{\bm{b} \odot (\bm{r}-\bm{1})}_1 < \varepsilon .
}
Then $\bm{d} \in \ell^p_{\mathsf{A}}(\cF)$ and $\nm{\bm{d}}_{p,\mathsf{A}} \leq C(\bm{b},\varepsilon,\bm{r},p,c,\gamma,\xi)$.
}

\prf{
We proceed in a similar manner, but with some key differences. Let $s \in (0,1)$ be such that $\nm{\bm{b} \odot (\bm{r}-\bm{1})}_1 \leq s^2 \varepsilon$. Now define $a_j$ as in \ef{aj-def} once more and notice that $a_{\min} = \min \{ a_j \} > 1$. We also have
\bes{
\sum^{\infty}_{j=1} (a_j - 1) b_j \leq s \varepsilon,
}
as before.
We now write $\bbN = E \cup F$, where $E = [d]$ and $F = \bbN \backslash [d]$ once more and assume that $d \in \bbN$ is now sufficiently large so that
\bes{
\sum_{j > d} (r_j-1) b_j  \leq \frac{(1-s) \varepsilon}{2 \beta}
}
for some $\beta > 0$ that will be chosen later.
Given $\bm{b} \in \cF$, we now let $\tilde{\bm{\rho}} = \tilde{\bm{\rho}}(\bm{b},\beta,\varepsilon,\bn) \geq \bm{1}$ be defined by
\bes{
\frac{\tilde{\rho}_j + \tilde{\rho}^{-1}_j}{2} = \begin{cases} a_j & j \in \mathrm{supp}(\bn_E) \\ a_j + \beta r_j + \frac{(1-s) \varepsilon n_j}{2b_j \nm{\bn_F}_1 }& j \in  \mathrm{supp}(\bn_F) \\ 1 & j \notin \mathrm{supp}(\bn) \end{cases}
}
Observe that $\{ k : \rho_k > 1 \} = \mathrm{supp}(\bn)$ and $\tilde{\bm{\rho}}$ satisfies \ef{rho-admissible}. Note that \ef{dn-Bn-bd} also holds for this choice of $\tilde{\bm{\rho}}$. The rest of the proof involves constructing an upper bound $\tilde{B}(\bn) \geq B(\bn)$ for which $(\tilde{B}(\bn))_{\bn \in \cF}$ is monotonically nonincreasing, anchored and $\ell^p$-summable. This will immediately imply the result.

Write $B(\bn) = B_E(\bn) B_F(\bn)$ as before, and consider $B_E(\bn)$.
Let
\bes{
\kappa = 1 + \frac{(1-s)s \varepsilon}{2 \nm{\bm{b}}_1 } > 1
}
and observe that there is a constant $D = D(\tilde{c},\tilde{\gamma},\kappa) \geq 1$ such that
\bes{
(\tilde{c} n + 1)^{\tilde{\gamma}} \leq D \left ( \frac{1+\kappa}{2} \right )^n,\quad \forall n \in \bbN_0.
}
Notice from \ef{aj-def} that
\bes{
\tilde{\rho}_j \geq a_j \geq 1 + \frac{(1-s)s \varepsilon}{2 \nm{\bm{b}}_1 } = \kappa,\quad \forall j \in \mathrm{supp}(\bm{n}_E).
}
Hence
\be{
\label{tilde-BE-def}
B_{E}(\bn) \leq D^d \prod_{j \in \mathrm{supp}(\bn_E) } \left (\frac{1+\kappa}{2 \kappa} \right )^{n_j} = : \tilde{B}_E(\bn).
}
Now consider $B_F(\bn)$. Let $\bm{t} = \widetilde{\bm{b} \odot \bm{r}}$ be a minimal monotone majorant of $\bm{b} \odot \bm{r}$. Then
\eas{
B_F (\bn) & \leq \prod_{j \in \mathrm{supp}(\bn_F)} (\tilde{c} n_j + 1)^{\tilde{\gamma}} r^{n_j}_{j} \left ( \beta r_j + \frac{(1-s) \varepsilon n_j}{2 b_j \nm{\bn_F}_1 } \right )^{-n_j}  
\\
& \leq \prod_{j \in \mathrm{supp}(\bn_F)} (\tilde{c} n_j + 1)^{\tilde{\gamma}}  \left ( \beta  + \frac{(1-s) \varepsilon n_j}{2 t_j \nm{\bn_F}_1 } \right )^{-n_j} = : \tilde{B}_F(\bn)
}
Hence, we define $\tilde{B}(\bn) = \tilde{B}_E(\bn) \tilde{B}_F(\bn)$.

We now show that $\tilde{B}(\bn)$ is monotonically nonincreasing. We do this by showing it separately for $\tilde{B}_E(\bn)$ and $\tilde{B}_F(\bn)$. For the former, we see this readily holds since $\kappa > 1$ and therefore $\frac{1+\kappa}{2 \kappa} < 1$.
Now consider $\tilde{B}_F(\bn)$. To show monotonicity, we show that for every $\bn \in \cF$ and every $i > d$, we have $\tilde{B}_F(\bn + \bm{e}_i) \leq \tilde{B}_F(\bn)$. Assume first that $n_i \neq 0$. Then, arguing as in the proof of \cite[Lem.\ 5.2]{adcock2025optimal}, we have
\bes{
\frac{\tilde{B}_F(\bn + \bm{e}_i)}{\tilde{B}_F(\bn)} \leq \frac{ 2^{\tilde{\gamma}} \E}{\beta}.
}
Now suppose that $n_i = 0$. Then 
\bes{
\frac{\tilde{B}_F(\bn + \bm{e}_i)}{\tilde{B}_F(\bn)} \leq \frac{(\tilde{c} + 1)^{\tilde{\gamma}} \E }{\beta }.
}
Therefore, it suffices to choose $\beta$ satisfying
\be{
\label{beta-ineq-1}
\beta \geq \max \left \{ 2^{\tilde{\gamma}} \E , (\tilde{c}+1)^{\tilde{\gamma}} \E \right \}.
}
For any such $\beta$, we deduce that $\tilde{B}_F(\bn)$, and therefore $\tilde{B}(\bn)$, is monotonically nonincreasing.

We now show that $\tilde{B}(\bn)$ is anchored, i.e., $\tilde{B}(\bm{e}_j) \leq \tilde{B}(\bm{e}_i)$ whenever $i \leq j$. Observe that
\bes{
\tilde{B}(\bm{e}_j) = \begin{cases} 
D^d \left ( \frac{1+\kappa}{2 \kappa} \right ) & j \in E \\
(\tilde{c} + 1)^{\tilde \gamma}  \left ( \beta + \frac{(1-s) \varepsilon}{2 t_j} \right )^{-1} & j \in F 
\end{cases}
}
Suppose first that $i,j \in E$. Then $\tilde{B}(\bm{e}_j) = \tilde{B}(\bm{e}_i)$, as required.
Now suppose that $i ,j \in F$. Then we require
\bes{
\left ( \beta + \frac{(1-s) \varepsilon}{2 t_j} \right )^{-1} \leq  \left ( \beta + \frac{(1-s) \varepsilon}{2 t_i} \right )^{-1}.
}
This holds, since $\bm{t}$ is a monotonically nonincreasing by assumption. Finally, if $i \in E$, $j \in F$ then we require
\bes{
(\tilde{c} + 1)^{\tilde \gamma}  \left ( \beta + \frac{(1-s) \varepsilon}{2 t_j} \right )^{-1} \leq D^d \left ( \frac{1+\kappa}{2 \kappa} \right ).
}
We see that this holds, provided
\bes{
\beta \geq \frac{(\tilde{c}+1)^{\tilde \gamma} }{D^d\left ( \frac{1+\kappa}{2 \kappa} \right )}.
}
Recall that $D \geq 1$ by assumption. Hence this condition is implied by
\bes{
\beta \geq (\tilde{c}+1)^{\tilde \gamma} \left (\frac{2 \kappa}{1+\kappa} \right ).
}
Combining this with \ef{beta-ineq-1}, we now pick
\bes{
\beta = \max \left \{ 2^{\tilde{\gamma}} \E , (\tilde{c} + 1)^{\tilde \gamma} \E , (\tilde{c}+1)^{\tilde \gamma} \left (\frac{2 \kappa}{1+\kappa} \right ) \right \}
}
This concludes the third case, and therefore we have shown that $\tilde{B}(\bn)$ is anchored.

It remains to show that $\bm{d} \in \ell^p_{\mathsf{A}}(\cF)$. By construction, we have
\bes{
\nm{\bm{d}}^p_{p,\mathsf{A}} \leq \sum_{\bn \in \cF} \tilde{B}(\bn)^{p} = \sum_{\bn \in \cF_E} \tilde{B}_E(\bn)^p \cdot \sum_{\bn \in \cF_F} \tilde{B}_F(\bn)^p = : \tilde{\Sigma}_E \cdot \tilde{\Sigma}_F.
}
Using \ef{tilde-BE-def}, we see that
\bes{
\tilde{\Sigma}_E = D^{dp} \left ( \sum^{\infty}_{n=0} \left(\frac{\kappa + 1}{2\kappa}\right )^{n} \right )^p < \infty,
}
since $\kappa > 1$.
Now consider $\tilde{\Sigma}_F$. As in the previous proof, we have
\bes{
\tilde{B}_F(\bn) \leq \prod_{j \in \mathrm{supp}(\bn_F)} \left ( \frac{2 t_j \nm{\bn_F}_1 }{(1-s) \varepsilon n_j} \right )^{n_j} (\tilde{c} n_j + 1)^{\tilde{\gamma}}.
}
We now argue in the same way as before, except with $\bm{b} \odot \bm{r}$ replaced by the sequence $2 \bm{t}$, using the assumption that $\bm{t} \in \ell^p(\bbN)$.
}

We are now ready to establish Theorems \ref{thm:weighted-summability-coeffics-1} and \ref{thm:weighted-summability-coeffics-2}.

\prf{
[Proof of Theorem \ref{thm:weighted-summability-coeffics-1}]
We use Lemma \ref{lem:abstract-summability-1}. Let $\bm{n} \in \cF$. Using \cite[Lem.\ 5.4]{adcock2025optimal}, we have
\be{
\label{Leg-coeffics-bd}
\nm{c_{\bn}}_{\cY} \leq  \prod_{k \in \mathrm{supp}(\bn) } \frac{\rho^{-n_k+1}_k}{(\rho_k-1)^2} (n_k+1),
}
for any $\bm{\rho}$ satisfying \ef{rho-admissible} and $\{ k : \rho_k > 1 \} \supseteq \mathrm{supp}(\bn)$, since $f\in \cH(\bm{b},\varepsilon ; \cY)$ by assumption. Now write
\bes{
\nm{\bm{c}}_{p,\bm{w} ; \cY} = \nm{\tilde{\bm{c}}}_{p ; \cY},
}
where $\tilde{\bm{c}} = (\tilde{c}_{\bn})_{\bn \in \cF}$ is given by $\tilde{c}_{\bn} = c_{\bn} w^{2/p-1}_{\bn}$. Using \ef{Leg-coeffics-bd} and \ef{Psi-OOD-bound}, we have
\bes{
\nm{\tilde{c}_{\bn}}_{\cY} \leq \nm{f}_{L^{\infty}(\cE_{\bm{\rho}} )} \prod_{k \in \mathrm{supp}(\bn) } \frac{\rho^{-n_k+1}_k}{(\rho_k-1)^2} (n_k+1) \prod_{k \in \mathrm{supp}(\bn)} (2n_k+1)^{1/p-1/2} \zeta(\omega_k)^{(2/p-1)n_k}
}
Therefore
\bes{
\nm{ \tilde{c}_{\bn} }_{\cY} \leq \nm{f}_{L^{\infty}(\cE_{\bm{\rho}} )} d_{\bm{n}},
}
where
\bes{
d_{\bn} = \prod_{k \in \mathrm{supp}(\bn) } \xi(\rho_k) \rho^{-n_k}_{k} r^{n_k}_{p,k} (c n_k + 1)^{\gamma}
}
and
\bes{
\xi(\rho) = \frac{\rho}{(\rho-1)^2},\quad r_{p,k} = \zeta(\omega_k)^{2/p-1},\quad c = 2,\quad \gamma = 1/p+1/2.
}
We now apply Lemma \ref{lem:abstract-summability-1} with $\bm{r}$ replaced by $\bm{r}_p = (r_{p,k})_{k \in \bbN}$ to see that $\bm{d} \in \ell^p(\cF)$ with 
\bes{
\nm{\bm{d}}_p \leq C(\bm{b} , \varepsilon , \bm{r}_p , p ,c , \gamma , \xi ) = C(\bm{b},\varepsilon,\bm{\omega} , p ),
}
since $\bm{b} \odot \bm{r}_p \in \ell^p(\bbN)$ and $\nm{\bm{b} \odot (\bm{r}_p-\bm{1}) }_1 < \varepsilon$ by assumption. This gives the first result.

For the second result, we use Lemma \ref{lem:weighted-stechkin}. This implies that $\sigma_k(\bm{c})_{q,\bm{w}} \leq \nm{\bm{c}}_{p,\bm{w} ; \cY} k^{1/q-1/p}$ for any $q \in (p,2]$. Moreover, inspecting its proof, we see that the upper bound is achieved by $\nm{\bm{c} - \bm{c}_S}_{q,\bm{w} ; \cY}$, where the set $S$ satisfies $|S|_{\bm{w}} \leq k$ and is chosen independently of $q$ (and $p$), and depending only on $\bm{c}$ and $\bm{w}$. 
}

\prf{
[Proof of Theorem \ref{thm:weighted-summability-coeffics-2}]
Let $\bm{d}$ be the sequence defined in the previous proof. We now appeal to Lemma \ref{lem:abstract-summability-2} with $\bm{r}$ replaced by $\bm{r}_p$ to see that $\bm{d} \in \ell^p_{\mathsf{A}}(\cF)$, since $\bm{b} \odot \bm{r}_p \in \ell^p_{\mathsf{M}}(\bbN)$ and $\nm{\bm{b} \odot (\bm{r}_p-\bm{1}) }_1 < \varepsilon$ by assumption.
The result now follows.
}

\subsection{Proof of Theorem \ref{t:OOD-s-term}}\label{ss:proofs-OOD-s-term}

We shall prove the following generalization of Theorem \ref{t:OOD-s-term}, which allows for Hilbert-valued functions. For this theorem, we define the Hilbert-valued version of the distribution shift constant \ef{dist-shift-const} for model classes of $\cY$-valued functions $P : D \rightarrow \cY$ simply by replacing the Lebesgue norms in \ef{dist-shift-const} by Lebesgue-Bochner norms. Given a set $S \subset \cF$, we also define the corresponding subspace of $\cY$-valued polynomials by
\bes{
\cP_{\cS ; \cY} = \left \{ \sum_{\bn \in S} c_{\bn} \Psi_{\bn} : c_{\bn} \in \cY,\ S \in \cS \right \}.
}

\thm{
\label{t:OOD-s-term-hilbert}
Let $\varrho$ be the uniform probability measure on $[-1,1]^{\bbN}$, $\bm{\omega} \geq \bm{1}$ and $D_{\br}$ be as in \ef{Dr-def}. Suppose that $f \in \cH(\bm{b},\varepsilon ; \cY)$, where $\bm{b}$ satisfies \ef{b-r-cond-main} for some $0 < p <1$. Then, for every $k > 0$, there exists a set $S$ depending on $\bm{b}$, $\varepsilon$, $\bm{\omega}$ and $k$ with $|S| \leq k$ such that
\be{
\label{dist-shift-const-bound}
\Delta(\cP_{S ; \cY} ; \varrho , \mu) \leq \sqrt{k}
}
and a $g \in \cP_{S ; \cY}$ such that
\bes{
\Delta(\cP_{S ; \cY} ; \varrho , \mu) \nm{f - g}_{L^2_{\varrho}(\bbR^{\bbN} ; \cY)} + \nm{f - g}_{L^{\infty}_{\mu}(\bbR^{\bbN} ; \cY)} \leq C(\bm{b},\varepsilon,\bm{\omega},p) k^{1-1/p}
}
for any probability measure $\mu$ supported in $D_{\bm{\omega}}$.
}

\prf{[Proof of Theorem \ref{t:OOD-s-term-hilbert}]
Define the weights $\bm{u} = (u_{\bn})_{\bn \in \cF}$ and  $\bm{w} = (w_{\bn})_{\bn \in \cF}$, where $u_{\bn}$ and $w_{\bn}$ are as in \ef{Psi-ID-bound} and \ef{Psi-OOD-bound}, respectively.
Let $S \subset \cF$, $|S| \leq s$, be an index set (we will choose $S$ later) and consider first the constant $\Delta(\cP_{S ; \cY} ; \varrho , \mu)$. For any $P = \sum_{\bn \in S} c_{\bn} \Psi_{\bn} \in \cP_{S;\cY}$, the triangle inequality, \ef{Psi-OOD-bound}, the Cauchy--Schwarz inequality and Parseval's identity give that
\bes{
\nm{P}_{L^{\infty}_{\mu}(\bbR^{\bbN})} \leq \sum_{\bn \in S} |c_{\bn} | w_{\bn} \leq \sqrt{|S|_{\bm{w}} } \nm{P}_{L^{2}_{\varrho}(\bbR^{\bbN})} .
}
Since $\cP_{S ; \cY}$ is a linear space and $P \in \cP_{S ; \cY}$ was arbitrary, it follows that
\bes{
\Delta(\cP_{S ; \cY} ; \varrho , \mu) \leq  \sqrt{|S|_{\bm{w}} } .
}
Now let $P = f_S$ be as in \ef{f-S-def}. Then, by Parseval's identity and the definition of the $\ell^2_{\bm{w}}$-norm,
\bes{
\nm{f - P}_{L^{2}_{\varrho}(\bbR^{\bbN} ; \cY)} = \nm{\bm{c} - \bm{c}_S }_{2 ; \cY} =  \nm{\bm{c} - \bm{c}_S}_{2,\bm{w} ; \cY},
}
and by the triangle inequality and \ef{Psi-OOD-bound},
\bes{
\nm{f - P}_{L^{\infty}_{\varrho}(\bbR^{\bbN} ; \cY)} \leq \sum_{\bn \notin S} w_{\bn} \nm{c_{\bn} }_{\cY} = \nm{\bm{c} - \bm{c}_S}_{1,\bm{w} ; \cY} .
}
Combining with \ef{dist-shift-const-bound}, we deduce that
\bes{
\Delta(\cP_{S ; \cY} ; \varrho , \mu)\nm{f - P}_{L^2_{\varrho}(\bbR^{\bbN} ; \cY)} + \nm{f - P}_{L^{\infty}_{\mu}(\bbR^{\bbN} ; \cY)} \leq \sqrt{|S|_{\bm{w}} } \nm{\bm{c} - \bm{c}_S}_{2,\bm{w} ; \cY} +  \nm{\bm{c} - \bm{c}_S}_{1,\bm{w} ; \cY} .
}
Now, the assumption \ef{b-r-cond-main} and Theorem \ref{thm:weighted-summability-coeffics-1} imply that there is a set $S$ with $|S|_{\bm{w}} \leq k$ such that $\nm{\bm{c} - \bm{c}}_{q,\bm{w} ; \cY} \leq C(\bm{b},\varepsilon,\bm{\omega},p) k^{1/q-1/p}$ for any $q \in (p,2]$. We now substitute this into the previous expression using $q = 1,2$ to get
\be{
\label{OOD-error-bound-k}
\Delta(\cP_{S ; \cY} ; \varrho , \mu)\nm{f - P}_{L^2_{\varrho}(\bbR^{\bbN} ; \cY)} + \nm{f - P}_{L^{\infty}_{\mu}(\bbR^{\bbN} ; \cY)} \leq \nm{\bm{c}}_{p,\bm{w} ; \cY} k^{1-1/p},\quad \text{where } P = f_S.
}
The weights $\bm{w} \geq \bm{1}$ and therefore $|S| \leq |S|_{\bm{w}} \leq k$.  The result now follows.
}

\subsection{Proof of Theorem \ref{thm:poly-optimized-least-squares}}\label{ss:sample-efficient-learning-proofs-II}

We now prove a generalization of Theorem \ref{thm:poly-optimized-least-squares}.
The arguments for this proof are inspired by existing compressed sensing approaches to polynomial approximation of high-dimensional functions \cite{rauhut2016interpolation,rauhut2017compressive,adcock2022sparse,adcock2024efficient,adcock2025optimal}, although with substantial differences in order to derive out-of-distribution generalization bounds. Furthermore, we also extend the setting of Theorem \ref{thm:poly-optimized-least-squares} to consider the approximation of Hilbert-valued functions, as opposed to scalar-valued functions. This will be of use later when we consider operator learning.

Given $k > 0$ (its precise value will be chosen in the proof), let
\be{
\label{Lambda-def}
\Lambda = \left \{ \bn = (n_j)^{\infty}_{j=1} \in \cF : \prod^{\infty}_{j=1} (n_j + 1) \leq \lceil k \rceil,\ n_j = 0,\ \forall j > \lceil k \rceil \right \}.
}
and define the set
\be{
\label{S-def}
\cS = \left\{ S \subseteq \Lambda : | S |_{\bm{w}} \leq k \right \}.
}
Another aspect we deal with in this section is the case where $\cY$ is discretized to some subspace $\hat{\cY}$, which we consider a Hilbert space with the same norm. This is relevant for computations -- as in practice, one can never work over $\cY$ when it is infinite dimensional -- and relevant later in the context of operator learning. To this end, we now define 
Let $\cP_{\cS ; \hat{\cY}} = \bigcup_{S \in \cS} \cP_{S ; \hat{\cY}}$. Now consider a Hilbert-valued function $f : D \rightarrow \cY$ with training data
\be{
\label{training-data-hilbert}
(\bm{x}_i , y_i),\ i = 1,\ldots,m,\quad \text{where } y_i = f(\bm{x}_i) + e_i \in \cY
}
and $e_i \in \cY$ represents measurement noise. We now consider the Hilbert-valued nonlinear least-squares fit
\be{
\label{poly-least-squares-hilbert}
\min_{P \in \cP_{\cS ; \hat{\cY}}} \frac1m \sum^{m}_{i=1} \nm{ y_i - P(\bm{x}_i) }^2_{\cY} . 
}
The main objective of this section is to establish the following extension of Theorem \ref{thm:poly-optimized-least-squares} (the latter corresponds to the case $\cY = \hat{\cY} = \bbR$). For this and subsequent results, we require several additional concepts. First, we let $\hat{\cB}_{\cY} : \cY \rightarrow \hat{\cY}$ be the orthogonal projection onto $\hat{\cY}$. Second, given an optimization problem $\min_t g(t)$ and scalars $\sigma \geq 1$ and $\tau \geq 0$, we say that $\hat{t}$ is a \textit{$(\sigma,\tau)$-minimizer} if $g(\hat{t}) \leq \sigma^2 \min_t g(t) + \tau^2$.

\thm{
\label{thm:poly-optimized-least-squares-hilbert}
Let $k \geq 1$, $\tau > 0$, $\sigma \geq 1$, $0 < \epsilon < 1$, $\bm{\omega} \geq 1$ and $D_{\bm{\omega}}$ be as in \ef{Dr-def}. Let $f \in \cH(\bm{b},\varepsilon ; \cY)$, where $\bm{b}$ satisfies \ef{b-r-cond-2-main} for some $0 < p <1$, draw $\bm{x}_1,\ldots,\bm{x}_m \sim_{\mathrm{i.i.d}} \varrho$, where $m$ satisfies
\bes{
m \geq c \cdot k \cdot \left ( \log^4(2k) + \log(1/\epsilon) \right )
}
for some universal constant $c > 0$. Then, with probability at least $1-\epsilon$, any $(\sigma,\tau)$-minimizer $\hat{f}$ of \ef{poly-least-squares-hilbert} satisfies
\eas{
\nm{f - \hat{f}}_{L^{\infty}_{\mu}(\bbR^{\bbN} ; \cY)} \lesssim &~ \frac{\sigma C(\bm{b},\varepsilon,\bm{\omega},p)}{\sqrt{\epsilon}} k^{1-1/p}
 + \frac{\sigma \sqrt{k}}{\sqrt{\epsilon}} \nm{f - \hat{\cB}_{\cY} \circ f}_{L^{2}_{\varrho}(\bbR^{\bbN};\cY)}
 + \frac{\sigma \sqrt{k} \nm{\bm{e}}_{2;\cY}}{\sqrt{m}} 
 + \sqrt{k} \tau
,
} 
where $\bm{e} = (e_i)^{m}_{i=1} \in \cY^m$, for any probability measure $\mu$ supported in $D_{\bm{\omega}}$.
}

To prove this result, we require several lemmas.

\lem{
\label{lem:poly-least-squares-alpha-bound}
Let $\bm{\omega} \geq 1$, $D_{\bm{\omega}}$ be as in \ef{Dr-def}, $f \in C(D_{\bm{\omega}}; \cY)$, $\bm{x}_1,\ldots,\bm{x}_m \in D$ and define
\be{
\label{alpha-def}
\alpha = \inf \left \{ \frac{\sqrt{\frac1m \sum^{m}_{i=1} \nm{P_1(\bm{x}_i) - P_2(\bm{x}_i) }^2_{\cY} }}{\nm{P_1 - P_2}_{L^2_{\varrho}(\bbR^{\bbN} ; \cY)}} : P_1,P_2 \in \cP_{\cS ; \hat{\cY}},\ P_1 \neq P_2 \right \}.
}
Suppose that $\alpha > 0$ and let $\hat{f}$ be a $(\sigma,\tau)$-minimizer of \ef{poly-least-squares-hilbert}. Then
\eas{
\nm{f - \hat{f}}_{L^{\infty}_{\mu}(\bbR^{\bbN} ; \cY)} \leq ~ \sqrt{2 k}  \Bigg [ &  \nm{f - \hat{\cB}_{\cY} \circ f}_{L^2_{\varrho}(\bbR^{\bbN} ; \cY)} + \alpha^{-1}(\sigma+1) \sqrt{\frac1m \sum^{m}_{i=1}\nm{f(\bm{x}_i) - \hat{\cB}_{\cY} \circ f(\bm{x}_i)}^2_{\cY} } 
\\
& + 2 \nm{ f - P}_{L^2_{\varrho}(\bbR^{\bbN} ; \cY)}  + \alpha^{-1}(\sigma+1) \sqrt{\frac{1}{m} \sum^{m}_{i=1} \nm{ f(\bm{x}_i) - P(\bm{x}_i) }^2_{\cY} } 
\\
& + \alpha^{-1}(\sigma+1) \frac{\nm{\bm{e}}_{2;\cY}}{\sqrt{m}} + \alpha^{-1} \tau \Bigg ] +  \nm{f - P}_{L^{\infty}_{\mu}(\bbR^{\bbN}; \cY)} .
}
for all $P \in \cP_{S;\cY}$ and any probability measure $\mu$ supported in $D_{\bm{\omega}}$.
}

Notice that this error bound is with respect to polynomials $P \in \cP_{\cS ; \cY}$, and not $P \in \cP_{\cS ; \hat{\cY}}$, which is the space over which the minimization problem \ef{poly-least-squares-hilbert} is formulated. This is critical in establishing Theorem \ref{thm:poly-optimized-least-squares-hilbert}.

\prf{
 Let $P \in \cP_{\cS ; \cY}$ be arbitrary and observe that $\hat{\cB}_{\cY} \circ P \in \cP_{\cS ; \hat{\cY}}$. Then
\eas{
\nm{f - \hat{f}}_{L^2_{\varrho}(\bbR^{\bbN} ; \cY)} \leq \nm{f - \hat{\cB}_{\cY} \circ f}_{L^2_{\varrho}(\bbR^{\bbN} ; \cY)} + \nm{\hat{\cB}_{\cY} \circ f - \hat{\cB}_{\cY} \circ P}_{L^2_{\varrho}(\bbR^{\bbN} ; \cY)} + \nm{\hat{\cB}_{\cY} \circ P - \hat{f}}_{L^2_{\varrho}(\bbR^{\bbN} ; \cY)}
}
For the middle term, we use the fact that $\hat{\cB}_{\cY}$ is an orthogonal projection, and for the third term, we use the fact that $\hat{\cB}_{\cY} \circ P , \hat{f} \in \cP_{\cS ; \cY}$ to get
\be{
\label{false-creek-1}
\begin{split}
\nm{f - \hat{f}}_{L^2_{\varrho}(\bbR^{\bbN} ; \cY)} \leq &~ \nm{f - \hat{\cB}_{\cY} \circ f}_{L^2_{\varrho}(\bbR^{\bbN} ; \cY)} + \nm{ f - P}_{L^2_{\varrho}(\bbR^{\bbN} ; \cY)} 
\\
& + \alpha^{-1} \sqrt{\frac{1}{m} \sum^{m}_{i=1} \nm{\hat{\cB}_{\cY} \circ P(\bm{x}_i) - \hat{f}(\bm{x}_i) }^2_{\cY} }.
\end{split}
}
Now write the training labels $y_i$ in \ef{training-data-hilbert} as
\bes{
y_i = \hat{\cB}_{\cY} \circ f(\bm{x}_i) + e'_i,\quad \text{where }e'_i = f(\bm{x}_i) - \hat{\cB}_{\cY} \circ f(\bm{x}_i) + e_i.
}
Then for the third term above, the triangle inequality and the fact that $\hat{f}$ is an approximate minimizer give
\eas{
&\sqrt{\frac{1}{m} \sum^{m}_{i=1} \nm{\hat{\cB}_{\cY} \circ P(\bm{x}_i) - \hat{f}(\bm{x}_i) }^2_{\cY} } 
\\
& \leq \sqrt{\frac{1}{m} \sum^{m}_{i=1} \nm{y_i - \hat{f}(\bm{x}_i) }^2_{\cY} } + \sqrt{\frac{1}{m} \sum^{m}_{i=1} \nm{\hat{\cB}_{\cY} \circ P(\bm{x}_i) - \hat{\cB}_{\cY} \circ f(\bm{x}_i) }^2_{\cY} } + \frac{1}{\sqrt{m}} \nm{\bm{e}'}_{2;\cY}
\\
& \leq \sigma \sqrt{\frac{1}{m} \sum^{m}_{i=1} \nm{y_i - \hat{\cB}_{\cY} \circ P (\bm{x}_i) }^2_{\cY} }  + \sqrt{\frac{1}{m} \sum^{m}_{i=1} \nm{\hat{\cB}_{\cY} \circ P(\bm{x}_i) - \hat{\cB}_{\cY} \circ f(\bm{x}_i) }^2_{\cY} } + \frac{1}{\sqrt{m}} \nm{\bm{e}'}_{2;\cY} + \tau
\\
& \leq (\sigma+1) \sqrt{\frac{1}{m} \sum^{m}_{i=1} \nm{P(\bm{x}_i) -  f(\bm{x}_i) }^2_{\cY} } + \frac{\sigma+1}{\sqrt{m}} \nm{\bm{e}'}_{2;\cY} + \tau
\\
& \leq (\sigma+1) \sqrt{\frac{1}{m} \sum^{m}_{i=1} \nm{P(\bm{x}_i) -  f(\bm{x}_i) }^2_{\cY} } + (\sigma + 1) \sqrt{\frac{1}{m} \sum^{m}_{i=1} \nm{f(\bm{x}_i) -  \hat{\cB}_{\cY} \circ f(\bm{x}_i) }^2_{\cY} }
\\
&~~~~~+\frac{\sigma+1}{\sqrt{m}} \nm{\bm{e}}_{2;\cY} + \tau .
}
Here, in the penultimate step we used the fact that $\hat{\cB}_{\cY}$ is an orthogonal projection, and in the final step we used the definition of $\bm{e}'$. Substituting this into \ef{false-creek-1} now gives that
\eas{
\nm{f - \hat{f}}_{L^2_{\varrho}(\bbR^{\bbN} ; \cY)} \leq &~  \nm{f - \hat{\cB}_{\cY} \circ f}_{L^{2}_{\varrho}(\bbR^{\bbN} ; \cY)} + \nm{ f - P}_{L^2_{\varrho}(\bbR^{\bbN} ; \cY)}
\\
& + \alpha^{-1}(\sigma+1) \sqrt{\frac{1}{m} \sum^{m}_{i=1} \nm{ P(\bm{x}_i) - f(\bm{x}_i) }^2_{\cY} }
\\
& +\alpha^{-1}(\sigma+1) \sqrt{\frac{1}{m} \sum^{m}_{i=1} \nm{ f(\bm{x}_i) - \hat{\cB}_{\cY} \circ f(\bm{x}_i) }^2_{\cY} } + \alpha^{-1}(\sigma+1) \frac{\nm{\bm{e}}_{2;\cY}}{\sqrt{m}} + \alpha^{-1} \tau.
}
Now observe that $\hat{f} - P \in \cP_{S' ; \cY}$ for some set $S'$ with $|S'|_{\bm{w}} \leq 2 k$. Using \ef{dist-shift-const-bound} with $S'$ we get $\Delta(\cP_{S' ; \cY} ; \varrho , \mu) \leq \sqrt{2 k}$. We now apply \ef{OOD-ID-bound} (or, more precisely, the Hilbert-valued extension of this bound, which is proved in exactly the same way) to deduce that
\bes{
\nm{f - \hat{f}}_{L^{\infty}_{\mu}(\bbR^{\bbN}; \cY)} \leq  \sqrt{2k} \left ( \nm{f - P}_{L^2_{\varrho}(\bbR^{\bbN}; \cY)} + \nm{f - \hat{f}}_{L^2_{\varrho}(\bbR^{\bbN}; \cY)} \right ) + \nm{f - P}_{L^{\infty}_{\mu}(\bbR^{\bbN}; \cY)}.
}
Combining this with the previous estimate now gives the result. 
}

\lem{
\label{lem:alpha-bound-prob}
Let $0 < \epsilon < 1$, $\bm{x}_{1},\ldots,\bm{x}_{m} \sim_{\mathrm{i.i.d.}} \varrho$ and consider the constant $\alpha$ defined in \ef{alpha-def}. Suppose that $k \geq 1$. Then $\bbP(\alpha < 1/2) \leq \epsilon$, provided
\be{
\label{m-cond-alpha-bound}
m \geq c \cdot k \cdot \left ( \log^4(2k) + \log(1/\epsilon) \right ),
}
where $c > 0$ is a universal constant.
}
\prf{
Let $P_1,P_2 \in \cP_{S ; \cY}$ be arbitrary and observe that $P_1 - P_2$ can be expressed as $P_1 - P_2 = \sum_{\bn \in S'} c_{\bn} \Psi_{\bn}$ where $S' \subseteq \Lambda$, $|S'|_{\bm{w}} \leq 2 k$. Define $c_{\bn} = 0 \in \cY$ for $\bn \notin S'$.
Now $\{ \psi_j \}_j$ be an orthonormal basis of $\cY$ and write $c_{\bn} = \sum_j d_{\bn,j} \psi_j$ for $d_{\bn,j} \in \bbR$. Let $\bm{d}_j = (d_{\bn,j})_{\bn \in \Lambda}$. Then, by Parseval's identity,
\bes{
\nm{P_1 - P_2}^2_{L^2(\bbR^{\bbN} ; \cY)} = \sum_{\bn \in \Lambda} \nm{c_{\bn}}^2_{\cY} =  \sum_{\bn \in \Lambda} \sum_j |d_{\bn,j} |^2 = \sum_j \nm{\bm{d}_j}^2_2,
}
and
\bes{
\frac1m \sum^{m}_{i=1} \nm{ P_1(\bm{x}_i) - P_2(\bm{x}_i) }^2_{\cY} = \sum_j \nm{\bm{A} \bm{d}_j}^2_2,
}
where $\bm{A} \in \bbR^{m \times N}$ is the matrix $\bm{A} = \frac{1}{\sqrt{m}} \left ( \Psi_{\bn}(\bm{x}_i) \right )_{i \in [m],\bn \in \Lambda}$. We deduce that
\bes{
\alpha \geq \inf_{ \bm{z} \in T}  \nm{\bm{A} \bm{z}}_2 ,\qquad \text{where }T  = \left \{  \bm{z} = (z_{\bn})_{\bn \in \Lambda} \in \bbR^N,\ \nm{\bm{z}}_2 = 1,\ | \mathrm{supp}(\bm{z}) |_{\bm{w}} \leq 2k \right \}.
}
Notice that $T \subseteq \{ \bm{z} : \nm{\bm{z}}_{1,\bm{w}} \leq \sqrt{2k} \}$ by the Cauchy--Schwarz inequality. Now define the random vector $\bm{X} \in \bbR^N$ by $\bm{X} = (\Psi_{\bn}(\bm{x}) )_{\bn \in \Lambda}$ for $\bm{x} \sim \varrho$ and observe that
\be{
\label{X-z-exp}
\bbE | \ip{\bm{X}}{\bm{z}} |^2 = \nm{\bm{z}}^2_2
}
by Parseval's identity, and also that
\bes{
| \ip{\bm{X}}{\bm{e}_{\bn}} | \leq \nm{\Psi_{\bn}}_{L^{\infty}_{\varrho}(\bbR^{\bbN})} = u_{\bn} \leq w_{\bn},\quad \forall \bn \in \Lambda,
}
almost surely, due to \ef{Psi-ID-bound} and \ef{Psi-OOD-bound}. We now apply \cite[Thm.\ 2.13]{brugiapaglia2021sparse}. This gives that there are univeral constants $\kappa,c_0,c_1>0$ such that, if
\be{
\label{m-rip-cond-1}
m \geq c_0 \cdot \delta^{-2} \cdot k \cdot \log(\E |\Lambda|) \cdot \log^2(k/\delta) \log^2(1/\delta)
} 
for some $\delta \in (0,\kappa)$, then
\bes{
\sup_{\bm{z} \in T} \left | \frac1m \sum^{m}_{i=1} | \ip{\bm{X}_i}{\bm{z}} |^2 - \bbE | \ip{\bm{X}}{\bm{z}} |^2 \right | \leq c_1 \delta \left ( 1 + \sup_{\bm{z} \in T} \bbE | \ip{\bm{X}}{\bm{z}} |^2 \right )
}
with probability at least $1-2 \exp(-c_2 \delta^2 m/k )$, where $\bm{X}_1,\ldots,\bm{X}_m$ are independent copies of $\bm{X}$. Notice that $\frac1m \sum^{m}_{i=1} | \ip{\bm{X}_i}{\bm{z}} |^2 = \nm{\bm{A} \bm{z}}^2_2$. Using this, \ef{X-z-exp} and rearranging, we see that
\bes{
\alpha^2 \geq \inf_{\bm{z} \in T} \nm{\bm{A} \bm{z}}^2_2 \geq 1 - 2 c_1 \delta,
}
with the same probability. Set $\delta = \min \{ \kappa / 2 , 3/(8 c_1) \}$. Then we see that
\bes{
\bbP(\alpha < 1/2) \leq 2 \exp(-c_2 \delta^2 m / k),
}
provided \ef{m-rip-cond-1} holds with this value of $\delta$. Therefore $\bbP(\alpha < 1/2) < \epsilon$, provided \ef{m-rip-cond-1} holds, along with
\bes{
m \geq c^{-1}_2 \cdot \delta^{-2} \cdot k \cdot \log(2/\epsilon)
}
for the given value of $\delta$.
In particular, it suffices for
\be{
\label{m-rip-cond-2}
m \geq c_3 \cdot k \cdot \left ( \log(\E |\Lambda|) \cdot \log^2(2 k) + \log(1/\epsilon) \right ),
}
provided $k \geq 1$. A standard bound (see, e.g., the proof of \cite[Lem.\ 6.4]{adcock2025near}) gives that $\log(E | \Lambda |) \leq 4 \log^2(\E \lceil k \rceil ) \leq 4 \log^2(2 \E k)$. Hence \ef{m-rip-cond-2} is implied by
\be{
\label{m-rip-cond-3}
m \geq c_4 \cdot k \cdot \left ( \log^4(2k) + \log(1/\epsilon) \right ).
}
This complete the proof.
}

We are now ready to establish Theorem \ref{thm:poly-optimized-least-squares-hilbert}.

\prf{
[Proof of Theorem \ref{thm:poly-optimized-least-squares-hilbert}]

For a suitable choice of the universal constant $c$ we see that \ef{m-cond-alpha-bound} holds with $\epsilon$ replaced by $\epsilon/3$. Hence Lemma \ref{lem:alpha-bound-prob} gives that $\bbP(\alpha < 1/2) \leq \epsilon/3$. We now apply Lemma \ref{lem:poly-least-squares-alpha-bound} to deduce that
\be{
\begin{split}\label{f-fhat-no-oracle}
\nm{f - \hat{f}}_{L^{\infty}_{\mu}(\bbR^{\bbN} ; \cY)} \lesssim  ~ \sqrt{ k} \Bigg [ & \nm{f - \hat{\cB}_{\cY} \circ f}_{L^{2}_{\varrho}(\bbR^{\bbN};\cY)} + \sigma \sqrt{\frac1m \sum^{m}_{i=1} \nm{ f(\bm{x}_i) - \hat{\cB}_{\cY} \circ f (\bm{x}_i) }^2_{\cY} }
\\
& + \nm{f - P}_{L^2_{\varrho}(\bbR^{\bbN}; \cY)} + \sigma \sqrt{\frac1m \sum^{m}_{i=1} \nm{ f(\bm{x}_i) - P(\bm{x}_i) }^2_{\cY} } 
 + \sigma \frac{\nm{\bm{e}}_{2;\cY}}{\sqrt{m}} +  \tau \Bigg ] 
\\
& + \nm{f - P}_{L^{\infty}_{\mu}(\bbR^{\bbN}; \cY)}
\end{split}
}
for all $P \in \cP_{\cS}$, with probability at least $1-\epsilon/3$.

We next fix a suitable $P$. Let $S \subset \cF$, $|S|_{\bm{w}} \leq k$, be the set specified in Theorem \ref{thm:weighted-summability-coeffics-1} (notice that \ef{b-r-cond-1} holds by assumption). Now write $S = S_1 \cup S_2$, where $S_1 = S \cap \Lambda$ and $S_2 = S \cap \Lambda^c$. Notice that $|S_i|_{\bm{w}} \leq k$, $i = 1,2$, and therefore $S_1 \in \cS$. Finally, let $P = f_{S_1}$.
Now consider the middle term of \ef{f-fhat-no-oracle}. Let $X$ be the random variable $\frac1m \sum^{m}_{i=1} \nm{P(\bm{x}_i) - f(\bm{x}_i) }^2_{\cY}$ and notice that $\bbE(X) = \nm{f - P}^2_{L^2_{\varrho}(\bbR^{\bbN} ; \cY)}$. Hence Markov's inequality gives that
\bes{
\bbP \left ( \sqrt{\frac1m \sum^{m}_{i=1} \nm{ P(\bm{x}_i) - f(\bm{x}_i) }^2_{\cY} } \geq \frac{\nm{f - P}_{L^2_{\varrho}(\bbR^{\bbN} ; \cY)}}{ \sqrt{\epsilon/3} } \right ) \leq \epsilon/3.
}
Similarly, Markov's inequality also gives that
\bes{
\bbP \left ( \sqrt{\frac1m \sum^{m}_{i=1} \nm{ f(\bm{x}_i) - \hat{\cB}_{\cY} \circ f (\bm{x}_i) }^2_{\cY} } \geq \frac{\nm{f - \hat{\cB}_{\cY} \circ f}_{L^2_{\varrho}(\bbR^{\bbN} ; \cY)}}{ \sqrt{\epsilon/3} } \right ) \leq \epsilon/3.
}
Using these two estimates, \ef{f-fhat-no-oracle} and the union bound, we deduce that
\be{
\label{f-fhat-no-oracle-2}
\begin{split}
\nm{f - \hat{f}}_{L^{\infty}_{\mu}(\bbR^{\bbN} ; \cY)} \lesssim &~ \sqrt{ k} \left ( \frac{\sigma}{\sqrt{\epsilon}} \nm{f - \hat{\cB}_{\cY} \circ f}_{L^2_{\varrho}(\bbR^{\bbN};\cY)} +\frac{\sigma}{\sqrt{\epsilon}} \nm{f - f_{S_1} }_{L^2_{\varrho}(\bbR^{\bbN}; \cY)}  + \sigma \frac{\nm{\bm{e}}_{2;\cY}}{\sqrt{m}} +  \tau \right )  
\\
&  + \nm{f - f_{S_1}}_{L^{\mu}_{\varrho}(\bbR^{\bbN}; \cY)}
 \end{split}
}
with probability at least $1-\epsilon$.

Finally, we bound the two terms involving $f - f_{S_1}$ in \ef{f-fhat-no-oracle-2}. Let $\bm{c} = (c_{\bn})_{\bn \in \cF}$ be the coefficients of $f$. Since $S_1 = S \cap \Lambda$, we have
\eas{
\nm{f - f_{S_1}}_{L^2_{\varrho}(\bbR^{\bbN} ; \cY)} & = \nm{\bm{c} - \bm{c}_{S_1}}_{2 ; \cY} 
\\
& \leq \nm{\bm{c} - \bm{c}_S}_{2 ; \cY} + \nm{\bm{c} - \bm{c}_{\Lambda}}_{2;\cY} 
\\
& = \nm{\bm{c} - \bm{c}_S}_{2,\bm{w} ; \cY} + \nm{\bm{c} - \bm{c}_{\Lambda}}_{2,\bm{w} ; \cY},
}
where $\bm{w} = (w_{\bn})_{\bn \in \cF}$ is as in \ef{Psi-OOD-bound} and 
\bes{
\nm{f - f_{S_1}}_{L^{\infty}_{\mu}(\bbR^{\bbN} ; \cY)} = \nm{\bm{c} - \bm{c}_{S_1}}_{1,\bm{w}} \leq \nm{\bm{c} - \bm{c}_S}_{1,\bm{w} ; \cY} + \nm{\bm{c} - \bm{c}_{\Lambda}}_{1,\bm{w} ; \cY}.
}
By definition of $S$, Theorem \ref{thm:weighted-summability-coeffics-1} gives that 
\bes{
\nm{\bm{c} - \bm{c}_S}_{q,\bm{w} ; \cY} \leq C(\bm{b},\varepsilon,\bm{\omega},p) k^{1/q-1/p}
}
for all $q \in (p,2]$. Theorem \ref{thm:weighted-summability-coeffics-2} (notice that \ef{b-r-cond-2} holds by assumption) implies that there exists an anchored set $T$ with $|T| \leq \lceil k \rceil$ such that
\bes{
\nm{\bm{c} - \bm{c}_{T}}_{q,\bm{w} ; \cY} \leq C(\bm{b},\varepsilon,\bm{\omega},p) \lceil k \rceil^{1/q-1/p}
}
for all $q \in (p,2]$. Moreover, $\Lambda$ contains all anchored sets of size at most $n$ \cite[Prop.\ 2.18]{adcock2022sparse}. Hence $T \subseteq \Lambda$ and we deduce that
\bes{
\nm{\bm{c} - \bm{c}_{\Lambda}}_{q,\bm{w} ; \cY} \leq C(\bm{b},\varepsilon,\bm{\omega},p) \lceil k \rceil^{1/q-1/p}
}
for all $q \in (p,2]$. Combining the various estimates, we conclude that
\bes{
\nm{f - f_{S_1}}_{L^2_{\varrho}(\bbR^{\bbN} ; \cY)} \leq C(\bm{b},\varepsilon,\bm{\omega},p) k^{1/2-1/p} ,\qquad \nm{f - f_{S_1}}_{L^{\infty}_{\mu}(\bbR^{\bbN} ; \cY)}  \leq C(\bm{b},\varepsilon,\bm{\omega},p) k^{1-1/p} .
}
Substituting this into \ef{f-fhat-no-oracle-2} we get
\be{
\label{f-fhat-partial-2}
\nm{f - \hat{f}}_{L^{\infty}_{\mu}(\bbR^{\bbN} ; \cY)} \leq \frac{\sigma}{\sqrt{\epsilon}} C(\bm{b},\varepsilon,\bm{\omega},p) k^{1-1/p} + \frac{\sigma \sqrt{k}}{\sqrt{\epsilon}} \nm{f - \hat{\cB}_{\cY} \circ f}_{L^2_{\varrho}(\bbR^{\bbN} ; \cY)}
+ \frac{\sigma \sqrt{k}\nm{\bm{e}}_2}{\sqrt{m}} + \sqrt{k} \tau .
}
with probability at least $1-\epsilon$, as required.
}

\prf{[Proof of Theorem \ref{thm:poly-optimized-least-squares}]
Let $c$ be the constant from Lemma \ref{lem:alpha-bound-prob} and assume, without loss of generality, that $c \geq 1/\log^4(2)$. Set
\be{
\label{k-def}
k = \frac{m}{c \left ( \log^4(m) + \log(1/\epsilon) \right )}.
}
Notice that $k \leq m / (c \log^4(2))$ whenever $m \geq 2$ and therefore $k \leq m$. Now let $\bar{m} = \bar{m}(\epsilon) \in \bbN$ be the smallest $m \in \bbN$, $m \geq 2$, such that $k \geq 1$. Let $\cY = \hat{\cY} = \bbR$ with the Euclidean inner product, so that, in particular, $\hat{\cB}_{\cY} = \cI_{\cY}$ is the identity. We now apply Theorem \ref{thm:poly-optimized-least-squares-hilbert} with this value of $k$, $\sigma = 1$ and $\tau = 0$. This yields the result.
}

\subsection{Proofs of Theorems \ref{thm:deep-learning-extrap} and \ref{thm:operator-learning-extrap}}\label{ss:deep-learning-proofs}

We now prove Theorems \ref{thm:deep-learning-extrap} and \ref{thm:operator-learning-extrap}. Our main focus in the proof of Theorem \ref{thm:operator-learning-extrap}, since Theorem \ref{thm:deep-learning-extrap} will then follow as a special case.

Our broad approach is based on well-established technique of emulation of polynomials via DNNs, or in this case, DNOs, \cite{adcock2025near,de-ryck2021approximation,guhring2021approximation,li2020better,lu2021deep,yarotsky2017error,schwab2019deep,mhaskar1996neural,opschoor2022exponential,schwab2023deepa} (see \cite[\S 7.1]{adcock2024learning} for a review). As noted earlier, emulation techniques are normally used to assert the existence of DNNs/DNOs with certain approximation guarantees. Conversely, we follow ideas developed in \cite{adcock2021deep,adcock2025near,adcock2024optimalb} that emulate a whole polynomial training problem as a DNO training problem, thus showing the existence of DL strategies with guaranteed error bounds. 

Our proof proceeds in three steps. These steps are broadly similar to those used in \cite{adcock2024optimalb}, which develops generalization bounds learning holomorphic operators, but considers only the in-distribution setting. Our results generalize these by considering arbitrary test distributions.
First, we introduce the polynomial training problem
\be{
\label{poly-DNN-training}
\min_{P \in \cP_{\cS ; \hat{\cY}} } \frac1m \sum^{m}_{i=1} \nm{Y_i - P \circ \hat{\cE}_{\cX}(X_i)}^2_{\cY},
}
where $\cS$ is as defined in \ef{S-def} with $\Lambda$ as in \ef{Lambda-def}. The second step is to then emulate the relevant polynomials in the hypothesis class a DNOs $\cN$ (\S \ref{sss:dnn-class-constr}). Next, we show that the weights corresponding to an approximate minimizer of the associated DNO training problem yield coefficients of a polynomial that is an approximate minimizer of \ef{poly-DNN-training} (\S \ref{sss:from-dnn-to-poly}). Using this, we establish a general out-of-distribution error bound for the DNO training problem, based on the analogous bound for the polynomial training problem. Finally, we specialize this to the settings of Theorems \ref{thm:deep-learning-extrap} and \ref{thm:operator-learning-extrap} to conclude the proof (\S \ref{sss:final-arguments}).

\subsubsection{Construction of $\cN$}\label{sss:dnn-class-constr}

Let $0 < k \leq d_{\cX}$ be a number whose value will be chosen later.  We first require the following lemma.

\lem{
Let $\Gamma \subseteq \Lambda$ and $m(\Gamma) = \max_{\bn \in \Gamma} \nm{\bn}_1 < \infty$. Then there exists a fully-connected family $\cN_o$ of tanh DNNs $\bbR^{\lceil k \rceil} \rightarrow \bbR$ with
\bes{
\mathrm{width}(\cN_o) \lesssim m(\Gamma),\quad \mathrm{depth}(\cN_o) \lesssim \log(m(\Gamma)),
}
such that, for any $0 < \delta < 1$, $\bm{\omega} \geq \bm{1}$ and $\bm{n} \in \Gamma$, there is a $N_{\bn} \in \cN_o$ satisfying
\bes{
\sup_{\bm{x} \in D_{\bm{\omega}} } | N_{\bn}(\bm{x}) - \Psi_{\bn}(\bm{x}) | \leq \delta.
}
Moreover, the zero network $0 : \bm{x} \mapsto 0$ also belongs to $\cN_o$ (trivially, since $\mathrm{tanh}(0) = 0$).
}

Note here we slightly abuse notation and consider the DNNs $N_{\bn}$ both as maps with domain $\bbR^{\lceil k \rceil}$ and maps with domain $\bbR^{\bbN}$ that depend only on the first $\lceil k \rceil$ entries of the given input. This is permissible, since the Legendre polynomials $\Psi_{\bn}$, $\bn \in \Lambda$, depend on the first $\lceil k \rceil$ variables only.
This lemma is based on \cite[Lem.\ D.9]{adcock2024optimalb}, which is itself based on \cite{de-ryck2021approximation,opschoor2022exponential}. Notice that \cite[Lem.\ D.9]{adcock2024optimalb} considers only the domain $D = D_{\bm{1}}$. However, the proof extends to $\bm{\omega} \geq \bm{1}$. In particular, the width and depth bounds do not depend on $\bm{\omega}$ (although the weights and biases of $N_{\bn}$ do). This is a particular feature of $\tanh$ DNNs, which makes them particularly desirable over, for instance, ReLU DNNs for out-of-distribution generalization.

Now fix $\delta > 0$ (its value will be chosen later in the proof), let $\Gamma = \cup_{S \in \cS} S$ and $\cN_o$ and $N_{\bn}$, $\bn \in \Gamma$, be as in the above lemma. We now estimate $m(\Gamma)$. Let $\bn \in \Gamma$. Since $\bn \in S$ for some $S \in \cS$ and, by construction, $|S|_{\bm{w}} \leq k$, we have
\bes{
\prod^{\infty}_{j=1} (2 n_j+1) \zeta(\omega_j)^{2 n_j} = w^2_{\bn} \leq k .
}
Since $\zeta(\omega) = \omega + \sqrt{\omega^2-1} \geq \omega$, it follows that
\bes{
\nm{\bn}_1 \exp(2 \nm{\bn}_1 \log(\omega_{\min}) ) \leq k,
}
where $\omega_{\min} = \min \{\omega_j \} > 1$ by assumption, and therefore
\bes{
m(\Gamma) \lesssim \frac{\log(k)}{\log(\omega_{\min})}.
}
Notice that $|S| \leq |S|_{\bm{w}} \leq k$ for any $S \in \cS$. We now define the family $\cN$ of DNNs $N : \bbR^{d_{\cX}} \rightarrow \bbR^{d_{\cY}}$ as
\bes{
\cN = \left \{ N = \bm{C} \begin{bmatrix} N_{\bn_1} \\ \vdots \\ N_{\bn_{|S|}} \\ 0 \\ \vdots \\ 0 \end{bmatrix} : S = \{ \bn_1,\ldots,\bn_{|S|} \} \in \cS,\ \bm{C} \in \bbR^{d_{\cY} \times \lceil k \rceil } \right \}.
}
Here $0$ denotes the zero network.
As before, we slightly abuse notation by allowing $N$ to be a function whose domain is either $\bbR^{\lceil k \rceil}$, $\bbR^{d_{\cX}}$ or $\bbR^{\infty}$ depending on the context. Notice that this family satisfies
\be{
\label{N-width-bd}
\mathrm{width}(\cN) \leq \mathrm{width}(\cN_o) \max_{S \in \cS} |S|  \lesssim k m(\Gamma) \lesssim \frac{k \log(k)}{\log(\omega_{\min})}
}
and
\be{
\label{N-depth-bd}
\mathrm{depth}(\cN) \leq \mathrm{depth}(\cN_o) \lesssim \log(m(\Gamma)) \lesssim \log \left ( \frac{\log(k)}{\log(\omega_{\min})} \right ).
}
Finally, given $N \in \cN$ and its associated matrix $\bm{C} \in \bbR^{d_{\cY} \times \lceil k \rceil}$. Let $\bm{c}_1,\ldots,\bm{c}_{\lceil k \rceil}$ denote the columns of $\bm{C}$, so that $N$ can be expressed as
\bes{
N = \sum^{|S|}_{i=1} \bm{c}_i N_{\bm{n}_i}.
}
Moreover, we have
\bes{
\hat{\cD}_{\cY} \circ N = \sum^{|S|}_{i=1} \hat{\cD}_{\cY}(\bm{c}_i) N_{\bm{n}_i} = \sum^{|S|}_{i=1} c_i N_{\bn_i},\quad \text{where } c_i \in \hat{\cY},
}
where $\hat{\cY} = \hat{\cD}_{\cY}(\bbR^{d_{\cY}}) \subseteq \cY$. Therefore, we can associate $\cN$ with the space of functions $\cQ_{\cS ; \hat{\cY}}$ given by 
\bes{
\cQ_{\cS ; \hat{\cY}} = \left \{ \sum_{\bn \in S} c_{\bn} N_{\bn} : c_{\bn} \in \hat{\cY},\ S \in \cS \right \}.
}

\subsubsection{Polynomial training problem and approximate DNN minimizers}\label{sss:from-dnn-to-poly}


We first require the following lemma.

\lem{
\label{lem:closeness-polys-dnns}
Let $\cS$, $\{ N_{\bm{n}} \}_{\bm{n} \in \Gamma}$ be as defined above and $(\cZ , \nm{\cdot}_{\cZ})$ be a Hilbert space. Let $S \in \cS$ and suppose that
\bes{
p = \sum_{\bn \in S} c_{\bn} \Psi_{\bn},\qquad \tilde{p} = \sum_{\bn \in S} c_{\bn} N_{\bn},
}
where $c_{\bn} \in \cZ$. Then
\bes{
\sup_{\bm{x} \in D_{\bm{\omega}}} \nm{p(\bm{x}) - \tilde{p}(\bm{x})}_{\cZ} \leq \delta \sqrt{|S|} \nm{p}_{L^2_{\varrho}(\bbR^{\bbN};\cZ)}.
}
}
\prf{
We have
\bes{
\nm{p(\bm{x}) - \tilde{p}(\bm{x})}_{\cZ} \leq \sum_{\bn \in S} \nm{c_{\bn}}_{\cZ} | \Psi_{\bn}(\bm{x}) - N_{\bn}(\bm{x}) | \leq \delta \sum_{\bn \in S} \nm{c_{\bn}}_{\cZ}.
}
The result now follows from the Cauchy--Schwarz inequality and Parseval's identity.
}

We now show that approximate minimizers of the DNO training problem yield approximate minimizers of the polynomial training problem.

\lem{
\label{lem:poly-dnn-minimizers}
Suppose that $\alpha > 0$, where $\alpha$ is as in \ef{alpha-def} with $\cY$ replaced by $\hat{\cY}$, let $\cN$ be the DNN model class defined above and $\hat{N}$ be any $(\sigma,\tau)$-approximate minimizer of \ef{DNN-training-prob-operators}. Write
\bes{
\cD_{\cY} \circ \hat{N} = \sum_{\bn \in S} \hat{c}_{\bn} N_{\bn} \in \cQ_{\cS ; \hat{\cY}}
}
where $S \in \cS$ and $\hat{c}_{\bn} \in \hat{\cY}$, 
and define the polynomial
\bes{
\hat{P} = \sum_{\bn \in S} \hat{c}_{\bn} \Psi_{\bn} \in \cP_{\cS ; \hat{\cY}}.
}
Then $\hat{P}$ is a $(\sigma',\tau')$-approximate minimizer of \ef{poly-DNN-training}, where
\be{
\label{sigma-tau-prime}
\begin{split}
\sigma' & \leq \sigma(1+\delta \sqrt{k}/\alpha)
\\
\tau' & \leq \tau + \frac{\sigma \delta \sqrt{k}}{\alpha} \left ( \nm{F}_{L^{\infty}_{\nu}(\cX ; \cY)} + \frac{1}{\sqrt{m}} \nm{\bm{E}}_{2;\cY} \right ) + \delta \sqrt{k} \nm{\hat{P}}_{L^2_{\varrho}(\bbR^{\bbN};\cY)}.
\end{split}
}
}

\prf{
Let $P = \sum_{\bm{n} \in S} c_{\bm{n}} \Psi_{\bm{n}} \in P_{\cS ; \hat{\cY}}$ be arbitrary and $N \in \cN$ be such that $\hat{\cD}_{\cY} \circ N = \sum_{\bm{n} \in S} c_{\bm{n}}  N_{\bm{n}}$. Then
\eas{
\sqrt{\frac1m \sum^{m}_{i=1} \nm{Y_i - \hat{P} \circ \hat{\cE}_{\cX}(X_i)}^2_{\cY}} \leq &~ \sqrt{\frac1m \sum^{m}_{i=1} \nm{Y_i -  \hat{\cD}_{\cY} \circ \hat{N} \circ \hat{\cE}_{\cX}(X_i)}^2_{\cY}} 
\\
& + \sqrt{\frac1m \sum^{m}_{i=1} \nm{\hat{P} \circ \hat{\cE}_{\cX}(X_i) -  \hat{\cD}_{\cY} \circ \hat{N} \circ \hat{\cE}_{\cX}(X_i)}^2_{\cY}}
\\
\leq &~ \sigma \sqrt{\frac1m \sum^{m}_{i=1} \nm{Y_i -  \hat{\cD}_{\cY} \circ N \circ \hat{\cE}_{\cX}(X_i)}^2_{\cY}} + \tau
\\
& + \sqrt{\frac1m \sum^{m}_{i=1} \nm{\hat{P} \circ \hat{\cE}_{\cX}(X_i) -  \hat{\cD}_{\cY} \circ \hat{N} \circ \hat{\cE}_{\cX}(X_i)}^2_{\cY}}
\\
\leq &~ \sigma \sqrt{\frac1m \sum^{m}_{i=1} \nm{Y_i -  P \circ \hat{\cE}_{\cX}(X_i)}^2_{\cY}} + \tau
\\
& +  \sqrt{\frac1m \sum^{m}_{i=1} \nm{\hat{P} \circ \hat{\cE}_{\cX}(X_i) -  \hat{\cD}_{\cY} \circ \hat{N} \circ \hat{\cE}_{\cX}(X_i)}^2_{\cY}}
\\
& + \sigma\sqrt{\frac1m \sum^{m}_{i=1} \nm{P \circ \hat{\cE}_{\cX}(X_i) -  \hat{\cD}_{\cY} \circ N \circ \hat{\cE}_{\cX}(X_i)}^2_{\cY}}
}
and therefore
\eas{
\sqrt{\frac1m \sum^{m}_{i=1} \nm{Y_i - \hat{P} \circ \hat{\cE}_{\cX}(X_i)}^2_{\cY}} \leq &~ \sigma \sqrt{\frac1m \sum^{m}_{i=1} \nm{Y_i -  P \circ \hat{\cE}_{\cX}(X_i)}^2_{\cY}} + \tau 
\\
& +  \nm{\hat{P} - \hat{\cD}_{\cY} \circ \hat{N}}_{L^{\infty}_{\mu}(\bbR^{\bbN} ; \cY)} +  \sigma \nm{P - \hat{\cD}_{\cY} \circ N}_{L^{\infty}_{\mu}(\bbR^{\bbN} ; \cY)} ,
}
where we recall that $\hat{\cE}_{\cX}(X_i) \in D$ by Assumption \ref{ass:opl}(ii). We now apply Lemma \ref{lem:closeness-polys-dnns} to the latter two terms to get, after recalling that $|S| \leq k$ since $S \in \cS$,
\eas{
\sqrt{\frac1m \sum^{m}_{i=1} \nm{Y_i - \hat{P} \circ \hat{\cE}_{\cX}(X_i)}^2_{\cY}} \leq &~ \sigma \sqrt{\frac1m \sum^{m}_{i=1} \nm{Y_i -  P \circ \hat{\cE}_{\cX}(X_i)}^2_{\cY}} + \tau 
\\
& + \delta \sqrt{k} \nm{\hat{P}}_{L^2_{\varrho}(\bbR^{\bbN};\cY)}+ \sigma \delta \sqrt{k} \nm{P}_{L^2_{\varrho}(\bbR^{\bbN};\cY)} .
}
It remains to estimate the final term. Applying \ef{alpha-def} with $P_1 = P$ and $P_2 = 0$, we get
\eas{
\nm{P}_{L^2_{\varrho}(\bbR^{\bbN};\cY)} \leq &~ \alpha^{-1} \sqrt{\frac1m \sum^{m}_{i=1} \nm{P \circ \hat{\cE}_{\cX}(X_i)}^2_{\cY} } 
\\
\leq & ~ \alpha^{-1} \sqrt{\frac1m \sum^{m}_{i=1} \nm{Y_i - P \circ \hat{\cE}_{\cX}(X_i)}^2_{\cY} } + \alpha^{-1} \sqrt{\frac1m \sum^{m}_{i=1} \nm{Y_i}^2_{\cY}}
\\
\leq &~\alpha^{-1} \sqrt{\frac1m \sum^{m}_{i=1} \nm{Y_i - P \circ \hat{\cE}_{\cX}(X_i)}^2_{\cY} } + \alpha^{-1} \nm{F}_{L^{\infty}_{\nu}(\cX ; \cY)} + \frac{\alpha^{-1}}{\sqrt{m}} \nm{\bm{E}}_{2;\cY}.
}
Substituting this into the previous expression and recalling that $P$ was arbitrary now gives the result.
}

With this to hand, we now establish the following general out-of-distribution error bound for the DNO training problem.

\thm{[Out-of-distribution generalization bound for the DNO training problem]
\label{thm:opl-ood-app}
Let $1 \leq k \leq d_{\cX}$, $\tau > 0$, $\sigma \geq 1$, $0 < \epsilon < 1$ and $\cN$ be the DNN model class defined above, where $\delta$ satisfies $\delta \leq (4\sqrt{k})^{-1}$. Suppose that Assumption \ref{ass:opl} holds and Assumption \ref{ass:opl-ood} holds with $\bm{\omega} \geq \bm{1}$ satisfying $\omega_{\min} = \min\{\omega_j\} > 1$. Let $F \in \cH(\bm{b},\varepsilon, \cX , \cY)$, where $\bm{b}$ satisfies \ef{b-r-cond-2-main} for some $0 < p <1$, and draw $X_1,\ldots,X_m \sim_{\mathrm{i.i.d.}} \nu$, where $m$ satisfies 
\bes{
m \geq c \cdot k \cdot \left ( \log^4(2k) + \log(1/\epsilon) \right )
}
for some universal constant $c > 0$. Then, with probability at least $1-\epsilon$, any $(\sigma,\tau)$-approximate minimizer $\hat{N}$ of \ef{DNN-training-prob-operators} is such that the approximation $\hat{F} = \hat{\cD}_{\cY} \circ \hat{N} \circ \hat{\cE}_{\cX}$ satisfies
\eas{
\nm{F - \hat{F} }_{L^{\infty}_{\mu}(\cX ; \cY)} \lesssim &~ \frac{\sigma C(\bm{b},\varepsilon,\bm{\omega},p)}{\sqrt{\epsilon}} k^{1-1/p}
 + \frac{\sigma \sqrt{k}}{\sqrt{\epsilon}} \nm{F - \hat{\cB}_{\cY}  \circ F }_{L^{2}_{\nu}(\cX;\cY)} + \frac{\sigma\sqrt{k} \nm{\bm{E}}_{2;\cY}}{\sqrt{m}} 
 \\
 & + \sqrt{k} \tau + \sigma \delta k,
}
where $E = (E_i)^{m}_{i=1} \in \cY^m$. Further, the class $\cN$ satisfies
\bes{
\mathrm{width}(\cN) \lesssim \frac{k \log(k)}{\log(\omega_{\min})}, \qquad \mathrm{depth}(\cN) \lesssim \log \left ( \frac{\log(k)}{\log(\omega_{\min})} \right ).
}
}
\prf{
The width and depth bounds for $\cN$ follow immediately from the construction and \ef{N-width-bd}-\ef{N-depth-bd}. We divide the remainder of the proof into a series of steps.

\pbk
\textit{Step 1: Reduction to polynomials.}
Let $\hat{P}$ be as in Lemma \ref{lem:poly-dnn-minimizers}. Then 
\be{
\label{easter-1}
\nm{F - \hat{F} }_{L^{\infty}_{\mu}(\cX ; \cY)} \leq \nm{F - \hat{P} \circ \hat{\cE}_{\cX} }_{L^{\infty}_{\mu}(\cX ; \cY)} + \nm{\hat{P} \circ \hat{\cE}_{\cX} - \hat{\cD}_{\cY} \circ \hat{N} \circ \hat{\cE}_{\cX}  }_{L^{\infty}_{\mu}(\cX ; \cY)}.
}
For the second term, we recall first that $\hat{P}$ and $\hat{N}$ depend on the first $\lceil k \rceil \leq d_{\cX}$ variables only. Therefore, we can write
\bes{
\nm{\hat{P} \circ\hat{\cE}_{\cX} - \hat{\cD}_{\cY} \circ \hat{N} \circ \hat{\cE}_{\cX}  }_{L^{\infty}_{\mu}(\cX ; \cY)} = \nm{\hat{P} - \hat{\cD}_{\cY} \circ \hat{N} }_{L^{\infty}_{\cE_{\cX} \sharp \mu}(\bbR^{\bbN}; \cY)}
}
We now apply Lemma \ref{lem:closeness-polys-dnns} to get
\be{
\label{easter-2}
\nm{\hat{P} \circ \hat{\cE}_{\cX} - \hat{\cD}_{\cY} \circ \hat{N} \circ \hat{\cE}_{\cX}  }_{L^{\infty}_{\mu}(\cX ; \cY)} \leq \delta \max_{S \in \cS} \sqrt{|S|} \nm{P}_{L^2_{\varrho}(\bbR^{\bbN};\cY)} \leq \delta \sqrt{k}  \nm{\hat{P}}_{L^2_{\varrho}(\bbR^{\bbN};\cY)}.
}
For the first term of \ef{easter-1}, we write $f = F \circ \cD_{\cX} : \bbR^{\infty} \rightarrow \cY$ and recall that $\hat{P} \circ \hat{\cE}_{\cX} = \hat{P} \circ \cE_{\cX}$ to get
\bes{
\nm{F - \hat{P} \circ \hat{\cE}_{\cX} }_{L^{\infty}_{\mu}(\cX ; \cY)} = \nm{f - \hat{P} }_{L^{\infty}_{\cE_{\cX} \sharp \mu}(\bbR^{\bbN} ; \cY)}.
}
Hence
\bes{
\nm{F - \hat{F} }_{L^{\infty}_{\mu}(\cX ; \cY)} \leq \nm{f - \hat{P}}_{L^{\infty}_{\cE_{\cX} \sharp \mu}(\bbR^{\bbN};\cY)} + \delta \sqrt{k} \nm{\hat{P}}_{L^2_{\varrho}(\bbR^{\bbN};\hat{\cY})}.
}

\pbk
\textit{Step 2: $\hat{P}$ is an approximate minimizer.} 
Since $\hat{P}$ comes from Lemma \ref{lem:poly-dnn-minimizers}, it is a $(\sigma',\tau')$-minimizer of \ef{poly-DNN-training}, where $(\sigma',\tau')$ are given by \ef{sigma-tau-prime}. Recall that any $P \in \cP_{\cS ; \hat{\cY}}$ depends on the first $\lceil k \rceil$ entries of an input $\bm{x}$. Since $k \leq d_{\cX}$, we see that $P \circ \hat{\cE}_{\cX}(X_i) = P(\bm{x}_i)$, where $\bm{x}_i = \cE_{\cX}(X_i)$ 
is a uniformly distributed random variable on $D = (-1,1)^{\bbN}$. We also have $F(X_i) = f(\bm{x}_i)$ and therefore $Y_i = F(X_i) + E_i$ is equal to $y_i = f(\bm{x}_i) + e_i$ with $E_i = e_i$.
Hence \ef{poly-DNN-training} is equivalent to \ef{poly-least-squares-hilbert}. 

The conditions of Theorem \ref{thm:poly-optimized-least-squares-hilbert} hold by assumption, hence we may apply it with $\epsilon$ to deduce that
\eas{
\nm{f - \hat{P}}_{L^{\infty}_{\mu}(\bbR^{\bbN} ; \cY)} \lesssim &~ \frac{\sigma' C(\bm{b},\varepsilon,\bm{\omega},p)}{\sqrt{\epsilon}} k^{1-1/p}
 + \frac{\sigma' \sqrt{k}}{\sqrt{\epsilon}} \nm{f - \hat{\cB}_{\cY} \circ f}_{L^{2}_{\varrho}(\bbR^{\bbN};\cY)}
\\
&  + \frac{\sigma' \sqrt{k} \nm{\bm{E}}_{2;\cY}}{\sqrt{m}} 
 + \sqrt{k} \tau'
,
}
where $\bm{E} = (\bbR^{\bbN}_i)^{m}_{i=1}$, with probability at least $1-\epsilon/2$. 
Plugging in the values \ef{sigma-tau-prime} for $\sigma',\tau'$ combining with the previous bound, we deduce that
\eas{
\nm{F - \hat{F} }_{L^{\infty}_{\mu}(\cX ; \cY)} \lesssim &~ \frac{\sigma (1+\delta \sqrt{k} / \alpha) C(\bm{b},\varepsilon,\bm{\omega},p)}{\sqrt{\epsilon}} k^{1-1/p}
\\
&  + \frac{\sigma (1+\delta \sqrt{k} / \alpha) \sqrt{k}}{\sqrt{\epsilon}} \nm{f - \hat{\cB}_{\cY} \circ f}_{L^{2}_{\varrho}(\bbR^{\bbN};\cY)}
 + \frac{\sigma (1+\delta \sqrt{k} / \alpha) \sqrt{k} \nm{\bm{E}}_{2;\cY}}{\sqrt{m}} 
\\
& + \sqrt{k} \tau + \frac{\sigma \delta k}{\alpha} \left ( \nm{F}_{L^{\infty}_{\nu}(\cX ; \cY)} 
 + \frac{1}{\sqrt{m}} \nm{\bm{E}}_{2;\cY} \right )
 \\
& + \delta k \nm{\hat{P}}_{L^2_{\varrho}(\bbR^{\bbN};\cY)} + \delta \sqrt{k} \nm{\hat{P}}_{L^2_{\varrho}(\bbR^{\bbN};\cY)},
}
where $\alpha$ is given by \ef{alpha-def}. Now recall from the proof of Theorem \ref{thm:poly-optimized-least-squares-hilbert} $\alpha \geq 1/2$ as part of the probabilistic event that holds with probability at least $1-\epsilon$ that guarantees its main error bound. Since $\delta \leq 1/(4\sqrt{k})$ by assumption, we obtain
\be{
\label{F-Fhat-split-mid}
\begin{split}
\nm{F - \hat{F} }_{L^{\infty}_{\mu}(\cX ; \cY)} \lesssim &~ \frac{\sigma C(\bm{b},\varepsilon,\bm{\omega},p)}{\sqrt{\epsilon}} k^{1-1/p}
 + \frac{\sigma \sqrt{k}}{\sqrt{\epsilon}} \nm{f - \hat{\cB}_{\cY} \circ f}_{L^{2}_{\varrho}(\bbR^{\bbN};\cY)}
 \\
 & + \frac{\sigma\sqrt{k} \nm{\bm{E}}_{2;\cY}}{\sqrt{m}} 
 + \sqrt{k} \tau + \sigma \delta k \left ( \nm{F}_{L^{\infty}_{\nu}(\cX ; \cY)} 
 + \frac{1}{\sqrt{m}} \nm{\bm{E}}_{2;\cY} \right )
 \\
& + \delta k \nm{\hat{P}}_{L^2_{\varrho}(\bbR^{\bbN};\cY)} .
\end{split}
}

\pbk
\textit{Step 3: Estimating $\nm{\hat{P}}_{L^2_{\varrho}(\bbR^{\bbN};\cY)}$}
The next step is to estimate the final term. Using \ef{alpha-def} and the fact that $0 , \hat{P} \in \cP_{\cS ; \hat{\cY}}$, we get
\eas{
\nm{\hat{P}}_{L^2_{\varrho}(\bbR^{\bbN};\cY)} & \leq \alpha^{-1} \sqrt{\frac1m \sum^{m}_{i=1} \nm{\hat{P}(\bm{x}_i) }^2_{\cY} }
 \\
& \leq \alpha^{-1} \sqrt{\frac1m \sum^{m}_{i=1} \nm{\hat{\cD}_{\cY} \circ \hat{N}(\bm{x}_i) }^2_{\cY} } + \alpha^{-1}\nm{\hat{\cD}_{\cY} \circ \hat{N} - \hat{P}}_{L^{\infty}_{\varrho}(\bbR^{\bbN};\cY)}.
}
We now use Lemma \ref{lem:closeness-polys-dnns} and the fact that $\hat{F} = \hat{\cD}_{\cY} \circ \hat{N} \circ \hat{\cE}_{\cX}$ to obtain
\bes{
\nm{\hat{P}}_{L^2_{\varrho}(\bbR^{\bbN};\cY)} \leq \alpha^{-1}  \sqrt{\frac1m \sum^{m}_{i=1} \nm{\hat{F}(\bm{x}_i) }^2_{\cY} } + \alpha^{-1} \delta \sqrt{k} \nm{\hat{P}}_{L^2_{\varrho}(\bbR^{\bbN};\cY)}.
}
We now rearrange and use the facts that $\hat{F}$ is as $(\sigma,\tau)$-minimizer and the zero network belongs to $\cN$, to get
\eas{
\nm{\hat{P}}_{L^2_{\varrho}(\bbR^{\bbN};\cY)}  &\leq \frac{\alpha^{-1}}{1-\alpha^{-1} \delta \sqrt{k}} \left ( \sqrt{\frac1m \sum^{m}_{i=1} \nm{\hat{F}(X_i) - Y_i }^2_{\cY} } + \frac{1}{\sqrt{m}} \nm{\bm{Y}}_{2;\cY} \right )
\\
& \leq \frac{\alpha^{-1}}{(1-\alpha^{-1} \delta \sqrt{k}) } \left ( \frac{1+\sigma}{\sqrt{m}}\nm{\bm{Y}}_{2;\cY} + \tau \right ).
}
We now recall that $\alpha \geq 1/2$, $\delta \leq 1/(4 \sqrt{k})$ and $Y_i = F(X_i) + E_i$. This gives
\bes{
\nm{\hat{P}}_{L^2_{\varrho}(\bbR^{\bbN};\cY)} \lesssim (1+\sigma) \nm{F}_{L^{\infty}_{\nu}(\cX ; \cY)} + \frac{1+\sigma}{\sqrt{m}} \nm{\bm{E}}_{2;\cY} + \tau.
}
Substituting this into \ef{F-Fhat-split-mid} and recalling that $\delta \leq 1/(4 \sqrt{k})$ and $\sigma \geq 1$ we deduce that
\eas{
\nm{F - \hat{F} }_{L^{\infty}_{\mu}(\cX ; \cY)} \lesssim &~ \frac{\sigma C(\bm{b},\varepsilon,\bm{\omega},p)}{\sqrt{\epsilon}} k^{1-1/p}
 + \frac{\sigma \sqrt{k}}{\sqrt{\epsilon}} \nm{f - \hat{\cB}_{\cY} \circ f}_{L^{2}_{\varrho}(\bbR^{\bbN};\cY)} + \frac{\sigma\sqrt{k} \nm{\bm{E}}_{2;\cY}}{\sqrt{m}} 
 \\
 & + \sqrt{k} \tau + \sigma \delta k  \nm{F}_{L^{\infty}_{\nu}(\cX ; \cY)}.
}
After rewriting
\bes{
\nm{f - \hat{\cB}_{\cY} \circ f}_{L^{2}_{\varrho}(\bbR^{\bbN};\cY)} = \nm{F - \hat{\cB}_{\cY} \circ F}_{L^2_{\nu}(\cX ; \cY)}
}
and recalling that $\nm{F}_{L^{\infty}_{\nu}(\cX ; \cY)} \leq 1$ since $F \in \cH(\bm{b},\varepsilon ; \cX , \cY)$, we obtain the desired bound.
}

\subsubsection{Final arguments}\label{sss:final-arguments}

\prf{[Proof of Theorem \ref{thm:operator-learning-extrap}]
Let $c$ be the constant from Theorem \ref{thm:opl-ood-app} and suppose, without loss of generality, that $c \geq 1 /\log^4(2)$. We argue as in the proof of Theorem \ref{thm:poly-optimized-least-squares} and set $k$ as in \ef{k-def} and let $\bar{m} = \bar{m}(\epsilon)$ be the smallest $m \in \bbN$, $m \geq 2$, such that $k \geq 1$. We now set $\delta = \min \{ (4 \sqrt{k})^{-1} , 2^{-m} / k \}$ and apply Theorem \ref{thm:opl-ood-app} with $\sigma = 1$ and $\tau = 0$. This completes the proof.
}

\prf{[Proof of Theorem \ref{thm:deep-learning-extrap}]
We will apply Theorem \ref{thm:operator-learning-extrap}. Let $\cY = \bbR$ with the Euclidean inner product and $d_{\cY} = 1$, so that $\hat{\cD}_{\cY} = \hat{\cB}_{\cY} = \cI_{\cY}$ is the identity map. Let $\cX$ be any separable, infinite-dimensional Hilbert space, $(\lambda_i)_{i \in \bbN}$ be any summable sequence with $\lambda_1 \geq \lambda_2 \geq \cdots > 0$ and $\{ \phi_i \}_{i \in \bbN} \subset \cX$ be any orthonormal basis. Define the probability measure $\nu$ on $\cX$ as the pushforward $\nu = \cD_{\cX} \sharp \varrho$.
Define the encoder $\hat{\cE}_{\cX} : X \mapsto (\ip{X}{\phi_i}_{\cX} / \sqrt{\lambda_i})^{d_{\cX}}_{i=1}$. It is clear that Assumption \ref{ass:opl} holds. Given $\bm{\omega} > \bm{1}$ with $\omega_{\min} > 1$, let $\cN$ be the DNN model class with $n_0 = d_{\cX} = \lceil c m / L \rceil$ and $n_{M+1} = d_{\cY} = 1$ whose existence is established by Theorem \ref{thm:operator-learning-extrap}. Now let $f \in \cH(\bm{b},\varepsilon)$ for some $\bm{b}$ satisfying \ef{b-r-cond-2-main} with $0 < p < 1$ and $\bm{x}_{1},\ldots,\bm{x}_m \sim_{\mathrm{i.i.d.}} \varrho$, where $\varrho$ is the uniform probability measure on $(-1,1)^{\bbN}$. Let $X_i = \cD_{\cX}(\bm{x}_i)$ and notice that $X_{1},\ldots,X_m \sim_{\mathrm{i.i.d.}} \nu$. Further, let $F = f \circ \cE_{\cX}$ and observe that $F \in \cH(\bm{b},\varepsilon ; \cX , \cY)$. Then we can write
\bes{
y_i = f(\bm{x}_i) + e_i = F(X_i) + e_i  = : F(X_i) + E_i = : Y_i
}
and observe that any minimizer $\hat{N}$ of \ef{training-data-intro} is also a minimizer of \ef{DNN-training-prob-operators}. Now let $\mu$ be an arbitrary probability measure on $\bbR^{\bbN}$ that is supported in $D_{\bm{\omega}}$. Define the measure $\upsilon = \cD_{\cX} \sharp \mu$ on $\cX$ and notice that $\upsilon$ satisfies Assumption \ref{ass:opl-ood}.  Observe that $\hat{F} = \hat{\cD}_{\cY} \circ \hat{N} \circ \hat{\cE}_{\cX} = \hat{N} \circ \cE_{\cX}$ and therefore
\bes{
\nm{f - \hat{N}}_{L^{\infty}_{\mu}(\bbR^{\bbN})} = \nm{f \circ \cE_{\cX} - \hat{N} \circ \cE_{\cX}}_{L^{\infty}_{\cD_{\cX} \sharp \mu}(\cX ; \bbR)} = \nm{F - \hat{F} }_{L^{\infty}_{\upsilon}(\cX ; \bbR)} ,
}
since $\hat{\cE}_{\cX} \circ \cD_{\cX}$ is the identity on $\bbR^{\bbN}$. We now apply the error bound from Theorem \ref{thm:operator-learning-extrap} and notice that the second term vanishes since $\hat{\cB}_{\cY} = \cI_{\cY}$.
}

\section{Further details on the numerical experiments}\label{app:experiments}

We now give full details on the numerical experiments. First, we remark that all the experiments in the paper compute a certain error versus number of training samples $m$. We use multiple trials and, following \cite[\S A.1.3]{adcock2022sparse}, we display the geometric mean and one geometric standard deviation. The number of trials is described below.

\subsection{Polynomial estimators}\label{app:ALS}

For the experiments in \S \ref{ss:num-exp-poly}, we compute the relative $L^{\infty}_{\mu}(\bbR^d)$-norm error. Estimating this error accurately in high dimensions may be challenging. We use a grid of size $10^4$, drawn randomly and independently from $\mu$. However, since the functions in this section are sharply peaked near the bottom left corner of any hyperrectangular region, we also augment this grid with the point $\bm{x} = (-\omega_1,\ldots,-\omega_d)$. We use a total of 30 trials in our experiments. 

As mentioned, we employ an adaptive Least Squares (ALS) estimator. This follows \cite{migliorati2015adaptive} (see also \cite{migliorati2019adaptive}, which is based on techniques for the constructing adaptive sparse grid quadratures \cite{gerstner2003dimension}. ALS iteratively computes a sequence of least-squares polynomial approximations $\hat{f}^{(1)},\hat{f}^{(2)},\ldots$ and corresponding (nested) multi-index sets $S^{(1)} \subseteq S^{(2)} \subseteq \cdots \subset \bbN^d_0$, where $\hat{f}^{(l)} \in \cP_{S^{(l)}}$, $\forall l \in \bbN$ and $\cP_{S^{(l)}}$ is as in \ef{f-S-def}. 

We now describe how the one-step update in ALS is performed. This requires several definitions. First, we define the \textit{margin} of a multi-index set $S \subseteq \bbN^d_0$ as
\bes{
\cM(S) = \left \{ \bm{n} \in \bbN^d_0 \backslash S : \exists j \in [d] : \bm{n} - \bm{e}_j \in S \right \},
}
where $\bm{e}_j$ is $j$th canonical multi-index, and the \textit{reduced margin} of $S$ as the set
\bes{
\cR(S) =  \{ \bm{n} = (\nu_j)^{d}_{j=1} \in \cM(S) : \forall j \in [d], \nu_j \neq0\ \Rightarrow\ \bm{n} - \bm{e}_j \in S \}.
}
Given a current multi-index set $S \subset \bbN^d_0$, ALS constructs a new set by selecting multi-indices $\cR(S)$ using a \textit{bulk chasing procedure}, defined as follows. Let $e : \cR(S) \rightarrow \bbR$ and $0 < \beta \leq 1$ be a parameter. The function $e$ serves as an estimate for the true polynomial coefficients of the function being approximated. We describe a precise choice of $e$ below. With this, the procedure $\texttt{bulk}(\cR(S) , e , \beta )$ computes a set $T \subseteq \cR(S)$ of minimal positive cardinality such that
\bes{
\sum_{\bm{n} \in T} e(\bm{n}) \geq \beta \sum_{\bm{n} \in \cR(S)} e(\bm{n}).
}
Now suppose at some iteration we have an multi-index set $S \subset \cF$ and corresponding least-squares approximation $\hat{f} \in \cP_S$ based on sample points $\bm{x}_1,\ldots,\bm{x}_m$. We define $e : \cR(S) \rightarrow \bbR$ as the discrete residual
\be{
\label{bulk-estimator}
e(\bm{n}) = \left | \frac1m \sum^{m}_{i=1} \left ( f(\bm{x}_i) - \hat{f}(\bm{x}_i) \right) \Psi_{\bm{n}}(\bm{x}_i) \right |^2,\quad \forall \bm{n} \in \cR(S).
}
We now compute $T = \mathrm{bulk}(\cR(S),e,\beta)$, construct the new index set $S_{\mathsf{new}} = S \cup T$ and compute a new estimator $\hat{f}_{\mathsf{new}}$ as the least-squares polynomial approximation to $f$ from $\cP_{S_{\mathsf{new}}}$. 

The full ALS approximation is presented in Algorithm \ref{a:ALS}. Note that, following \cite{migliorati2019adaptive}, we use the value $\beta = 0.5$ in our experiments. We also make the choice $m^{(l)} = c_d \max \{ n^{(l)}+1 , \lceil n^{(l)} \cdot \log(n^{(l)}) \rceil \}$, where $c_d = 2$ when $d = 8$ and $c_d = 1$ otherwise. The reason for the larger constant in lower dimensions is to avoid instability of the least-squares system.

\begin{algorithm}[t]
\caption{Adaptive least-squares approximation}
\label{a:ALS}
\begin{algorithmic}
\REQUIRE{Function to approximate $f \in L^2_{\varrho}([-1,1]^d)$, bulk chasing parameter $0 < \beta \leq 1$}
\ENSURE{Sequence of multi-index sets $S^{(1)} \subseteq S^{(2)} \subseteq \cdots \subset \bbN^d_0$ and approximations $\hat{f}^{(1)},\hat{f}^{(2)},\ldots$}
\STATE{Set $S^{(1)} = \{ \bm{0} \}$}
\FOR{$l = 1,2,\ldots$}
\STATE{Set $n^{(l)} = |S^{(l)} |$}
\STATE{Determine a number of samples $m^{(l)} \geq n^{(l)}$}
\STATE{Draw $\bm{x}_1,\ldots,\bm{x}_{m^{(l)}} \sim_{\mathrm{i.i.d.}} \varrho$}
\STATE{Compute $\hat{f}^{(l)} \in \argmin{p \in \cP^{(l)}} \frac{1}{m^{(l)}} \sum^{m^{(l)}}_{i=1} | f(\bm{x}_i) - p(\bm{x}_i) |^2$ }
\STATE{Compute the estimator $e = e^{(l)}$ via \ef{bulk-estimator} with $m = m^{(l)}$, $S = S^{(l)}$ and $\hat{f} = \hat{f}^{(l)}$}
\STATE{Compute $T^{(l)} = \texttt{bulk}(\cR(S^{(l)}) , e^{(l)} , \beta )$}
\STATE{Set $S^{(l+1)} = S^{(l)} \cup T^{(l)}$}
\ENDFOR
\end{algorithmic}
\end{algorithm}

\subsection{DNN estimators}
\label{app:DNN}
For the experiments in \S~\ref{ss:ex-DNNs}, we compute the relative $L^{2}_{\mu}(\bbR^d)$-norm error. We do this on a grid of size $10^4$, drawn randomly and independently from $\mu$. 
We use fully connected DNNs with the architecture and parameters specified in Table~\ref{tab:dnn_parameters}. For each training-set size, we perform 30 independent trials. 

\begin{table}[t]
\centering
\begin{tabular}{ll}
\hline
\textbf{Parameter} & \textbf{Value} \\
\hline
Hidden layers & 15 \\
Hidden layer width & 150 \\
Activation & $\tanh$ \\
Training samples & $\{10,20,40,80,160,320,640,1280\}$ \\
Test samples & 10,000 \\
Batch size & Full batch (all training samples) \\
Epochs & 2000 \\
Optimizer & Adam \\
Weight decay & 0 \\
Initial learning rate & $10^{-4}$ \\
Learning-rate schedule & Multiplied by $0.999$ every epoch \\
\hline
\end{tabular}
\caption{DNN architecture and training parameters for benchmark
functions from the \textit{Virtual Library of Simulation Experiments}.}
\label{tab:dnn_parameters}
\end{table}

\subsection{DNO estimators}
\label{app:DNO} We now describe the DNO setup for the experiments in \S \ref{ss:ex-DNO}. As noted we use FNOs. Our setup and implementation is based on \cite{li2021fourier, li2024physics}.

\paragraph{Darcy flow} 

\begin{table}[t]
\centering
\begin{tabular}{ll}
\hline
\textbf{Parameter} & \textbf{Value} \\
\hline
Training resolution & $64 \times 64$ \\
Testing resolution & $64 \times 64$ \\
Training samples & $\{10,20,40,80,160,320,640,1280\}$ \\
Test samples & 200 \\
Input channels & 1 \\
Output channels & 1 \\
Fourier layers & 4 \\
Fourier modes & 12 \\
Hidden channels & 32 \\
Projection channels & 2 \\
Activation & GeLU \\
Batch size & 32 \\
Epochs & 100 \\
Optimizer & AdamW \\
Weight decay & $10^{-4}$ \\
Initial learning rate & $10^{-3}$ \\
Learning-rate schedule & Halved every 100 epochs \\
\hline
\end{tabular}
\caption{FNO architecture and training parameters for the Darcy Flow problem.}
\label{tab:fno_parameters}
\end{table}
















Details of the architecture and training parameters are given in Table \ref{tab:fno_parameters}. In our experiments, we show the discrete relative $L^2_{\mu}$-norm error, computed over the computational grid. For a fixed input-output pair $(a,u)$, this is defined as
\[
    e_{\mathrm{rel}}
    =
    \frac{
        \left(
        \sum_{p,q}
        |\widetilde{u}_{p,q}-u_{p,q}|^2
        \right)^{1/2}
    }{
        \left(
        \sum_{p,q}|u_{p,q}|^2
        \right)^{1/2}
    },
\]
where $\widetilde{u}$ denotes the FNO prediction, and $u_{p,q}$ and $\widetilde{u}_{p,q}$ are the grid values of $u$ and $\tilde{u}$, respectively.

The test distribution $\mu^{\mathrm{cookies}}$ is based on a benchmark commonly
known as the \textit{fixed-radius cookie problem}
\cite{ballani2015hierarchical,back2011stochastic,chkifa2015discrete}.  Let \(D=[0,1]^2\), and let \(\Omega_i\), \(i=1,\ldots,8\), be
non-overlapping disks of radius \(r_{\mathrm{cookie}}\), with centres
$(0.5 \pm 0.3, 0.5 \pm 0.3), (0.5, 0.5 \pm 0.3), (0.5 \pm 0.3, 0.5)$. Then $a \sim \mu^{\mathrm{cookies}}$ if
\be{
\label{cookies}
a(x)
=
1-\sum_{i=1}^{8}
\mathbf{1}_{\Omega_i}(x)
\left(C_1+C_2y_i\right),\quad \text{where } y_1,\ldots,y_{8} \sim_{\mathrm{i.i.d.}} \cU([-1,1]).
}
In our experiments, we set $r_{\mathrm{cookie}}=0.14$, $C_1= 0.625, C_2=0.375$.

\begin{figure}
\centering    \includegraphics[width=\linewidth]{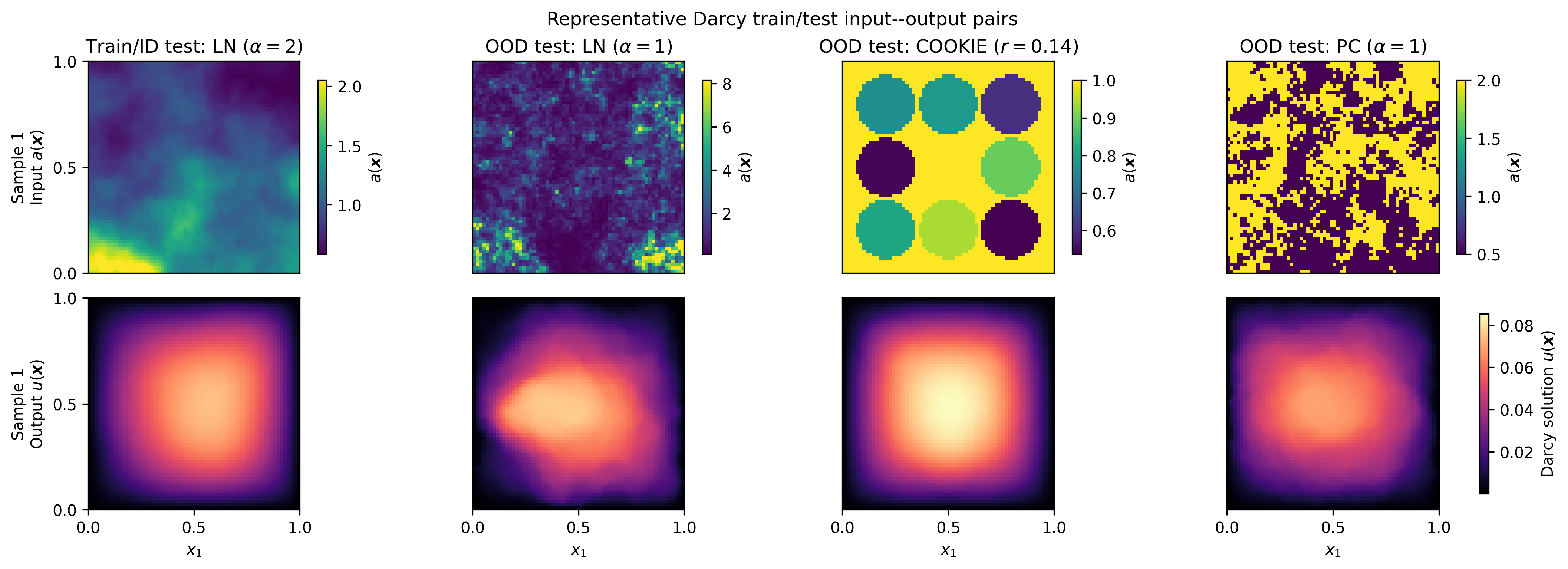}
    \caption{Representative input-output pairs from the Darcy flow
training and test distributions. The top row shows the diffusion coefficients $a$, and the bottom row shows the corresponding
solutions $u$. The distributions are $\mu^{\mathrm{LN}}_{3,2}$, $\mu^{\mathrm{LN}}_{3,1}$, $\mu^{\mathrm{cookies}}$ and $\mu^{\mathrm{PC}}_{3,1}$ (left to right).}
    \label{fig:input-realizations}
\end{figure}

Fig.\ ~\ref{fig:input-realizations} illustrates selected input--output
pairs drawn from the different training and test distributions used in our experiments.

\paragraph{Incompressible Navier--Stokes}

\begin{table}[t]
\centering
\begin{tabular}{ll}
\hline
\textbf{Parameter} & \textbf{Value} \\
\hline
Training resolution & $64 \times 64 \times 50$ \\
Testing resolution & $64 \times 64 \times 50$ \\
Training samples & $[150,300,600,1200,2400,4800]$ \\
Test samples & 200 \\
Input channels & 4 \\
Output channels & 1 \\
Fourier layers & 4 \\
Fourier modes & 12 \\
Hidden channels & 32 \\
Projection channels & 2 \\
Activation & GeLU \\
Batch size & 32 \\
Epochs & 200 \\
Optimizer & Adam \\
Weight decay & 0 \\
Initial learning rate & $10^{-3}$ \\
Learning-rate schedule & Halved every 50 epochs \\
\hline
\end{tabular}
\caption{FNO architecture and training parameters for the Navier-Stokes problem.}
\label{tab:fno_parameters_ns}
\end{table}

Details of the architecutre and training parameters are given in Table \ref{tab:fno_parameters_ns}. In our experiments, we show the discrete relative $L^2_{\mu}$-norm error, computed jointly spatio-temporal computational grid. For a fixed input-output pair $(w_0,w)$, this is defined as
\[
    e_{\mathrm{rel}}
    =
    \frac{
        \left(  
        \sum_{p,q,l}
        \left|
        \widetilde{w}_{p,q,l}
        -
        w_{p,q,l}
        \right|^2
        \right)^{1/2}
    }{
        \left(
        \sum_{p,q,l}
        \left|
        w_{p,q,l}
        \right|^2
        \right)^{1/2}
    },
\]
where $\widetilde{w}$ denotes the FNO prediction, and $w_{p,q,l}$ and $\widetilde{w}_{p,q,l}$ are the grid values of $w$ and $\widetilde{w}$, respectively.

\bibliographystyle{siamplain}
\bibliography{OODrefs}